\pdfoutput=1

\documentclass[11pt,a4paper,twoside]{article}
\usepackage[UKenglish]{babel}
\usepackage{amsmath,amsxtra,amsthm,amssymb,amsfonts,wasysym,
bbm,pifont,color,subfig,wrapfig,epsfig,multirow,enumerate,tensor,bm,
epstopdf,graphicx,cite,mathrsfs,setspace,fancyhdr,mathtools}

\usepackage[all]{xy}
\usepackage[margin=2cm,inner=2cm]{geometry}
\usepackage[bottom]{footmisc}
\DefineFNsymbolsTM{myfnsymbols}{
  \textdagger    \dagger
  \textdaggerdbl \ddagger
  \textsection   \mathsection
  \textbardbl    \|
  \textparagraph \mathparagraph
  \textasteriskcentered *
}
\setfnsymbol{myfnsymbols}
\newtheorem*{theorem'}{Theorem}
\newtheorem*{lemma'}{Lemma}
\newtheorem*{proposition'}{Proposition}
\newtheorem*{conjecture'}{Conjecture}
\newtheorem*{claim'}{Claim}
\newtheorem*{overview'}{Overview}
\newtheorem*{question}{Question}

\newtheorem{theorem}{Theorem}[section]
\newtheorem{lemma}[theorem]{Lemma}

\newtheorem{claim}[theorem]{Claim}

\newtheorem{conjecture}[theorem]{Conjecture}
\theoremstyle{definition}
\newtheorem{definition}[theorem]{Definition}
\newtheorem{remark}[theorem]{Remark}

\newcommand{\nl}{\\$ _{}$\\}
\newcommand{\dd}{\mathrm{d}}
\newcommand{\ee}{\mathrm{e}}
\newcommand{\RR}{\mathbbm{R}}
\newcommand{\pr}{\mathbbm{P}}
\newcommand{\WW}{\mathbbm{W}}
\newcommand{\QQ}{\mathbbm{Q}}
\newcommand{\EE}{\mathcal{E}}
\newcommand{\ZZ}{\mathbbm{Z}}
\newcommand{\ex}{\mathbbm{E}}

\newcommand{\BESQ}[2]{\mathrm{BES}Q^{{#1}}({#2})}
\newcommand{\bp}{\vspace{-0.5\baselineskip}\begin{proof}}
\newcommand{\bpc}{\vspace{-0.5\baselineskip}\begin{proof}[Proof of Claim]}
\newcommand{\ep}{\renewcommand{\qedsymbol}{$\square$}\end{proof}}
\newcommand{\epc}{\renewcommand{\qedsymbol}{$\square$}\end{proof}}

\newcommand{\myparagraph}{\vspace{-\baselineskip}\paragraph}

\numberwithin{equation}{section}
\numberwithin{figure}{section}
\numberwithin{table}{section}
\begin{document}

%\setstretch{1.3}
%\raggedbottom

\title{A universal exponent for Brownian entropic repulsion}
\author{Edward Mottram
\footnote{University of Cambridge, Centre for Mathematical Sciences, Wilberforce Road, Cambridge. CB3 0WA. UK. $_{}$ $_{}$ $_{}$ $_{}$ $_{}$ $_{}$ $_{}$ $_{}$
E-mail: e.j.mottram@maths.cam.ac.uk}}
\date{\today}
\maketitle

\begin{abstract}
\vspace{-2\baselineskip} \[\]
We investigate the extent to which the phenomenon of Brownian entropic repulsion is universal.  Consider a Brownian motion conditioned on the event $\EE$ -- that its local time is bounded everywhere by 1.  This event has probability zero and so must be approximated by events of positive probability.  We prove that several natural quantities, in particular the speed of the process, are highly sensitive to the approximation procedure, and hence are not universal.  However, we also propose an exponent $\kappa$ -- which measures the strength of the entropic repulsion by evaluating the probability that a particular point comes close to violating the condition $\EE$.  We show that $\kappa=3$ for several natural approximations of $\EE$, and conjecture that $\kappa=3$ is universal in a sense that we make precise.
%We investigate the ballistic behaviour of a Brownian motion $(W_t)_{t\geq 0}$ conditioned to have bounded local time.  In particular we condition on the event $\mathcal{E}=\{L_x(t)\leq 1\text{ for all }x\text{ and }t\}$.  Since $\mathcal{E}$ is an event with probability 0 it has to be realised as a limit of events with positive probability.  In this paper it is shown that if we condition on the events $\mathcal{E}^{\bullet}_T=\{L_x(T)\leq 1 \text{ for all } x\in\RR\}$ then $\WW(\cdot\,|\,\mathcal{E}^{\bullet}_T)$ converges weakly to a measure $\mathbbm{Q}^{\bullet}$ as $T\longrightarrow\infty$.  Furthermore $W_t$ has a $\QQ^{\bullet}$-almost sure limiting speed $\gamma^{\bullet}$ as $t\longrightarrow\infty$.  Since $\gamma^{\bullet}$ is strictly less than $\gamma^*$ -- the limiting speed obtained in a paper of Benjamini and Berestycki where a different decomposition of $\EE$ is used, \cite{b&b1} -- we see that the speed of $W_t$ is sensitive to the particular way in which $\EE$ is decomposed. However, in both cases we have $\mathbbm{P}(L_x(\infty)>1-\varepsilon)\asymp \varepsilon^3$, suggesting that there is a sense to which Brownian entropic repulsion is universal.
\end{abstract}

\setlength{\parindent}{0pt}
\setlength{\parskip}{0.5\baselineskip}



\section{Introduction}\label{sec:introduction}

Suppose $(W_t)_{t\geq 0}$ is a Brownian motion conditioned to have bounded local time, $L_x(t)\leq 1$ for all $x\in\RR$ and all $t\geq 0$ say.  Under this conditioning $W_t$ has a self avoiding nature, and so intuitively one would expect $W_t$ to escape to infinity with positive speed.  Moreover, we would expect this speed to be at least equal to 1 since that is precisely what it means to spend less than 1 unit of local time at a given level.  Because it is relatively expensive for a Brownian motion to have positive velocity we might expect the limiting speed to be 1.  However, since the local time of a Brownian motion can fluctuate wildly, the effect of \emph{entropic repulsion} comes into play, and so the easiest way for the process to meet the global constraint $L_x(t)\leq 1$ is for it to have an average local time which is significantly less than 1.  This means that the speed of the process must be strictly greater than 1.

In \cite[Theorem 2]{b&b1} Benjamini and Berestycki make this argument precise in the following way.  Set $\tau_a=\inf\{t : W_t \geq a\}$, and let
\begin{align}
\EE^*_a = \{L_x(t) \leq 1 \text{ for all } x\in\RR \text{ and all } t\leq\tau_a  \} = \{L_x(\tau_a) \leq 1 \text{ for all } x\in\RR \} , \label{eq:Ea*}
\end{align}
then the conditioned Wiener measures $\WW(\cdot\,|\,\EE^*_a)$ converge weakly to a measure $\QQ^*$ as $a\longrightarrow\infty$.  Moreover, $\QQ^*$-almost surely we have
\begin{align}
\lim_{t\rightarrow\infty}\frac{W_t}{t} = \gamma^* = \frac{3}{1-2j_0^{-2}}=4.586\ldots\, ,
\end{align}
where $j_0$ is the first zero of the Bessel function of the first kind, $\mathcal{J}_0(x)$.

It is clear that $\bigcap_{a>0}\EE^*_a = \{L_x(t)\leq 1\text{ for all }x\text{ and }t\}=\EE$.  However, the the event $\EE$ can be realised as a limit of events with positive $\WW$-probability in many other ways, and one could argue that it is more natural to fix $T$ and condition on the events
\begin{align}
\tilde{\EE}_T^{\bullet}=\{L_x(t)\leq 1 \text{ for all }x\in\mathbbm{R} \text{ and all }t\leq T\} = \{L_x(T)\leq 1 \text{ for all }x\in\mathbbm{R}\}. \label{eq:Eabu0}
\end{align}
Notice that again we have $\bigcap_{T>0} \tilde{\EE}_T^{\bullet} = \EE$.

Conditionally on $\EE_a^{*}$ we have $W_{\tau_a}=a>0$, and so $W_t\longrightarrow + \infty$, $\QQ^{*}$-almost surely.  However, when we condition on $\tilde{\EE}^{\bullet}_T$ there is no preferred direction for $W_t$.  For simplicity we now restrict our attention to the case $W_T\geq 0$, and therefore we replace $\tilde{\EE}_T^{\bullet}$ by
\begin{align}
\EE_T^{\bullet}=\{L_x(T)\leq 1 \text{ for all }x\in\mathbbm{R}\text{ and }W_T\geq 0\}. \label{eq:Eabul}
\end{align}
Having done this we can now prove the following.

\begin{theorem}\label{thm:reduced_speed}
The conditioned Wiener measures $\WW(\cdot\,|\,\EE_T^{\bullet})$ converge weakly to a measure $\QQ^{\bullet}$ as $T\longrightarrow\infty$.  Moreover, $\QQ^{\bullet}$-almost surely we have
\begin{align}
\lim_{t\rightarrow\infty}\frac{W_t}{t} = \gamma^{\bullet} ,\label{eq:thm_speed}
\end{align}
where $1<\gamma^{\bullet}<\gamma^{*}$.
\end{theorem}

\begin{remark}
By using symmetry we can reconstruct $\WW(\cdot\,|\,\tilde{\EE}^{\bullet}_T)$ from $\WW(\cdot\,|\,{\EE}^{\bullet}_T)$.  Therefore from Theorem \ref{thm:reduced_speed} we can also deduce that the limit $\tilde{\QQ}^\bullet =\displaystyle{ \lim_{T\rightarrow\infty} \WW(\cdot\,|\,\tilde{\EE}^\bullet_T)}$ exists and that $\displaystyle{\lim_{t\rightarrow\infty}\frac{W_t}{t} \in \{-\gamma^{\bullet},\gamma^{\bullet}\}}$ $\tilde{\QQ}^\bullet$-almost surely.

%Replacing $\tilde{\EE}^{\bullet}_{T}$ by $\EE_{T}^{\bullet}$ results in no loss of generality since we can easily reconstruct $\WW(\cdot\,|\,\tilde{\EE}^{\bullet}_T)$ and $\displaystyle{\lim_{T\rightarrow\infty}\WW(\cdot\,|\,\tilde{\EE}^{\bullet}_T)}$ by using the identity $\WW(A\,|\,\tilde{\EE}^{\bullet}_T)=\frac{1}{2}\WW(A\,|\,\EE^{\bullet}_T) + \frac{1}{2}\WW(-A\,|\,\EE^{\bullet}_T)$ for each measurable $A\in\mathcal{F}$.
\end{remark}

Theorem \ref{thm:reduced_speed} shows us that the limiting speed of the process is sensitive to the particular way that we condition on $\EE$.  Therefore it is clear that for a general set of events $\{\EE^{\prime}_a\}_{a>0}$, with $\EE=\bigcap_{a} \EE^{\prime}_a$, no limiting process need exist.  In Section \ref{sec:discussion} we shall suggest a general framework where -- provided $v$ is not too small -- it is possible to construct a sequence of stopping times $\{\tau_a^v\}_{a > 0}$ with $\tau_a^v\longrightarrow\infty$ $\WW$-almost surely, and a sequence of events
\begin{align}
\EE_a^v=\{L_x(\tau_{a}^v)\leq 1 \text{ for all }x\in\mathbbm{R} \text{ and }W_T\geq 0\},  \label{eq:E^v}
\end{align}
for which the limit $\WW(\cdot\,|\,\EE^{v}_a)\longrightarrow\QQ^{v}$ does exist.  Moreover, in each of theses cases we conjecture that $\displaystyle{\lim_{t\rightarrow\infty}\dfrac{W_t}{t} = v}$ in $\QQ^v$-probability.

A calculation of the speed $\displaystyle{v=\lim_{t\rightarrow\infty}\frac{W_t}{t}}$ might perhaps seem like the most natural way of measuring the entropic repulsion phenomenon.  However, we see from Theorem \ref{thm:reduced_speed} that the value of $v$ is highly sensitive to the approximation of $\EE$ by events of positive probability.
Therefore, in order to find a more universal way of quantifying Brownian entropic repulsion, we shall also consider how likely it is for $L_x(\infty)=\displaystyle{\lim_{T\rightarrow\infty} }L_x(T)$ -- the local time at level $x$ -- to be close to $1$.  We can then prove the following.

\begin{theorem}\label{thm:epsilon3a}
There exists constants $C^{*}$ and $C^{\bullet}$ such that
\begin{align}
\lim_{x\rightarrow\infty}\QQ^{*}(L_x(\infty)>1-\varepsilon)\sim C^{*} \varepsilon^3
\quad\text{and}\quad
\lim_{x\rightarrow\infty}\QQ^{\bullet}(L_x(\infty)>1-\varepsilon)\sim C^{\bullet} \varepsilon^3 
\end{align}
as $\varepsilon\longrightarrow 0$.
\end{theorem}

Similar behaviour is also present in the general framework of Section \ref{sec:discussion}, and therefore we conjecture that this is a universal property of Brownian entropic repulsion.

\subsection{Measuring the entropic repulsion of the Gaussian free field}

In \cite{bdg} and \cite{bdz} the authors consider the Gaussian free field, $\phi$, with a hard wall at 0.  It is shown that if the field is conditioned to be positive on some open set $D\subseteq(\mathbbm{Z}/N\mathbbm{Z})^d$, then the value of the field on $D$ is typically of order $\log N$.  This behaviour occurs because of the effect of entropic repulsion.  In fact the easiest way for the field to satisfy the condition $\phi_x\geq 0$ for all $x\in D$ is for the field to experience a global shift of a size equal to the size of the largest fluctuations.

Since this paper shows that it is possible to measure the strength of Brownian entropic repulsion by looking at how likely it is for $L_x(\infty)$ to be close to 1, we also ask if there is anything similar that can be said for the Gaussian free field.  In particular we pose the following question.

\begin{question}
Can we find scaling functions $f(N)$ and $g(N)$ such that
\begin{align}
\frac{1}{f(N)} \,\pr\left(\phi_y < \varepsilon\, g(N) \,|\, \phi_x> 0\text{ for all }x\in D\right)
\end{align}
converges to a non-trivial function of $\varepsilon$ as $N\longrightarrow\infty$.  If so, then how does this function behave as $\varepsilon\longrightarrow 0$?
%Is there an exponent $\alpha$ such that
%\begin{align}
%\frac{1}{N^\alpha} \,\pr\left(\phi_y < \varepsilon \log N\,|\, \phi_x> 0\text{ for all }x\in D\right)
%\end{align}
%converges to a non-trivial function of $\varepsilon$ as $N\longrightarrow\infty$, and if so then how does this function behave as $\varepsilon\longrightarrow 0$?
\end{question}

\subsection{Outline of the paper}\label{sec:outline}

%Section \ref{sec:intuition} is devoted to preliminaries.  In particular we give an outline of the Ray--Knight Theorems and the Donsker--Varadhan Theorem along with two useful corollaries.  The end of Section \ref{sec:intuition} contains expressions for the Donsker--Varadhan rate functions in integral form.

Section \ref{sec:calc} is devoted to preliminary lemmas.  In particular we find integral forms for the Donsker--Varadhan rate functions of $\mathrm{BES}Q^2$ and $\mathrm{BES}Q^0$ processes, and then use these to find unique minimising measures of the local time process.  From these calculations it is then possible to identify the limiting speed $\gamma^{\bullet}$.

Having done this preliminary work, the proofs of Theorem \ref{thm:reduced_speed} and Theorem \ref{thm:epsilon3a} are then contained in Section \ref{sec:proof1} and Section \ref{sec:proof2} respectively. 
It is worth noting that, once the limiting measure $\QQ$ has been constructed, the proof of Theorem \ref{thm:epsilon3a} is relatively simple.  On the other hand computing of the limiting speed in Theorem \ref{thm:reduced_speed} proves to be a much more arduous task.
%It is interesting to note that whilst $\QQ(L_x(\infty)>1-\varepsilon)$ may seem like a strange quantity to investigate, proving Theorem \ref{thm:epsilon3a} is actually relatively straightforward.  On the other hand, there are many technical estimates required to even show that $\displaystyle{\lim_{t\rightarrow\infty}\frac{W_t}{t}}$ exists, and so the proof of Theorem \ref{thm:reduced_speed} is much more arduous.

The paper closes with a discussion of how our results could be extended to a more general setting, see Section \ref{sec:discussion}.

\section{Preliminaries}\label{sec:intuition}

%This paper relies on the Ray--Knight Theorems and the Donsker--Varadhan Theorem.  This section begins by giving an outline of these results.  We also prove two corollaries which allow us to apply the Ray--Knight and Donsker--Varadhan Theorems in the proof of Theorem \ref{thm:reduced_speed}.

This paper relies on the Ray--Knight Theorems and the Donsker--Varadhan Theorem.  We begin this section by giving an outline of these results.  We also prove a corollary which allows us to apply the Ray--Knight Theorems in the proof of Theorem \ref{thm:reduced_speed}.  
In general the rate function given by the Donsker--Varadhan Theorem is somewhat impenetrable.  However, under certain conditions it can be realised in integral form.  The latter part of this section is devoted to finding integral forms for the rate functions of $\mathrm{BES}Q^2$ and $\mathrm{BES}Q^0$ processes.

\subsection{Local times via the square Bessel process}\label{sec:ray--knight}

Recall that if $(W_t)_{t\geq 0}$ is a $\RR$-valued Brownian motion with Wiener measure $\WW$,  then $\WW$-almost surely the jointly continuous version of the \emph{local time process} is given by
\begin{align}
L_x(t) = \lim_{\varepsilon\rightarrow 0}\frac{1}{2\varepsilon}\int_{0}^{t}\mathbbm{1}_{\{|W_s-x|<\varepsilon\}}\,\dd s,
\end{align} 
for all $x\in\RR$ and all $t\geq 0$.  Roughly speaking the local time tells us how long a Brownian motion has spent at a given point .  We refer the reader to \cite[Chapter VI]{r&y} for an overview of the local time process of Brownian motion.

Recall also that for every $d\geq 0$ there is a unique strong solution of
\begin{align}
\dd Y_x = 2\sqrt{Y_x}\,\dd B_x + d\, \dd x \label{eq:BESQ}
\end{align}
which we call the \emph{square Bessel process of dimension} $d$.  In the case where $Y_0 = c$, we shall write that $(Y_x)_{x\geq 0}$ is a $\mathrm{BES}Q^d(c)$ process.
If we suppose further that $(Y_x)_{0\leq x \leq a}$ is conditioned on the event $Y_a = b$, then $(Y_x)_{0\leq x\leq a}$ becomes a \emph{square Bessel bridge of dimension}\label{page:bridge} $d$ \emph{and length} $a$.  Here we write that $(Y_x)_{0\leq x\leq a}$ is a $\mathrm{BES}Q_a^d(c,b)$ bridge.  Note that although the event $Y_a=b$ has zero probability, we can make sense of such a restriction by conditioning on the event $|Y_a-b|<\varepsilon$ and letting $\varepsilon\longrightarrow 0$.  One can then show that the conditioned measures converge weakly to the law of the square Bessel bridge.  See \cite[Chapter XI, \S 3]{r&y} for details.
%The law of the square Bessel bridge can be realised in other ways, and an explicit law can be given by using \emph{Doob's h-transform}.  See \cite[Chapter XI, \S 3]{r&y} for details.

In this paper we are interested in the cases where $d=2$ or $d=0$ since the local times of a Brownian motion can be related to $\mathrm{BES}Q^2$ and  $\mathrm{BES}Q^0$ processes via the Ray--Knight Theorems.  These apply to a Brownian motion stopped at certain stopping times.  For $a\in\RR$ and $b\geq 0$ we define
\begin{align}
\tau^a = \tau^a_0 = \inf\{t:W_t = a\}\quad \text{and}\quad\tau^a_b=\inf\{t : L_a(t) > b\} , 
\end{align}
and then the Ray--Knight Theorems tell us the following:

\begin{theorem}[First Ray--Knight Theorem]\label{thm:ray_knight1}
Let $(W_t)_{t\geq 0}$ be a Brownian motion, fix $a>0$ and define $Y_x=L_{a-x}(\tau^a_0)$ for $0\leq x\leq a$.  The process $(Y_x)_{0\leq x\leq a}$ is then equal in law to a $\mathrm{BES}Q^2(0)$ process.
\end{theorem}

\begin{theorem}[Second Ray--Knight Theorem]\label{thm:ray_knight2}
Let $(W_t)_{t\geq 0}$ be a Brownian motion, fix $b\geq 0$ and define $Y^+_x=L_x(\tau_b^0)$ for $x\geq 0$ and $Y^-_x=L_{-x}(\tau_b^0)$ for $x\geq 0$.  The processes $(Y^+_x)_{x\geq 0}$ and $(Y^-_x)_{x\geq 0}$ are then equal in law to two independent $\mathrm{BES}Q^0(b)$ processes.
\end{theorem}

For a reference to the Ray--Knight Theorems see \cite[Chapter XI, \S 2]{r&y}.
Theorem 1.2 of \cite[Chapter XI]{r&y} tells us that a $\BESQ{d_1}{b_1}$ process plus a $\BESQ{d_2}{b_2}$ process is equal in law to a $\BESQ{d_1+d_2}{b_1+b_2}$ process.  Therefore it is possible to combine Theorem \ref{thm:ray_knight1} and Theorem \ref{thm:ray_knight2} to describe $(L_x(\tau_b^a))_{x\in\RR}$ for all $a\in\RR$ and each $b\geq 0$.  However, later on we will want to be able to describe $(L_x(T))_{x\in\RR}$ at a fixed (rather than random) $T>0$.  If $W_T\geq 0$ then we can define
\begin{align*}
\displaystyle{S^-_T = \int_{0}^T \mathbbm{1}_{\{W_s < 0\} }\,\dd s}=\int_{-\infty}^0 L_x(T)\,\dd x \quad\text{and}\quad\displaystyle{S^+_T = \int_{0}^T \mathbbm{1}_{\{W_s > W_T\} }\,\dd s}=\int_{W_T}^\infty L_x(T)\,\dd x .
\end{align*}
Using $(Y_x)_{x\geq 0}$ to denote a square Bessel process, we let
\vspace{-0.5\baselineskip}\begin{itemize}
\item $q(a,c,\cdot)$ be the density of $Y_a$ with respect to the law of a $\BESQ{2}{c}$ process.
\item  $f(c,\cdot)$ be the density of ${\int^\infty_0 Y_x \,\dd x}$ with respect to the law of a $\BESQ{0}{c}$ process.
\item $g(a,c,b,\cdot)$ be the density of ${\int_{0}^{a} Y_x\,\dd x}$ with respect to the law of a $\mathrm{BES}Q^2_a(c,b)$ bridge.\label{page:bridge_density}
\end{itemize}
In \cite{Leuridan} we are told that the density $f$ has a relatively simple expression
\begin{align}
f(c,s) = \frac{c}{\sqrt{8\pi }s^{\frac{3}{2}}}\exp\left(-\frac{c^2}{8s}\right),\label{eq:f_density}
\end{align}
whereas the expressions for $q$ and $g$ turn out to be more complicated.  From a result of Leuridan, \cite[Theorem 1]{Leuridan}, we now have the following:

%Here $q_a(c,\cdot)$ denotes the density of $Y_a$ with respect to the law of a $\BESQ{2}{c}$ process, $f(c,\cdot)$ is the density of the random variable $\int^\infty_0 Y_x \,\dd x$ with respect to the law of a $\BESQ{0}{c}$ process, and $g_a(c,b,\cdot)$ is the density of $\int_{0}^{a} Y_x\,\dd x$ with respect to the law of a $\mathrm{BES}Q^2_a(c,b)$ bridge.

\begin{theorem}[Leuridan, 1998]\label{thm:ray_knight_fixed_t}
Let $T>0$ be fixed and suppose $(W_t)_{t\geq 0}$ is a Brownian motion conditioned on the event $\{W_T\geq 0\}$.  The joint distribution of $W_T$ and $(L_x(T))_{x\in\RR}$ is characterised by the following properties.
\vspace{-0.5\baselineskip}\begin{itemize}
\item The 5-tuple $(W_T,L_{W_T}(T),L_0(T),S^-_T,S^+_T)$ admits a probability density on $[0,\infty)^5$,
\begin{align*}
(a,b,c,s^-,s^+) \longmapsto 2\, q(a,c,b) f(c,s^-) f(b,s^+) g(a,c,b,T-s^--s^+).
\end{align*}
\item Conditionally on $(W_T,L_{W_T}(T),L_0(T),S^-_T,S^+_T)=(a,c,b,s^-,s^+)$
\vspace{-0.5\baselineskip}\begin{itemize}
\item $(L_{-x}(T))_{x\geq 0}$, $(L_{a+x}(T))_{x\geq 0}$ and $(L_{x}(T))_{0\leq x\leq a}$ are independent.
\item $(L_{-x}(T))_{x\geq 0}$ and $(L_{a+x}(T))_{x\geq 0}$ are equal in law to $\BESQ{0}{c}$ and $\BESQ{0}{b}$ processes conditioned on the events ${\int_{0}^{\infty} Y_x\,\dd x = s^-}$ and ${\int_{0}^{\infty} Y_x\,\dd x = s^+}$.
\item $(L_{x}(T))_{0\leq x\leq a}$ is equal in law to a $\mathrm{BES}Q^2_a(c,b)$ bridge which has been conditioned on the event that ${\int_{0}^{a} Y_x\,\dd x = T-s^--s^+}$.
\end{itemize}
\end{itemize}
\end{theorem}
Of course the events  $\left\{\int_{0}^{\infty} Y^-_x\,\dd x = s^-\right\}$, $\left\{\int_{0}^{\infty} Y^+_x\,\dd x = s^+\right\}$ and $\left\{\int_{0}^{a} Y_x\,\dd x = T-s^--s^+\right\}$ all have probability 0.  However, as with the construction of the Brownian bridge, the conditioned processes can all be realised as a weak limit.  In \cite[Theorem 1]{Leuridan} we are given explicit generators for these conditioned processes.
\vspace{-0.5\baselineskip}\begin{itemize}
\item The joint Markov processes $\left(L_{-x}(T),\int_{-\infty}^{-x} L_{y}(T)\,\dd y \right)_{x\geq 0}$ and $\left(L_{a+x}(T),\int_{a+x}^{\infty} L_{y}(T)\,\dd y \right)_{x\geq 0}$ both have infinitesimal generator
\begin{align}
2z_1\frac{\partial^2}{\partial z_1^2} + \left(4- \frac{z_1^2}{z_2}\right)\frac{\partial}{\partial z_1} - z_1 \frac{\partial}{\partial z_2 } . \label{eq:BESQ0SDE}
\end{align}
\item The joint process $\left(L_{x}(T),\int_x^a L_y(T)\,\dd y,x\right)_{0\leq x\leq a}$ is Markovian with infinitesimal generator
\begin{align}
2z_1\frac{\partial^2}{\partial z_1^2} + \left(2+4z_1\left(\frac{\partial_1 q_{z_3}}{q_{z_3}}(z_1,b)+\frac{\partial_1 g_{z_3} }{g_{z_3}}(z_1,b,z_2)\right)\right) \frac{\partial}{\partial z_1} - z_1\frac{\partial}{\partial z_2} +\frac{\partial}{\partial z_3} . \label{eq:BESQ2SDE}
\end{align}
\end{itemize}
In \cite{Leuridan} the generators (\ref{eq:BESQ0SDE}) and (\ref{eq:BESQ2SDE}) are both be deduced from \emph{Doob's h-transform}.  See \cite{r&w2} for a detailed explanation of the h-transform. 
In 1999 Pitman also proved a similar result to Theorem \ref{thm:ray_knight_fixed_t} via a branching process approximation.  See \cite{pitman1999}.

\subsection{The invariant density of a diffusion process}

Suppose $M$ is a compact metric space and $(X_t)_{t\geq 0}$ is a diffusion process on $M$ with infinitesimal generator $\mathcal{L}$, and starting point $X_0 = y$.  For each $A\in\mathcal{B}(M)$ let
\begin{align}
L((X_t),T,A)=\frac{1}{T}\int_0^T \mathbbm{1}_{\{X_t\in A\}}\,\dd t .\label{eq:occupation}
\end{align}
Then for a given $(X_t)_{t\geq 0}$ and $T$, $L((X_t),T,A)$ gives the proportion of time that $(X_t)_{0\leq t\leq T}$ spends in $A$, and $L((X_t),T,\cdot)$ defines a probability measure on $M$.  We call $L((X_t),T,\cdot)$ the \emph{occupation measure of $(X_t)$ at time $T$}. Now for each probability measure $\mu\subseteq\mathcal{P}(\RR)$ define
\begin{align}
I(\mu) = -\inf_{u\in\mathcal{D},u>0}\int_{M}\frac{\mathcal{L}(u)}{u}\,\dd \mu , \label{eq:Ifunction1}
\end{align}
where $\mathcal{D}$ is the domain of $\mathcal{L}$.  A theorem of Donsker and Varadhan, \cite[Theorem 1]{d&v}, now gives us the following.

\begin{theorem}[Donsker--Varadhan]\label{thm:donsker_varadhan}
Let $\mathcal{P}(\RR)$ be the space of all probability measures on $\RR$ equipped with the weak topology, and suppose $\pr_y$ is the measure for a process generated by $\mathcal{L}$ and started at $y$.  For all $C\subseteq\mathcal{P}(\RR)$ closed and all $O\subseteq\mathcal{P}(\RR)$ open, if we also assume that $y\in\mathrm{support}(\mu)$ for each $\mu\in O$, then we have
\begin{align}
\limsup_{T\rightarrow\infty}\frac{1}{T}\log\pr_y(L((Y_x),T,\cdot)\in C)&\leq -\inf_{\mu\in C} I(\mu) \label{eq:d&vclosed}\\
\liminf_{T\rightarrow\infty}\frac{1}{T}\log\pr_y(L((Y_x),T,\cdot)\in O)&\geq -\inf_{\mu\in O} I(\mu) \label{eq:d&vopen} .
\end{align}
In other words the occupation measure of the process generated by $\mathcal{L}$ satisfies a \emph{large deviations principle} with rate function $I$.
\end{theorem}

\begin{remark} \label{rem:donsker-varadhan}
By following the proof of the Donsker--Varadhan Theorem in \cite[Section 2]{d&v} we see that (\ref{eq:d&vclosed}) and (\ref{eq:d&vopen}) also hold in the case where $y$ is a random variable on $\RR$.  In this situation we substitute $\pr_y(L((Y_x),T,\cdot)\in \cdot)$ with $\ex^y\big\{\pr_y(L((Y_x),T,\cdot)\in \cdot)\big\}$, and replace the condition that $\big\{y\in\mathrm{support}(\mu)$ for each $\mu \in O\big\}$ with the condition that $\big\{\mathrm{support}(y)\subseteq\mathrm{support}(\mu)$ for each $\mu \in O\big\}$.
\end{remark}

In general the formula for $I(\mu)$ given by (\ref{eq:Ifunction1}) is rather impenetrable.  However, in \cite[Theorem 5]{d&v} it is shown that when $M$ is the real line equipped with the usual metric, and $\mathcal{L}$ is self-adjoint with respect to Lebesgue measure, then
\begin{align}
I(\mu) = \left\lbrace\begin{array}{cl}
\displaystyle{\left\Vert \sqrt{-\mathcal{L}}g\right\Vert_2^2 },& \text{where }g=\sqrt{\dfrac{\mathrm{d}\mu}{\mathrm{d}x}}\text{ exists and is in the domain of }\sqrt{-\mathcal{L}}
\\ \infty , & \text{otherwise}
\end{array}\right. .\label{eq:Ifunction2}
\end{align}
In Chapter 4 we will only be interested in the case when $\mathcal{L}$ is a second order differential operator of the form
\begin{align}
\mathcal{L}f(x) = \frac{1}{2} \frac{\dd}{\dd x} \left( a(x)\frac{\dd}{\dd x}\right)f(x) + b(x) \frac{\dd}{\dd x}f(x) ,
\end{align}
where $a$ is continuous and $b$ is continuously differentiable.  In this setting we can use (\ref{eq:Ifunction2}) to express (\ref{eq:Ifunction1}) as an integral.

First suppose that the drift function $b$ can be written as $b(x) = a(x) \dfrac{\dd}{\dd x}Q(x)$, for a continuously differentiable function $Q(x)$.  In this case $\mathcal{L}$ can be realised as a self-adjoint operator with respect to a measure $\mu_{\mathrm{rev}}$, defined by $\dd\mu_{\mathrm{rev}} = \ee^{2Q}\dd x$.  Therefore, provided $\dfrac{\dd \mu}{\dd \mu_{\mathrm{rev}}}$ exists and $ g=\sqrt{\dfrac{\dd \mu}{\dd \mu_{\mathrm{rev}}}}$ belongs to the domain of $\sqrt{-\mathcal{L}}$, we get
\begin{align}
I(\mu) = \left\Vert \sqrt{-\mathcal{L}} g \right\Vert^2_{2,\mu_{\mathrm{rev}}} = \frac{1}{2} \int a(x) \left(\frac{\dd}{\dd x} g(x) \right)^2 \ee^{2 Q(x)}\,\dd x . \label{eq:Imu}
\end{align}
See \cite{pinsky3} for more details.  It is clear from (\ref{eq:Ifunction2}) that if $\mathcal{L}$ is self-adjoint then $I:\mathcal{P}(\RR)\longrightarrow [0,\infty]$ is lower semi-continuous.
%In fact it follows from a result of Pinsky, \cite{pinsky2}, that $I$ is lower semi-continuous for all infinitesimal generators $\mathcal{L}$.
Note however that $I$ is not continuous since any measure $\mu$ with $g$ in the domain of $\sqrt{-\mathcal{L}}$ can be approximated arbitrary closely by measures without a Radon-Nikodym derivative.  We also remark that $I$ is in fact lower semi-continuous for all infinitesimal generators $\mathcal{L}$.  This follows from \cite{pinsky2} where Pinsky shows that (\ref{eq:Ifunction1}) can be replaced by
\begin{align}
I(\mu) = -\inf_{u\in\mathcal{D},u>0}\int_{M}\frac{\mathcal{L}(u)}{u}\,\dd \mu = -\inf_{u\in C^2(\RR),u>0}\int_{M}\frac{\mathcal{L}(u)}{u}\,\dd \mu . \label{eq:ratefunction}
\end{align}

\section{Calculating the rate function for the occupation measure of a square Bessel process}\label{sec:calc}

Recall that Theorem \ref{thm:ray_knight_fixed_t} tells us that the local time of a Brownian motion can be described in terms of conditioned $\mathrm{BES}Q^2$ and $\mathrm{BES}Q^0$ processes. If $(Y_x)_{x\geq 0}$ is such a process then (\ref{eq:occupation}) tells us that the occupation measure of $(Y_x)_{x\geq 0}$ at time $T$ is given by
\begin{align*}
L((Y_x),T,\cdot) =\frac{1}{T}\int_0^T \mathbbm{1}_{\{Y_x\in\cdot\}}\,\dd x .
\end{align*}
The Donsker--Varadhan Theorem, Theorem \ref{thm:donsker_varadhan}, gives us a powerful tool for estimating the large deviations of the occupation measure in terms of a rate function $I$, see (\ref{eq:Ifunction1}).  From now on we shall use $I_2$ and $I_0$ to denote the Donsker--Varadhan rate functions for the occupation measures of $\mathrm{BES}Q^2$ and $\mathrm{BES}Q^0$ processes.
%Recall that the Ray--Knight Theorems, Theorem \ref{thm:ray_knight1} and Theorem \ref{thm:ray_knight2}, tell us that the local time of Brownian motion can be described in terms of $\mathrm{BES}Q^2$ and $\mathrm{BES}Q^0$ processes.  It will now be useful for us to make various computations to give bounds on the Donsker--Varadhan rate functions in each case.
We shall also use the abuse of notation $\ex(\mu)=\int x\,\dd\mu(x)$ to denote the expectation of the identity with respect to the measure $\mu$.
This section is then devoted to the proofs of the following three lemmas.

\begin{lemma}\label{lem:lemma9}
There is a unique probability measure $\mu^{*}$ supported on $[0,1]$ which minimises $I_2(\mu)$ over all measures $\mu$ with $\mathrm{support}(\mu)\subseteq[0,1]$.  Furthermore,  $\mu^{*}$ is such that $\ex(\mu^{*})=(\gamma^{*})^{-1}=\frac{1}{3}(1-2j_0^{-2})<1$ and $I_2(\mu^{*})=\frac{1}{2}j_0^2>0$.  Here $j_0$ denotes the first zero of the Bessel function $\mathcal{J}_0$.
\end{lemma}

\begin{lemma}\label{lem:I_0}
Suppose $\mu$ is a probability measure which is supported on $(0,1]$, then
\begin{align}
I_0(\mu)\geq 2\pi^2 \,\ex(\mu) .
\end{align}
Furthermore the measure $\mu^{\circ}$ defined by $\dfrac{\dd\mu^{\circ}}{\dd x} =\dfrac{1}{Z} \dfrac{1}{x}\sin(\pi x)^2$, where $\displaystyle{Z=\int_{0}^{1}\frac{1}{x}\sin(\pi x)^2\,\dd x}$, is the unique probability measure for which equality is achieved.
\end{lemma}

\begin{lemma}\label{lem:mu_bullet}
There is a unique probability measure $\mu^{\bullet}$ with $\mathrm{support}(\mu^\bullet)\subseteq[0,1]$ which minimises $\ex(\mu)^{-1}I_2(\mu)$ over all measures $\mu$ with $\mathrm{support}(\mu)\subseteq[0,1]$.  Furthermore, $\mu^{\bullet}$ satisfies $\ex(\mu^\bullet)^{-1}I_2(\mu^\bullet)<2\pi^2$ and $\ex(\mu^{*})<\ex(\mu^{\bullet})<1$.  Here $\mu^*$ is as defined by Lemma \ref{lem:lemma9}.
\end{lemma}

At this point we set $\gamma^\bullet = \ex(\mu^\bullet)^{-1}$, $\gamma^* = \ex(\mu^*)^{-1}$ and $\Gamma^\bullet = \ex(\mu^\bullet)^{-1}I_2(\mu^\bullet)$.  We also make the following definition. 

\begin{definition}\label{def:J}
Define $J:[0,1]\longrightarrow[0,\infty)\cup\{\infty\}$ by
\begin{align}
J(\alpha) = \inf \left\{I_2(\mu):\mathrm{support}(\mu)\subseteq [0,1]\text{ and }\ex(\mu) =\alpha\right\}. \label{eq:J}
\end{align}
%where $I_2$ is the rate function of a $\mathrm{BES}Q^2$ process.  See (\ref{eq:Ifunction3}).
\end{definition}

The proofs of the three lemmas are mainly calculation, and therefore can be omitted in a first reading of the chapter.

\subsection{The rate function of $\mathrm{BES}Q^2$ and $\mathrm{BES}Q^0$ processes}

In order to prove Lemma \ref{lem:lemma9}, Lemma \ref{lem:I_0} and Lemma \ref{lem:mu_bullet} we shall first use (\ref{eq:Imu_int}) and (\ref{eq:Imu}) to find integral forms of $I_2$ and $I_0$.
%the Donsker--Varadhan rate functions of $\mathrm{BES}Q^2$ and $\mathrm{BES}Q^0$ processes.  
Since the square Bessel process of dimension $d$ satisfies (\ref{eq:BESQ}) then its infinitesimal generator is given by
\begin{align}
\mathcal{L}_d f(x) = \frac{1}{2}\frac{\mathrm{d}}{\mathrm{d}x}\left(4x\,\frac{\mathrm{d}}{\mathrm{d}x}\right)f(x) +(d-2) \frac{\dd}{\dd x}f(x) .\label{eq:BESQgen}
\end{align}
When $d=2$ then $\mathcal{L}_d$ is self-adjoint with respect to Lebesgue measure.  Therefore we can use (\ref{eq:Ifunction2}) to write down the Donsker--Varadhan rate function $I_2$ as
\begin{align}
I_2(\mu) = \left\lbrace\begin{array}{cl}
\displaystyle{\left\Vert \sqrt{-\mathcal{L}}g \right\Vert_2^2 }, & \text{where }\dfrac{\mathrm{d}\mu}{\mathrm{d}x}\text{ exists and }g=\sqrt{\dfrac{\mathrm{d}\mu}{\mathrm{d}x}}\in \mathcal{D}_2
\\ \infty , & \text{otherwise}
\end{array}\right. .\label{eq:Ifunction3}
\end{align}
Here $\mathcal{D}_2$ donates the domain of $\mathcal{L}_2$.   Since $g\in\mathcal{D}_2$ implies that $g(x)\longrightarrow 0$ as $x\longrightarrow\infty$, we can integrate by parts to get
\begin{align}
\left\Vert \sqrt{-\mathcal{L}_2}g \right\Vert_2^2 = \left\langle g, -\mathcal{L}_2 g \right\rangle_2 &= -\frac{1}{2} \int_0^{\infty} g(x) \frac{\dd}{\dd x} \left(4x \frac{\dd}{\dd x}\right) g(x) \, \dd x  = \int_0^{\infty} 2x \left( \frac{\dd}{\dd x}g(x)\right)^2\,\dd x . \label{eq:Imu_int}
\end{align}

However, when $d=0$ then $\mathcal{L}_0$ is no longer self-adjoint and so we have to perform a change of measure and use equation (\ref{eq:Imu}) in order to calculate the rate function.  We have the drift term $b(x)=-2$ and so solving $b(x)=a(x)\dfrac{\dd}{\dd x}Q(x)$ gives
\begin{align*}
Q(x)=\int^x -\frac{2}{4y}\,\dd y = -\frac{1}{2}\log x+ c.
\end{align*}
Now set $\dd\mu_{\mathrm{rev}}=\ee^{2Q(x)}\dd x=\ee^{2c}\dfrac{1}{x}\dd x$ for $x>0$.  Note that it does not matter that we have an unknown constant of integration since this cancels later.  Suppose $\mu$ is a probability measure supported on $(0,\infty)$ for which $g=\sqrt{\dfrac{\dd\mu}{\dd x}}$ exists, we then have 
\begin{align*}
h(x)=\sqrt{\dfrac{\dd\mu}{\dd\mu_{\mathrm{rev}}}}(x)=\sqrt{\dfrac{\dd\mu}{\dd x}\dfrac{\dd x}{\dd\mu_{\mathrm{rev}}}}(x)=g(x) \sqrt{\dfrac{\dd x}{\dd\mu_{\mathrm{rev}}}}(x) = g(x)\ee^{-c}\sqrt{x}.
\end{align*}
Therefore from (\ref{eq:Imu}) we get
\begin{align*}
\left\Vert \sqrt{-\mathcal{L}_0} h \right\Vert_{2,\mu_{\mathrm{rev}}}^2 &= \ee^{2c} \frac{1}{2}\int_{0}^{\infty} 4x\left( \frac{\dd}{\dd x} h(x) \right)^2 \frac{1}{x}\,\dd x = \int_{0}^{\infty}2\left( \frac{\dd}{\dd x} \sqrt{x} g(x)\right)^2\,\dd x .
\end{align*}
Hence if we write $\mathcal{D}_0$ for the domain of $\mathcal{L}_0$ then for each $\mu$ supported on $(0,\infty)$ we have
\begin{align}
I_0(\mu) = \left\{ \begin{array}{ll} \displaystyle{\int_{0}^{\infty}2\left( \frac{\dd}{\dd x} \sqrt{x} g(x)\right)^2\,\dd x}, & \text{if } \dfrac{\dd \mu}{\dd x} \text{ exists and }g=\sqrt{\dfrac{\dd \mu}{\dd x}}\in\mathcal{D}_0 \\ \infty , & \text{otherwise} \end{array} \right. . \label{eq:I_0}
\end{align}

\begin{remark}
$Y_x=0$ is an absorbing state for a $\mathrm{BES}Q^0$ process and therefore we shall only be interested in occupation measures supported on $(0,\infty)$.
\end{remark}

\subsection{Proof of Lemma \ref{lem:lemma9}}

Since (\ref{eq:Imu_int}) gives us an explicit form for $I_2(\mu)$, then it is clear that minimising $I_2(\mu)$ over $\{\mu : \mathrm{support}(\mu)\subseteq[0,1]\}$ is equivalent to finding a $g\in C^1([0,1])$ with $\Vert g \Vert_2 = 1$ and
\begin{align}
\int_{0}^{1} 2x \left(\frac{\mathrm{d}}{\mathrm{d}x}g(x)\right)^2\dd x = \inf\left\{ \int_{0}^{1} 2x \left(\frac{\mathrm{d}}{\mathrm{d}x}h(x)\right)^2\dd x : h\in C^1([0,1]) \text{ and }\Vert h \Vert_2=1 \right\}.
\end{align}
Such a problem can be solved using the Euler--Lagrange equation.  The interested reader is directed to \cite[Chapter 9]{boas} for an overview of the calculus of variations.  In this case the additional constraint that $\Vert g \Vert_2 = 1$ can be included by adding the Lagrangian multiplier
\begin{align}
\lambda\left(\int_{0}^{1}g(x)^2\,\dd x -1\right) .
\end{align}
Thus we have
\begin{align}
\cfrac{\partial \mathcal{F}}{\partial g} - \cfrac{\mathrm{d}}{\mathrm{d} x}\left(\cfrac{\partial \mathcal{F}}{\partial g'}\right)= 0 ,
\end{align}
where $\mathcal{F}[x,g(x),g'(x),\lambda]=2x \,g'(x)^2+\lambda(g(x)^2-1)$, and so we get
\begin{align}
2x \frac{\mathrm{d}^2}{\mathrm{d}x^2}g(x)+2\frac{\mathrm{d}}{\mathrm{d}x}g(x)-\lambda g(x) = 0 \label{eq:de1}.
\end{align}
The solutions of this are precisely the eigenfunctions of $\mathcal{L}_2$ (with eigenvalue $\lambda$).  We also have the additional constraints that $g(x)\geq 0$ for $x\in [0,1]$, because $g(x)$ is the positive square root of the Radon--Nikodym derivative of $\mu$, and that $g(1)=0$, because $g(x)$ must be continuous on $[0,\infty)$.  Therefore there is a unique solution
\begin{align}
g(x) = \left\{ \begin{array}{ll}
\dfrac{\mathcal{J}_0(j_0 \sqrt{x})}{\mathcal{J}_1(j_0)}  & \quad \text{for } 0\leq x\leq 1 \\
0 & \quad \text{otherwise}
\end{array} \right. . \label{eq:minmu}
\end{align}
Here $\mathcal{J}_n(x)$ are the Bessel functions of the first kind, and $j_0$ is the first zero of $\mathcal{J}_0(x)$.  Thus by defining $\mu^*$ by $\dfrac{\dd \mu^*}{\dd x} = g(x)^2$ we get the unique minimiser of $I_2$ over $\{\mu:\mathrm{support}(\mu)\subseteq [0,1]\}$.

By calculation we can now check that
\begin{align*}
\ex(\mu^{*})=\int_{0}^{1} x\,\dd\mu(x) = \int_{0}^{1} x\left(\frac{\mathcal{J}_0(j_0 \sqrt{x})}{\mathcal{J}_1(j_0)}\right)^2\,\dd x =\frac{1}{3}(1-2 j_0^{-2})=\frac{1}{\gamma^*} ,
\end{align*}
and
\begin{align*}
\quad I_2(\mu^{*}) =\int_{0}^{1} 2x \left(\frac{\dd}{\dd x}\frac{\mathcal{J}_0(j_0 \sqrt{x})}{\mathcal{J}_1(j_0)}\right)^2\,\dd x = \frac{1}{2} j_0^2 ,
\end{align*}
as claimed.
\hfill$\square$

\begin{remark}
By using Lemma \ref{lem:lemma9} and the Donsker--Varadhan Theorem it is now possible to show the following.

Let $(Y_x)_{x\geq 0}$ be a $\BESQ{2}{y}$ process for some $0\leq y < 1$, and let $\mathbbm{Y}_y^T$ be the law of $(Y_x)_{x\geq 0}$ conditioned on the event $\{Y_x\leq 1\text{ for all }x\leq T\}$.  If we use $\ex_{\mathbbm{Y}_y^T}(X)$ to denote the expectation under this measure of a random variable $X\geq 0$.  Then
\begin{align}
\lim_{t\rightarrow \infty} \lim_{T\rightarrow\infty} \ex_{\mathbbm{Y}_y^T}(Y_t)=\frac{1}{\gamma^*}.
\end{align}

Therefore this method provides an alternative way of proving \cite[Lemma 9]{b&b1} -- a central component in the proof of ballistic behaviour in Benjamini and Berestycki's paper.
\end{remark}

\subsection{Proof of Lemma \ref{lem:I_0}}

We may assume that $g(x)= \sqrt{\dfrac{\dd \mu}{\dd x}}(x)$ exists, is piecewise differentiable on $(0,1]$ and that $\displaystyle{\lim_{x\rightarrow 0} g(x) =0}$ and $g(1)=0$.  If not then $I_0(\mu)=\infty$.  From (\ref{eq:I_0}) we then have
\begin{align*}
I_0(\mu) = \int_0^1 2\left(\frac{\dd}{\dd x} \sqrt{x} g(x)\right)^2\,\dd x \quad\text{and}\quad \ex(\mu) = \int_{0}^{1} x g(x)^2\,\dd x .
\end{align*}
Letting $h(x)=\sqrt{x}g(x)$  this becomes
\begin{align*}
I_0(\mu) = 2 \int_0^1 h'(x)^2\,\dd x \quad\text{and}\quad \ex(\mu) = \int_{0}^{1} h(x)^2\,\dd x .
\end{align*}
We now claim that if $f:[0,1]\longrightarrow\RR$ is a $C^1$ function such that $f(0)=f(1)=0$ then
\begin{align}
\pi^2 \int_0^1 |f(x)|^2\,\dd x\leq \int_0^1 |f'(x)|^2\,\dd x, \label{eq:Wirtinger}
\end{align}
with equality if and only if $f(x)=c\sin\left(\pi x\right)$ for some $c\neq 0$.  To prove this we note that because $f$ is periodic and $C^1$ then it can be written as a Fourier series,
\begin{align*}
f(x) = \sum_{n=1}^\infty a_n \sin(n\pi x) \quad \text{and} \quad f'(x) = \sum_{n=1}^\infty n\pi a_n \sin(n\pi x).
\end{align*}
Parseval's identity now tells us that
\begin{align*}
\int_0^1 |f(x)|^2\,\dd x = \frac{1}{2} \sum_{n=1}^{\infty} a_n^2 \quad \text{and}\quad \int_0^1 |f'(x)|^2\,\dd x = \frac{1}{2} \pi^2 \sum_{n=1}^{\infty} n^2 a_n^2 ,
\end{align*}
from which (\ref{eq:Wirtinger}) is clear.  Furthermore we can only have equality when $a_1\neq 0$ and $a_n= 0$ for all $n>1$, implying that $f(x)=c\sin\left(\pi x\right)$ for some $c\neq 0$.

The lemma now follows since $g(0)$ is finite implying that $h(0)=0$.  Therefore
\begin{align*}
I_0(\mu) = 2 \int_0^1 h'(x)^2\,\dd x \geq 2\pi^2\int_{0}^{1} h(x)^2\,\dd x = 2\pi^2 \ex(\mu),
\end{align*}
with equality if and only if $\dfrac{\dd\mu}{\dd x}=c^2\dfrac{1}{x} \sin(\pi x)^2$.  Here the requirement that $\displaystyle{\int_0^1\dfrac{\dd\mu}{\dd x}\,\dd x =1}$ tells us that $\displaystyle{c^2=\frac{1}{Z}=\left(\int_0^1 \frac{1}{x} \sin (\pi x)^2\,\dd x \right)^{-1}}$.
\hfill$\square$

\begin{remark}
The inequality given by (\ref{eq:Wirtinger}) is a special case of the Poincar\'{e} inequality, and is sometimes known as Wirtinger's inequality.
\end{remark}

\subsection{Proof of Lemma \ref{lem:mu_bullet}}\label{sec:J}

%In Section \ref{sec:closetoP} we encountered the loosely defined function $\tilde{J}$, see (\ref{eq:J_tilde}), which represented the cost of a Brownian motion travelling a unit time at a given speed whilst and having its local time bounded by $1$. 
%Having now formally defined the Donsker--Varadhan rate function for a $\mathrm{BES}Q^2$ process, $I_2$, we can now formally define the related function $J$.

Lemma \ref{lem:mu_bullet} is proved by studying the properties of the function $J$.  See Definition \ref{def:J}.

\begin{lemma}\label{lem:J_props}
$J$ has the following properties.
\begin{enumerate}\vspace{-0.5\baselineskip}
\item $J(\alpha)<\infty$ if and only if $\alpha\in(0,1)$.
\item For each $\alpha\in[0,1]$ there is a unique $\mu_{\alpha}$ with $\mathrm{support}(\mu_{\alpha})\subseteq [0,1]$, $\ex(\mu_{\alpha}) =\alpha$ and $I_2(\mu_{\alpha})=J(\alpha)$.
\item $J$ is strictly convex on $[0,1]$.
\item $J$ is continuously differentiable on $(0,1)$.
\item $J$ has a unique minimum at $\alpha=(\gamma^{*})^{-1}$, with $J\left( (\gamma^{*})^{-1}\right)>0$.
\item $J(\alpha)\longrightarrow\infty$ as $\alpha\longrightarrow 1$.
%\item $J(\alpha)\geq \frac{1}{2}\alpha^{-1}$ for all $\alpha>0$, and thus $J(\alpha)\longrightarrow\infty$ as $\alpha\longrightarrow 0$.
\item The function $v J\left(v^{-1}\right)$ has a unique minimum at $v=\gamma^\bullet$, with $1<\gamma^{\bullet}<\gamma^*$.  What is more this minimum satisfies $\gamma^{\bullet}J\left((\gamma^{\bullet})^{-1}\right)=\Gamma^\bullet<2\pi^2$.
\vspace{-0.5\baselineskip}\end{enumerate}
\bp
We shall prove these claims sequentially.
\begin{enumerate}\vspace{-0.5\baselineskip}
\item The only probability measures supported on $[0,1]$ with expectation 0 or 1 are the point masses $\delta_0$ and $\delta_1$.  Neither of these have a Radon--Nikodym derivative so we have $J(0)=I_2(\delta_0)=\infty$ and $J(1)=I_2(\delta_1)=\infty$.  To show $J(\alpha)<\infty$ for each $\alpha\in(0,1)$ if suffices to find a measure $\mu_{\alpha}$ with $\mathrm{support}(\mu_{\alpha})\subseteq [0,1]$, $\ex(\mu_{\alpha}) =\alpha$ and $I_2(\mu_{\alpha})<\infty$.

For a fixed $\alpha$ let $\tilde{\alpha}=\min\{\alpha,1-\alpha\}$ and define $g_{\alpha}:\RR\longrightarrow[0,\infty)$ to be the piecewise linear function with $g_{\alpha}(x)=0$ for $x\leq\alpha-\tilde{\alpha}$ and $x\geq\alpha+\tilde{\alpha}$, and $g_{\alpha}(\alpha)=\displaystyle{\sqrt{\frac{3}{2\tilde{\alpha}}}}$.  From this we can define the probability measure $\mu_{\alpha}$ by $\displaystyle{\frac{\dd \mu_{\alpha}}{\dd x}=g^2_{\alpha}}$.  This measure satisfies our conditions since $\mathrm{support}(\mu_{\alpha})=[\alpha-\tilde{\alpha},\alpha+\tilde{\alpha}]\subseteq [0,1]$, $\ex(\mu_{\alpha}) =\alpha$ and
\begin{align*}
I_2(\mu_{\alpha})=\int_0^1 2x \left(\frac{\dd}{\dd x}g_{\alpha}(x)\right)^2 \,\dd x \leq 2 \left(\frac{1}{\tilde{\alpha}}\sqrt{\frac{3}{2\tilde{\alpha}}}\right)^2 <\infty,
\end{align*}
as required.

\item The claim is trivial for $\alpha=0$ or $\alpha=1$, so fix $\alpha\in(0,1)$ and find a sequence of measures $\{\mu_{\alpha,j}\}_{j\geq 1}$ with $\mathrm{support}(\mu_{\alpha,j})\subseteq [0,1]$,  $\ex(\mu_{\alpha,j}) =\alpha$ for each $j$, and such that $I_2(\mu_{\alpha,j})\longrightarrow J(\alpha)$ as $j\longrightarrow\infty$.  Because $[0,1]$ is compact we can take a weakly convergent subsequence $\mu_{\alpha,j_k}\longrightarrow\mu_{\alpha}$ say.  Since $I_2$ is lower semi-continuous we must have $I_2(\mu_{\alpha})\leq J(\alpha)$.  Therefore because we must also have $\mathrm{support}(\mu_{\alpha})\subseteq [0,1]$ and $\ex(\mu_{\alpha}) =\alpha$, then in fact $I_2(\mu_\alpha)=J(\alpha)$.  Thus it follows that for each $\alpha\in[0,1]$ the infimum of $I_2(\mu)$ over the set
\begin{align*}
\{\mu\in\mathcal{P}(\RR):\mathrm{support}\{\mu\}\subseteq [0,1] \text{ and }\ex(\mu) =\alpha\}
\end{align*}
is attained.  Now to show uniqueness we start with the claim that for any two probability measures $\mu_1$ and $\mu_2$ which are supported on $[0,1]$, and for any $\lambda,\tilde{\lambda}\in(0,1)$ with $\lambda+\tilde{\lambda}=1$ we have
\begin{align}
I_2(\lambda\mu_1+\tilde{\lambda}\mu_2)\leq \lambda I_2(\mu_1)+\tilde{\lambda}I_2(\mu_2) ,\label{eq:convex}
\end{align}
with equality if and only if $\mu_1=\mu_2$.

To proves this we first assume that $\mu_1$ and $\mu_2$ have respective Radon--Nikodym derivatives $m_1$ and $m_2$, otherwise the right hand side is infinite, and then use the form of $I_2$ given by (\ref{eq:Imu_int}) to calculate explicitly.
\begin{align}
I_2(\lambda\mu_{1}+\tilde{\lambda}\mu_{2}) &= \int_{0}^{1} 2x \left(\frac{\dd}{\dd x} \sqrt{\lambda m_1(x) + \tilde{\lambda}m_2(x)} \right)^2\,\dd x \nonumber \\
& =\int_{0}^{1} \frac{x}{2} \frac{\left(\lambda  m_1'(x)+\tilde{\lambda} m_2'(x)\right){}^2}{ \lambda  m_1(x)+\tilde{\lambda} m_2(x)}\,\dd x , \label{eq:convex_1}
\end{align}
and
\begin{align}
\lambda I_2(\mu_{1})+\tilde{\lambda}I_2(\mu_{2}) & = \int_{0}^{1} 2x\left( \lambda\left(\frac{\dd}{\dd x} \sqrt{m_1(x)}\right)^2 + \tilde{\lambda}\left(\frac{\dd}{\dd x} \sqrt{m_2(x)} \right)^2\right)\,\dd x \nonumber \\
& = \int_{0}^{1} \frac{x}{2} \frac{\lambda  m_2(x) m_1'(x)^2  + \tilde{\lambda}m_1(x) m_2'(x)^2 }{m_1(x) m_2(x)}\,\dd x. \label{eq:convex_2}
\end{align}
By subtracting (\ref{eq:convex_1}) from (\ref{eq:convex_2}) we get
\begin{align}
\lambda I_2(\mu_{1})+\tilde{\lambda}I_2(\mu_{2}) -I_2(\lambda\mu_{1}+\tilde{\lambda}\mu_{2}) = \int_{0}^{1} \frac{x}{2} \frac{\lambda  \tilde{\lambda} \Big(m_2(x) m_1'(x)-m_1(x) m_2'(x)\Big)^2}{m_1(x) m_2(x) \left(\lambda  m_1(x)+\tilde{\lambda} m_2(x)\right)} \,\dd x, \label{eq:convex_3}
\end{align}
which is greater than or equal to 0 with equality if and only if $m_2(x) m_1'(x)-m_1(x) m_2'(x)=0$ for all $x\in[0,1]$.  This would imply that $m_2$ is a constant multiple of $m_1$ and thus $\mu_1=\mu_2$.  Therefore the claim follows.

Having shown this, now suppose that $\mu_{1}$ and $\mu_{2}$ are two measures with $\ex(\mu_1)=\ex(\mu_2)=\alpha$ and $I_2(\mu_{1})=I_2(\mu_2)=J(\alpha)$.  Then because $\frac{1}{2}\mu_1+\frac{1}{2}\mu_2$ also has $\ex\left(\frac{1}{2}\mu_1+\frac{1}{2}\mu_2\right)=\alpha$, we then know that
\begin{align*}
J(\alpha)\leq I_2\left(\tfrac{1}{2}\mu_1+\tfrac{1}{2}\mu_2\right)\leq \tfrac{1}{2}I_2(\mu_1)+\tfrac{1}{2}I_2(\mu_2)=J(\alpha) .
\end{align*}
Hence we must have equality throughout, implying that $\mu_1=\mu_2$.

\item Let $0\leq\alpha_1<\alpha_2\leq 1$.  Then we need to show that $J(\lambda\alpha_1+\tilde{\lambda}\alpha_2)<\lambda J(\alpha_1)+\tilde{\lambda}J(\alpha_2)$ for each $\lambda,\tilde{\lambda}\in(0,1)$ with $\lambda+\tilde{\lambda}=1$.  To do this we first find $\mu_{1}$ and $\mu_{2}$ with $I_2(\mu_{1})=J(\alpha_1)$, $I_2(\mu_{2})=J(\alpha_2)$, $\ex(\mu_{1})=\alpha_1$ and $\ex(\mu_{2})=\alpha_2$.  Now $\ex(\lambda\mu_1+\tilde{\lambda}\mu_2)=\lambda\alpha_1+\tilde{\lambda}\alpha_2$, and so from (\ref{eq:convex}) we know that 
\begin{align*}
J(\lambda\alpha_1+\tilde{\lambda}\alpha_2)\leq I_2(\lambda\mu_1+\tilde{\lambda}\mu_2)< \lambda I_2(\mu_1)+\tilde{\lambda}I_2(\mu_2) = \lambda J(\alpha_1)+\tilde{\lambda}J(\alpha_2).
\end{align*}
The second inequality is strict since $\mu_1$ and $\mu_2$ have different means and therefore are not equal. As this holds for each $\lambda,\tilde{\lambda}\in(0,1)$ with $\lambda+\tilde{\lambda}=1$ and for each $0\leq\alpha_1<\alpha_2\leq 1$ the claim is proved.

\item Because $J$ is convex and finite for each $\alpha\in(0,1)$, then it must also be continuous on $(0,1)$.  Moreover, convexity implies that $J$ has left and right derivatives $\partial_- J$ and $\partial_+ J$.  In order to prove that $J$ is differentiable it suffices to show that these are always equal.  We shall do this by showing that at each point $\alpha\in(0,1)$ there exists a neighbourhood $\mathcal{N}_{\alpha}$ on which we can construct a differentiable function $f_{\alpha}$, with $f_{\alpha}(\alpha)=J(\alpha)$ and $f_{\alpha}(\beta)\geq J(\beta)$ for all $\beta\in\mathcal{N}_{\alpha}$.  Having done this we must then have $\partial_-J(\alpha)\geq f_{\alpha}'(\alpha)$ and $\partial_+J(\alpha)\leq f_{\alpha}'(\alpha)$.  However, because $J$ is convex we also have $\partial_-J(\alpha)\leq\partial_+J(\alpha)$, and so $\partial_-J(\alpha)=\partial_+J(\alpha)=f_{\alpha}'(\alpha)$.  Whence we see that $J$ is differentiable at $\alpha$ with $J'(\alpha)=f_{\alpha}'(\alpha)$.  Since the convexity of $J$ also implies that $J'(\alpha)$ is monotonically increasing then it must also be the case that $J'$ is continuous.

Given a fixed $\alpha\in(0,1)$ we now find the unique probability measure $\mu_{\alpha}$ supported on $[0,1]$ with $\ex(\mu_{\alpha})=\alpha$ and such that $I_2(\mu_{\alpha})=J(\alpha)$.  Because $I_2(\mu_{\alpha})<\infty$ then we know $\mu_{\alpha}$ has a Radon--Nikodym derivative, $\dfrac{\dd\mu_{\alpha}}{\dd x}=m_{\alpha}(x)$ say.  For each $\xi\in(-1,1)$ we also define a bijection $\varphi_{\xi}:[0,1]\longrightarrow [0,1]$ by $\varphi_{\xi}(x)= x^{1+\xi}$.  Using these we can then define a collection of measures $\{\mu_{\alpha,\xi}\}_{\xi\in(-1,1)}$ by $\mu_{\alpha,\xi}(A)=\mu_{\alpha}(\varphi_{\xi}^{-1}(A))$ for each measurable set $A\subseteq [0,1]$.  Note that each $\mu_{\alpha,\xi}$ has a Radon-Nikodym derivative given by $\dfrac{\dd \mu_{\alpha,\xi}}{\dd x}(x)=\dfrac{1}{(\varphi_{\xi}^{-1})'(x)}m_{\alpha}(\varphi_{\xi}^{-1}(x))$.

Now observe that $\ex(\mu_{\alpha,\xi})$ is a continuously differentiable function of $\xi$, and that
\begin{align*}
\frac{\dd}{\dd \xi} \ex (\mu_{\alpha,\xi} ) = \int_{0}^{1} \frac{\dd}{\dd \xi} x^{\frac{1}{1+\xi}}\,\dd\mu_{\alpha}(x) = \int_{0}^{1} \frac{1}{(1+\xi)^2} x^{\frac{1}{1+\xi}}\log \frac{1}{x}\,\dd\mu_{\alpha}(x) .
\end{align*}
Since $\displaystyle{\left.\frac{\dd}{\dd \xi} \ex (\mu_{\alpha,\xi} )\right|_{\xi=0}>0}$, then there is some neighbourhood $\mathcal{N}_{\alpha}$ of $\alpha$ on which a continuously differentiable inverse to $\xi\longmapsto\ex(\mu_{\alpha,\xi})$ exists.  We call this inverse $i_{\alpha}:\mathcal{N}_{\alpha}\longrightarrow(-1,1)$, and define $f_\alpha:\mathcal{N}_{\alpha}\longrightarrow\RR\cup\{\infty\}$ by $f_\alpha(\beta) = I_2(\mu_{\alpha,i_{\alpha}(\beta)})$.  Now observe that $f_{\alpha}(\alpha)=I_2(\mu_{\alpha,0})=J(\alpha)$, and because $\ex ( \mu_{\alpha,i_{\alpha}(\beta)} )=\beta$ then $f_{\alpha}(\beta) = I_2(\mu_{\alpha,i_{\alpha}(\beta)})\geq J(\beta)$.  Moreover, for $\beta\in\mathcal{N}_{\alpha}$ we have 
\begin{align*}
\frac{\dd}{\dd \beta} f_{\alpha}(\beta) = \frac{\dd}{\dd \beta} I_2(\mu_{\alpha,i_{\alpha}(\beta)})=\int_{0}^{1}\frac{\dd}{\dd \beta}\,2x \left(\frac{\dd}{\dd x}\sqrt{\frac{m_{\alpha}(\varphi^{-1}_{i_{\alpha}(\beta)} (x))}{(\varphi_{i_{\alpha}(\beta)}^{-1})'(x)}}\right)^2\,\dd x.
\end{align*}
Thus $f_{\alpha}$ is differentiable on $\mathcal{N}_{\alpha}$, and so we can conclude that $J$ is differentiable at $\alpha$.  Since $\alpha$ was arbitrary it then follows that $J$ is continuously differentiable on $(0,1)$ as claimed.

%\item Define $\mu^*$ by $\dfrac{\dd \mu^*}{\dd x} (x) = \left(\dfrac{\mathcal{J}_0(j_0\sqrt{x})}{\mathcal{J}_1(j_0)}\right)^2$ for $x\in[0,1]$.  Here $\mathcal{J}_{0}$ is a Bessel function of the first kind and $j_0$ is the first zero of $\mathcal{J}_0$.  One can then show that $\mu^*$ provides the unique infimum of $I_2(\mu)$ over the set of $\mu$ such that $\mathrm{support} ( \mu) \subseteq [0,1]$, and therefore there must be a unique infimum to $J$ at $\alpha=\ex\{\mu^*\}$.  Further calculation shows that $\ex(\mu^*)=(\gamma^*)^{-1}$, and that $J\left((\gamma^*)^{-1} \right)=\frac{1}{2}j_0^2\approx 2.892$. See Lemma \ref{lem:lemma9} for details.

\begin{figure}[p]
\begin{center}
\includegraphics[width=\textwidth]{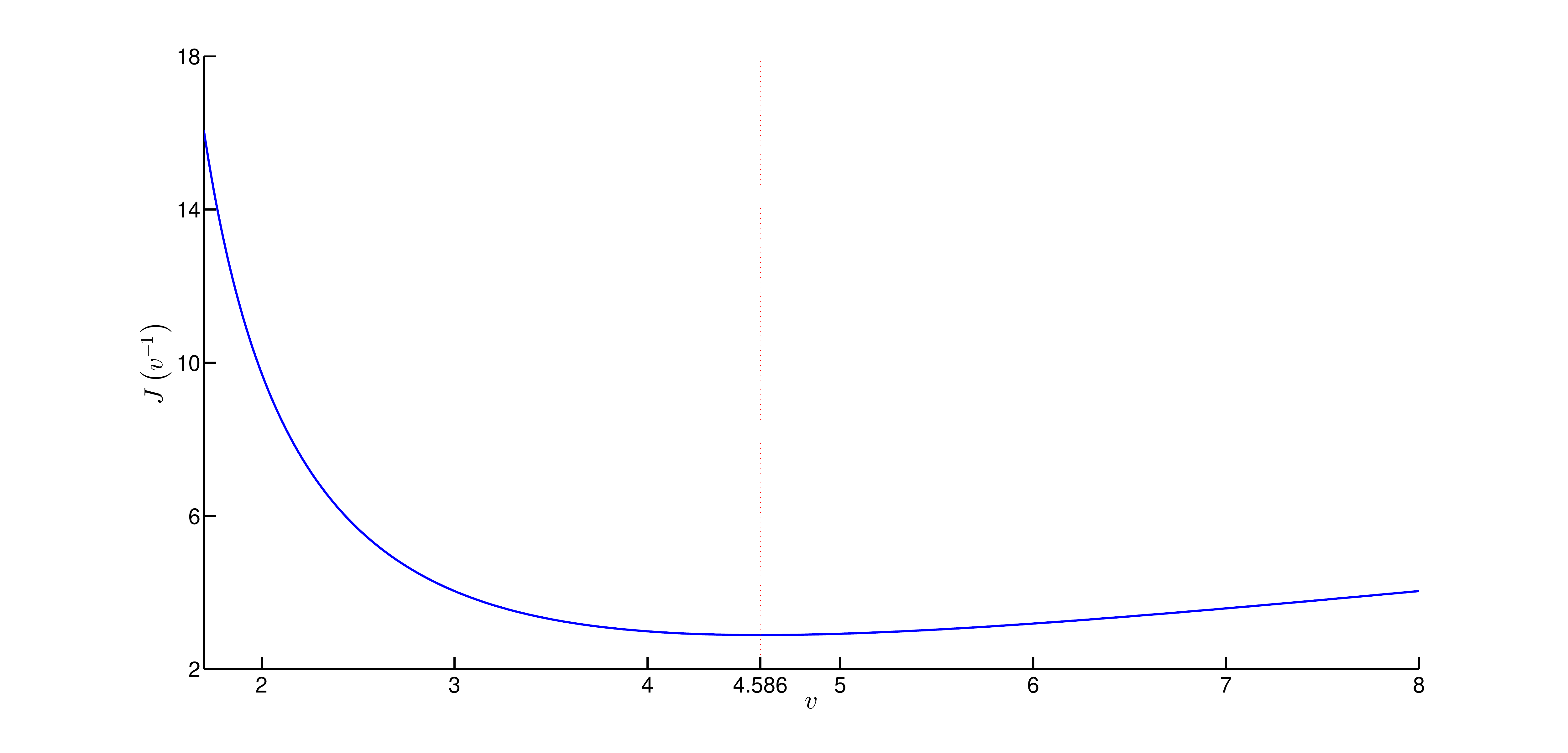}
\includegraphics[width=\textwidth]{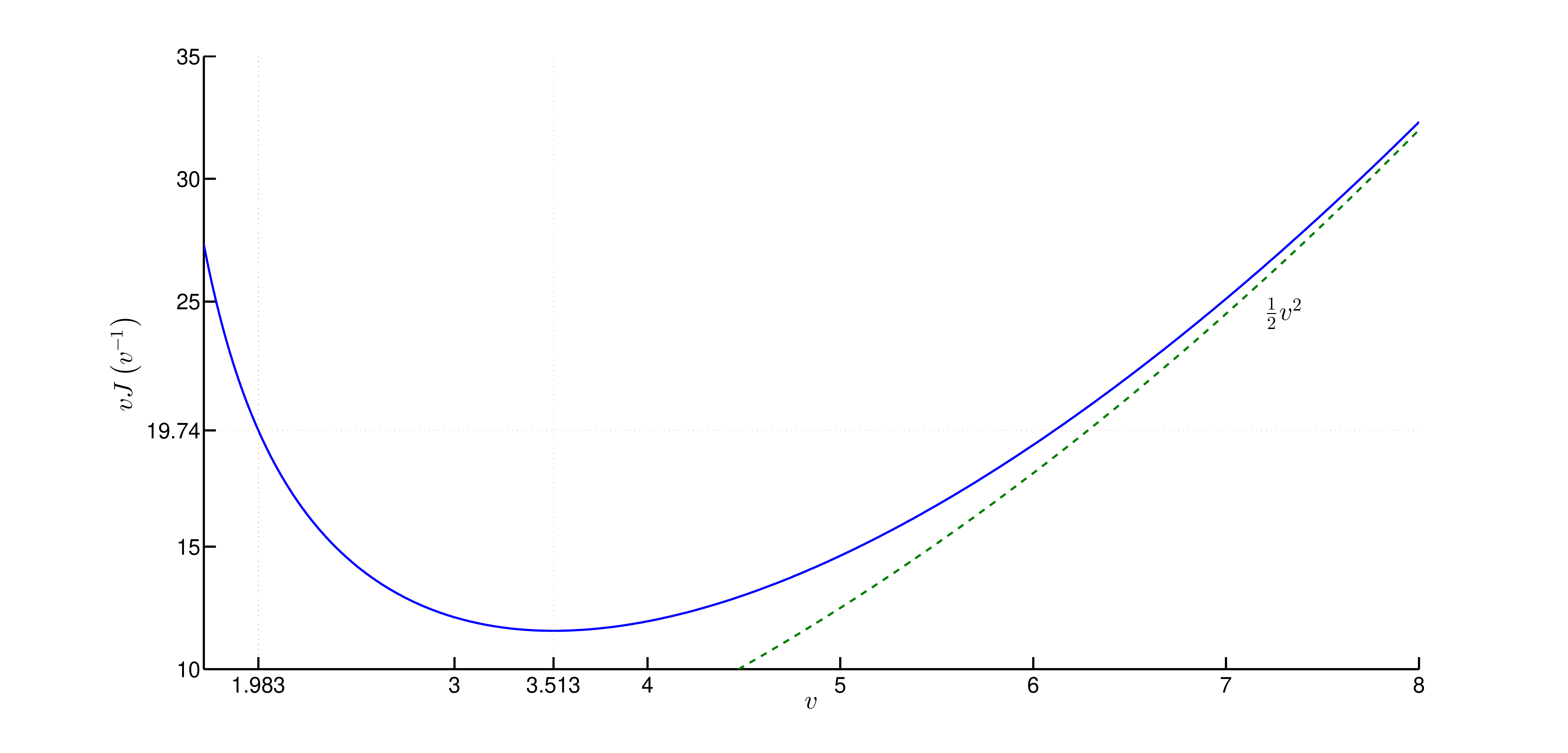}
\caption{Plots of $J\left(v^{-1}\right)$ and $v J\left(v^{-1}\right)$ against $v$.  Although we can not find a closed form for $J(\alpha)$, by using the by solving the Euler--Lagrange equations it is possible to plot $J(\alpha)$ by numerical methods.  Observe that $J\left(v^{-1}\right)$ has a unique minimum at $\gamma^{*}\approx 4.586$ and that the unique minimum of $vJ\left(v^{-1}\right)$ is attained at $\gamma^{\bullet}\approx 3.513 < \gamma^{*}$.
%\nl
%On the second graph we have included the curve $\frac{1}{2}v^2$ (dashed line), which Claim 7 of Lemma \ref{lem:J_props} shows is a lower bound for $vJ\left(v^{-1}\right)$.
Note also that the smallest root of $vJ\left(v^{-1}\right)=2\pi^2\approx 19.74$ is at $\gamma^{\circ}\approx 1.983$.  This value represents the minimal velocity for which we believe Conjecture \ref{con:rate} holds.
\nl
In Section \ref{sec:discussion} we describe $v J\left(v^{-1}\right)$ as being the exponential cost of a Brownian motion spending a unit of time traveling at speed $v$ and maintaining $L_x(t)\leq 1$.  Because (\ref{eq:speedcost}) shows that the unit time cost of a Brownian motion traveling at speed $v$ is $\frac{1}{2}v^2$, then this must be a lower bound for $v J\left(v^{-1}\right)$.}\label{fig:J}
\end{center}
\end{figure}

\item Because of the way that $J$ is defined then this is an immediate consequence of Lemma \ref{lem:lemma9}.

\item Fix $\varepsilon>0$ and suppose we have a probability measure $\mu$ with $\mathrm{support}(\mu)\subseteq[0,1]$ and $\ex(\mu)\geq 1-\varepsilon^2$.   Since $\ex(\mu)\leq (1-\varepsilon)\mu([0,1-\varepsilon])+\mu((1-\varepsilon,1])$ we must have $\mu([1-\varepsilon,1])\geq 1-\varepsilon$.  We can now use this fact to give a lower bound for $I_2(\mu)$.

For each $\delta>0$ define $u_{\varepsilon,\delta}:[0,1]\rightarrow(0,\infty)$ by
\begin{align*}
u_{\varepsilon,\delta}(x) = \left\{ \begin{array}{ll} 1+\delta & \text{for }0\leq x \leq 1-\varepsilon \\
\mathrm{sin}\left(\dfrac{\pi}{2\varepsilon}(1-x)\right) +\delta & \text{for }1-\varepsilon<x\leq 1 \end{array} \right. ,
\end{align*}
and now recall that the infinitesimal generator for a $\mathrm{BES}Q^2$ process is given by
\begin{align*}
\mathcal{L}_2 f(x) = 2x \frac{\dd^2}{\dd x^2}f(x) + 2 \frac{\dd}{\dd x}f(x).
\end{align*}
Now observe that if we set $u_{\varepsilon}(x)=\displaystyle{\lim_{\delta\rightarrow 0}}\,u_{\varepsilon,\delta}(x)$ then
\begin{align*}
\frac{\mathcal{L}_2 u_{\varepsilon}(x)}{u_{\varepsilon}(x)} \leq \left\{ \begin{array}{ll} 0 &\text{for }0\leq x \leq 1-\varepsilon \\
-2(1-\varepsilon)\left(\dfrac{\pi}{2\varepsilon}\right)^2 & \text{for }1-\varepsilon<x\leq 1 \end{array}\right. .
\end{align*}
Therefore because each $u_{\varepsilon,\delta}$ is in the domain of $\mathcal{L}_2$ and each $u_{\varepsilon,\delta}$ is strictly positive on $[0,1]$ we can use (\ref{eq:Ifunction1}) to get that
\begin{align*}
I_2(\mu)&= -\inf_{u\in,\mathcal{D},u>0}\int_{0}^1 \frac{\mathcal{L}_2(u)}{u}\dd\mu \geq \lim_{\delta\rightarrow 0} -\int_{0}^1 \frac{\mathcal{L}_2(u_{\varepsilon,\delta})}{u_{\varepsilon,\delta}}\dd\mu \\
&\geq \int_{1-\varepsilon}^1 2(1-\varepsilon) \left(\frac{\pi}{2\varepsilon}\right)^2\dd\mu(x) \geq 2 (1-\varepsilon)^2 \left(\frac{\pi}{2\varepsilon}\right)^2.
\end{align*}
Because this holds for each $\mu$ with $\mathrm{support}(\mu)\subseteq[0,1]$ and $\ex(\mu)\geq 1-\varepsilon^2$ we must have $J(1-\varepsilon^2)\geq \dfrac{\pi^2}{2}\dfrac{(1-\varepsilon)^2}{\varepsilon^2} $ for each $\varepsilon>0$.  It therefore follows that $J(\alpha)\longrightarrow\infty$ as $\alpha\longrightarrow 1$.

\item Suppose that $\alpha^{-1} J(\alpha)$ has minima at distinct points $0\leq\alpha_1<\alpha_2\leq 1$, and set $\tilde{\alpha}=\frac{1}{2}(\alpha_1+\alpha_2)$.  Now because $J$ is strictly convex we have
\begin{align*}
\frac{1}{\tilde{\alpha}}J(\tilde{\alpha}) &= \frac{2}{\alpha_1+\alpha_2}J\left(\frac{\alpha_1+\alpha_2}{2}\right)<\frac{1}{\alpha_1+\alpha_2}(J(\alpha_1)+J(\alpha_2)) \\
&= \frac{1}{\alpha_1+\alpha_2}\left(J(\alpha_1)+\frac{\alpha_2}{\alpha_1}J(\alpha_1)\right)=\frac{1}{\alpha_1}J(\alpha_1),
\end{align*}
contradicting the fact that $\alpha_1$ is a minimum.  Therefore the map $v\longmapsto v J\left(v^{-1}\right)$ must attain a unique minimum at some point $1\leq\gamma^{\bullet}\leq\infty$.  Furthermore, since claim 6 tells us that $J(\alpha)\longrightarrow\infty$ as $\alpha\nearrow 1$ then we know that $v J\left(v^{-1}\right)\longrightarrow\infty$ as $v\searrow 1$ implying that $\gamma^{\bullet}>1$.

To check that $\gamma^{\bullet}<\gamma^{*}$ we now use the fact that $J$ is continuously differentiable (claim 4) and so
\begin{align*}
\frac{\dd}{\dd v}v J\left(v^{-1}\right) = J\left(v^{-1}\right) + v \frac{\dd}{\dd v} J\left(v^{-1}\right) =J\left(v^{-1}\right)-\frac{1}{v}J'\left(v^{-1}\right) .
\end{align*}
Because $J$ is convex then we know that $\dfrac{\dd}{\dd v} J\left(v^{-1}\right)$ is increasing in $v$.  Therefore, as $\gamma^{*}$ is the minimum of $J$, it follows that $\dfrac{\dd}{\dd v}v J\left(v^{-1}\right)$ is increasing for $v\geq \gamma^{*}$.  Evaluating at $\gamma^{*}$ gives
\begin{align*}
\left.\frac{\dd}{\dd v}v J\left(v^{-1}\right)\right|_{v=\gamma^{*}} = J\left( (\gamma^{*})^{-1} \right)-\frac{1}{\gamma^{*}} J'\left( (\gamma^{*})^{-1} \right) = J\left( (\gamma^{*})^{-1} \right) > 0.
\end{align*}
Here $J'\left((\gamma^{*})^{-1}\right)=0$ because $(\gamma^{*})^{-1}$ is the unique minimum of $J$.  Thus the derivative of $v J\left(v^{-1}\right)$ is strictly positive for $v\geq\gamma^{*}$, and so its minimum must occur at a point $\gamma^{\bullet}<\gamma^{*}$.

Finally, to bound $\gamma^{\bullet} J\left((\gamma^{\bullet})^{-1}\right) $ we recall the explicit values for $\gamma^*$ and $J\left((\gamma^*)^{-1}\right)$ from Lemma \ref{lem:lemma9}.  From these we get
\begin{align*}
\gamma^{\bullet} J\left((\gamma^{\bullet})^{-1}\right) < \gamma^{*} J\left((\gamma^{*})^{-1}\right) =\frac{3}{1-2j_0^{-2}}\frac{1}{2}j_0^2 <2\pi^2 .
\end{align*}
Here $j_0$ denotes the first zero of the Bessel function $\mathcal{J}_0$, and $\dfrac{3}{1-2j_0^{-2}}\dfrac{1}{2}j_0^2 \approx 13.26$.
\qedhere
\end{enumerate}
\ep
\end{lemma}

To prove Lemma \ref{lem:mu_bullet} we now observe that if we put $\mu^{\bullet}=\mu_\alpha$ with $\alpha=(\gamma^{\bullet})^{-1}$ then for each $\mu$ supported on $[0,1]$ we have
\begin{align*}
\ex(\mu)^{-1}I_2(\mu)\geq \ex(\mu)^{-1} J(\ex(\mu))\geq\gamma^{\bullet} J\left((\gamma^{\bullet})^{-1}\right) = \ex(\mu^{\bullet})^{-1}I_2(\mu^{\bullet}) =\Gamma^{\bullet},
\end{align*}
with equality if and only if $\mu=\mu^{\bullet}$.
\hfill$\square$
\nl\vspace{-2\baselineskip}

\section{Proof of Theorem \ref{thm:reduced_speed}}\label{sec:proof1}

The proof of Theorem \ref{thm:reduced_speed} divides into two halves.  In the first we find estimates for $\WW(W_t\,|\,\EE^{\bullet}_T)$ -- from which we can deduce that the limiting process is ballistic. The second half then deals with the weak convergence of $\WW(\cdot\,|\,\EE_T^{\bullet})$ as $T\longrightarrow\infty$.

\subsection{Ballistic behaviour}\label{sec:ballistic}

Lemma \ref{lem:mu_bullet} tells us that there is a unique measure $\mu^\bullet$ which minimises $\ex(\mu)^{-1}I_2(\mu)$ over all probability measures supported on $[0,1]$.  We know that $\ex(\mu^{\bullet})^{-1}=\gamma^{\bullet}$ where $1<\gamma^\bullet<\gamma^*$, and so our aim is to use this fact to show that for each $\varepsilon>0$ there is a $T_{\varepsilon}>0$ such that for all $T_\varepsilon \leq t \leq T$ we have
\begin{align}
\WW\left(\left| \frac{W_t}{t}-\gamma^{\bullet}  \right|\geq\varepsilon \left. \right| \EE_T^{\bullet} \right) < \varepsilon .\label{eq:speed2}
\end{align}
Provided the limit $\displaystyle{\QQ^{\bullet}=\lim_{T\rightarrow\infty}\WW(\cdot\,|\,\EE_T)}$ exists -- which we will show in Section \ref{sec:conv} -- we can then deduce that 
\begin{align}
\lim_{t\rightarrow\infty}\frac{W_t}{t} = \gamma^{\bullet} \label{eq:speed}
\end{align}
in $\QQ^{\bullet}$-probability.
Since the left hand side of (\ref{eq:speed2}) can be written as
\begin{align}
\WW\left(\left| \frac{W_t}{t}-\gamma^\bullet \right|\geq\varepsilon \left. \right| \EE_T^{\bullet} \right) = \frac{\WW(\EE_T^{\bullet}\cap \{| W_t-\gamma^{\bullet}t|\geq\varepsilon t\})}{\WW(\EE^{\bullet}_T)} ,\label{eq:speedratio}
\end{align}
then our goal is to show that the ratio between $\WW(\EE_T^{\bullet}\cap\{| W_t-\gamma^{\bullet}t|\geq\varepsilon t\})$ and $\WW(\EE^{\bullet}_{T})$ is small when $t$ and $T$ are large.  For this we rely on the Ray--Knight Theorems and the calculations of Section \ref{sec:calc}.  Recall that in Section \ref{sec:ray--knight} we defined $S^-_T=\int_0^T \mathbbm{1}_{\{W_s<0\}}\,\dd s$ and $S^+_T=\int_0^T \mathbbm{1}_{\{W_s>W_T\}}\,\dd s$.  Theorem \ref{thm:ray_knight_fixed_t} tells us that the local time profile $(L_x(T))_{x\in\RR}$ can be described by two $\mathrm{BES}Q^0$ processes, with integrals $S^-_T$ and $S^+_T$, joined to a $\mathrm{BES}Q^2$ bridge with integral $T-S^-_T-S^+_T$.
More specifically, when we condition on $\{W_T\geq 0\}$ then the 5-tuple $(W_T,L_{W_T}(T),L_0(T),S^-_T,S^+_T)$ admits a density on $[0,\infty)^5$, and when we condition on $(W_T,L_{W_T}(T),L_0(T),S^-_T,S^+_T)=(a,b,c,s^-,s^+)$ we have
\vspace{-0.5\baselineskip}\begin{itemize}
\item $(Y^-_x)_{x\geq 0} = L_{-x}(T)$ is a $\mathrm{BES}Q^0(c)$ process conditioned to have integral equal to $s^-$.
\item $(Y_x)_{0\leq x\leq W_T} = L_{x}(T)$ is a $\mathrm{BES}Q^2_{a}(c,b)$ bridge conditioned to have integral equal to $T-s^- - s^+$.
\item $(Y^+_x)_{x\geq 0} = L_{a+x}(T)$ is a $\mathrm{BES}Q^0(b)$ process conditioned to have integral equal to $s^+$.
\end{itemize}
Now observe that $\EE^{\bullet}_T=\{L_x(T)\leq 1$ for all $x\in\RR\}$ is exactly the event that $Y^-$, $Y$ and $Y^+$ are all bounded above by 1.  Therefore, if we use $\mathbbm{Y}_c^0$ to denote the law of a $\mathrm{BES}Q^0(c)$ process and $\mathbbm{Y}^2_{a,c,b}$ to denote the law of a $\mathrm{BES}Q^2_{a}(c,b)$ bridge, then the first step to controlling (\ref{eq:speedratio}) is to understand how the probabilities
\begin{align}
Q_0(c,s) &= \mathbbm{Y}^0_c\left(Y_x\leq 1 \text{ for all } x\geq 0 \,|\, \textstyle{\int_0^\infty Y_x \,\dd x =s}\right) \label{eq:Q0} \\
\text{and}\quad Q_2(a,c,b,s) &= \mathbbm{Y}^2_{a,c,b}\left(Y_x\leq 1 \text{ for all }0\leq x\leq a\,|\,\textstyle{\int_0^a Y_x \,\dd x =s}\right) \label{eq:Q2}
\end{align}
depend on $c,b\in[0,1)$ and $a,s\geq 0$.  Before estimating these it will be useful for us to construct an auxiliary process which is parametrised by the integral of $(Y_x)_{x\geq 0}$.

\begin{definition}[Auxiliary process]\label{def:auxiliary}
Suppose we are given a $\mathrm{BES}Q^d(c)$ process $(Y_x)_{x\geq 0}$ with integral $S=\int^\infty_0 Y_x\,\dd x\in[0,\infty)\cup\{\infty\}$.  For each $0\leq t \leq S$ we set $\rho(t) = \inf\left\{u:\int_0^u Y_x\,\dd x = t\right\}$, then $\rho:[0,S)\longrightarrow[0,\infty)$ is a strictly increasing differentiable function.  Using $\rho$ we define the coupled auxiliary process $(Z_t)_{t\geq 0}$ by $Z_t = Y_{\rho(t \wedge S)}$, and denote its law by $\ZZ^d_c$.
\end{definition}

\begin{lemma}
$(Z_t)_{t \geq 0}$ is a Markov process whose infinitesimal generator, $\tilde{\mathcal{L}}_d$, is given by
\begin{align}
\tilde{\mathcal{L}}_d f(x) = \frac{1}{x}\mathcal{L}_d f(x) = 2\frac{\mathrm{d}^2}{\mathrm{d}x^2}f(x) +d \frac{1}{x} \frac{\dd}{\dd x}f(x) . \label{eq:aux_gen}
\end{align}
Here $\mathcal{L}_d$ denotes the generator of the $\BESQ{d}{c}$ process $(Y_x)_{x\geq0}$.  Recall (\ref{eq:BESQgen}).
\bp
We can only have $S<\infty$ if $d=0$.  In this situation we have $Z_t = Y_{\rho(S)} = 0$ for all $t\geq S$, and so from now on we will assume that $0\leq t< S$.  Because $\rho$ is an increasing function of $t$, then $(Z_t)_{t\geq 0}$ must inherit the Markov property from $(Y_x)_{x\geq 0}$. To show (\ref{eq:aux_gen}) we observe that $\displaystyle{\int^{\rho(t)}_0 Y_x \,\dd x = t}$, and differentiate both sides with respect to $t$ to get
\begin{align*}
Y_{\rho(t)}\,\frac{\dd \rho}{\dd t}(t) =1\quad \Longrightarrow \quad \dfrac{\dd\rho}{\dd t}(t) = \dfrac{1}{Z_t} .
\end{align*}
Since $\rho(0)=0$ we can approximate $\rho(t)$ near $0$ by $\rho(t)\sim  t\,\rho^{\prime}(0)$.  Therefore  making the substitution $s = t \,\rho^{\prime}(0)$ gives
\begin{align*}
\tilde{\mathcal{L}}_d f(x) &= \lim_{t\rightarrow 0}\frac{\ex^x\{ f(Z_t) \} - f(x)}{t} = \lim_{t\rightarrow 0}\frac{\ex^x\{ f(Y_{\rho(t)}) \} - f(x)}{t} \\ &= \lim_{t\rightarrow 0}\frac{\ex^x\{ f(Y_{t\,\rho^{\prime}(0)}) \} - f(x)}{t} =\rho^{\prime}(0) \lim_{s\rightarrow 0}\frac{\ex^x\{ f(Y_s) \} - f(x)}{s} = \frac{1}{x} \mathcal{L}_d f(x)
. \qedhere
\end{align*}
\ep
\end{lemma}

If $I_d$ is the Donsker--Varadhan rate function of the occupation measure of $(Y_x)_{x\geq 0}$, then we can deduce $\tilde{I}_d$, the rate function of the occupation measure of $(Z_t)_{t\geq 0}$, as follows:  Assume that a probability measure ${\mu}$ has $0<\int_0^\infty x^{-1}\,\dd {\mu}(x) < \infty$ (else $\tilde{I}_d({\mu})=\infty$).  Because the work of Donsker and Varadhan, \cite{d&v}, tells us that $\tilde{I}_d({\mu})$ is finite if and only if ${\mu}$ has a continuous Radon--Nikodym derivative $\dfrac{\dd {\mu}}{\dd x}$ in the domain of $\tilde{\mathcal{L}}_d$, then we may assume that $\dfrac{\dd {\mu}}{\dd x}$ exists.  Now for each such ${\mu}$ we can construct a tilted measure $\psi(\mu)$, with Radon--Nikodym derivative given by
\begin{align}
\frac{\dd\psi(\mu)}{\dd x}(x) = \frac{1}{\int_0^\infty x^{-1}\,\dd {\mu}(x)}\frac{1}{x}\frac{\dd {\mu}}{\dd x}(x) . \label{eq:tiltedmeasure}
\end{align}
One can then check that $\dfrac{\dd \psi(\mu)}{\dd x}$ is in the domain of $\mathcal{L}_d$ and that  
\begin{align*}
\ex(\psi(\mu)) = \int_{0}^{\infty} x\,\dd \psi(\mu)(x) = \dfrac{1}{\int_0^{\infty} x^{-1}\,\dd {\mu}(x)}.
\end{align*}   
By recalling the definition of the Donsker--Varadhan rate function given by (\ref{eq:ratefunction}), we can now find $\tilde{I}_d(\mu)$ is terms of $\psi(\mu)$.
\begin{align}
\tilde{I}_d({\mu}) = & -\inf_{{u\in C^2(0,\infty),  u>0}} \int_0^\infty \frac{1}{x}\frac{\mathcal{L}_d(u(x))}{u(x)}\frac{\dd {\mu}}{\dd x}(x)\,\dd x \nonumber \\= &-\inf_{{u\in C^2(0,\infty),  u>0}} \int_0^\infty \frac{\mathcal{L}_d(u(x))}{u(x)}\frac{1}{\ex(\psi(\mu))}\frac{\dd{\psi(\mu)}}{\dd x}(x)\,\dd x = \frac{1}{\ex(\psi(\mu))} I_d({\psi(\mu)}) .\label{eq:ratefunctiontilde}
\end{align}

Using the auxiliary process gives a useful tool for adapting the Donsker--Varadhan Theorem to the occupation measure of a square Bessel process stopped at the random time $\rho(s)$.

\begin{lemma}\label{lem:modifieddv}
Suppose that $(Y_x)_{x\geq 0}$ is a $\mathrm{BES}Q^d(c)$ process with law $\mathbbm{Y}_c^d$.
%Recall that we use $I_d$ to denote the Donsker--Varadhan rate function for the occupation measure of $(Y_x)_{x\geq 0}$, and  that $\rho(s)=\inf\left\{u:\int_0^u Y_x\,\dd x \geq s\right\}$.
We then have
\begin{align}
\limsup_{s\rightarrow\infty} \frac{1}{s} \log\mathbbm{Y}_c^d ( L((Y_x),\rho(s),\cdot)\in C)&\leq - \inf_{\mu\in C} \ex(\mu)^{-1} I_d(\mu) \label{eq:exIfunction1} \\
\liminf_{S\rightarrow\infty} \frac{1}{s} \log\mathbbm{Y}_c^d ( L((Y_x),\rho(s),\cdot)\in O)&\geq - \inf_{\mu\in O} \ex(\mu)^{-1} I_d(\mu) \label{eq:exIfunction2}
\end{align}
for each closed set $C\subseteq \mathcal{P}(\RR)$, and each open set $O\subseteq \mathcal{P}(\RR)$ such that $c\in\mathrm{support}(\mu)$ for each $\mu\in O$.
\end{lemma}

\begin{remark}
It was noted in Remark \ref{rem:donsker-varadhan} that the Donsker--Varadhan Theorem continues to hold when the start-point, $Y_0$, is a real valued random variable.  Since Lemma \ref{lem:modifieddv} is a consequence of Theorem \ref{thm:donsker_varadhan}, then (\ref{eq:exIfunction1}) and (\ref{eq:exIfunction2}) will also hold when $c$ is random.  In this case we substitute $\mathbbm{Y}_c^d(L((Y_x),\rho(s),\cdot)\in \cdot)$ with $\ex^c\big\{\mathbbm{Y}_c^d(L((Y_x),\rho(s),\cdot)\in \cdot)\big\}$, and replace the condition that $\big\{c\in\mathrm{support}(\mu)$ for each $\mu \in O\big\}$ with the condition that $\big\{\mathrm{support}(c)\subseteq\mathrm{support}(\mu)$ for each $\mu \in O\big\}$.
\end{remark}

\bp[Proof of Lemma \ref{lem:modifieddv}]
Fix $s>0$, then given $(Y_x)_{0\leq x\leq \rho(s)}$ let $(Z_t)_{0\leq t\leq s}$ be the auxiliary process to $Y$.  Suppose $\mu = L((Y_x),\rho(s),\cdot)$ is the occupation measure of $(Y_x)_{0\leq x\leq \rho(s)}$, and use $\varphi(\mu)=L((Z_t),s,\cdot)$ to denote the occupation measure of $(Z_t)_{0\leq t\leq s}$.  Since $Y$ and $Z$ are both continuous Markov processes then $\mu$ and $\varphi(\mu)$ have Radon--Nikodym derivatives almost surely.  Furthermore, by using $\dfrac{\dd\rho(t)}{\dd t} = \dfrac{1}{Y_{\rho(t)}}$ we can check that
\begin{align}
\frac{\dd\varphi({\mu})}{\dd x} (x) = \frac{1}{\int_0^\infty x\,\dd\mu(x)}x\frac{\dd\mu}{\dd x} (x) = \frac{1}{\ex(\mu)}x\frac{\dd\mu}{\dd x} (x). \label{eq:tiltedinverse}
\end{align}
Observe that (\ref{eq:tiltedinverse}) gives a continuous inverse to (\ref{eq:tiltedmeasure}), and thus $\mu\longmapsto\varphi(\mu)$ is a bicontinuous bijection with  inverse $\varphi(\mu)\longmapsto\psi(\varphi(\mu))=\mu$.  As $\varphi$ is continuous then $\varphi(C)$ is also a closed set, and so we can use the Donsker--Varadhan Theorem and the natural coupling between $Y$ and $Z$ to get
\begin{align*}
&\limsup_{s\rightarrow\infty} \frac{1}{s} \log\mathbbm{Y}_c^d \big( L((Y_x),\rho(s),\cdot)\in C\big) = \limsup_{s\rightarrow\infty} \frac{1}{s} \log\mathbbm{Z}_c^d \big( L((Z_x),s,\cdot)\in \varphi(C)\big) \\
&\qquad\leq -\inf_{{\mu}\in\varphi (C)} \tilde{I}_d({\mu}) =  -\inf_{{\mu}\in\varphi (C)} \ex(\psi(\mu))^{-1}{I}_d(\psi({\mu})) = -\inf_{\mu\in C} \ex(\mu)^{-1}I_d(\mu).
\end{align*}
Similarly when $O$ is open then so too is $\varphi(O)$ and therefore
\begin{align*}
&\liminf_{s\rightarrow\infty} \frac{1}{s} \log\mathbbm{Y}_c^d \big( L((Y_x),\rho(s),\cdot)\in O\big) = \liminf_{s\rightarrow\infty} \frac{1}{s} \log\mathbbm{Z}_c^d \big( L((Z_x),s,\cdot)\in \varphi(O)\big) \\
&\qquad\geq -\inf_{{\mu}\in\varphi (O)} \tilde{I}_d({\mu}) =  -\inf_{{\mu}\in\varphi (O)} \ex(\psi(\mu))^{-1}{I}_d(\psi({\mu})) = -\inf_{\mu\in O} \ex(\mu)^{-1}I_d(\mu). \qedhere
\end{align*}
\ep 

Because it is clear that $Y_x\leq 1$ for all $x$ if and only if $Z_t=Y_{\rho(t)}\leq 1$ for all $t$, then we can use the auxiliary process to  generalise (\ref{eq:Q2}) to the case where the bridge length $a$ is not fixed.  From (\ref{eq:tiltedmeasure}) we know that
\begin{align*}
\frac{\dd L(Y,a,\cdot)}{\dd x}(x) = \frac{1}{x} \frac{1}{\int x^{-1} \,\dd L(Z,s,x)}\frac{\dd L(Z,s,\cdot)}{\dd x}(x) .
\end{align*}
Therefore from the condition that  $\int^a_0 Y_x\,\dd x = s$ we get
\begin{align}
&\int_0^a Y_x\,\dd x = a \int x \,\dd L(Y,a,x) =  \frac{a}{\int x^{-1} \,\dd L(Z,s,x)} \int \frac{\dd L(Z,s,x)}{\dd x} =s \label{eq:ex_}\\
&\quad\Longrightarrow  s \int x^{-1} \,\dd L(Z,s,x)= \int_0^s \frac{1}{Z_t}\,\dd t  = a \nonumber.
\end{align}
This allows us to rewrite (\ref{eq:Q2}) in terms of the auxiliary process $(Z_t)_{0\leq t\leq s}$ as
\begin{align*}
{Q}_2(a,c,b,s) = \mathbbm{Z}^2_{s,c,b}\left(Z_t\leq 1 \text{ for all }0\leq t\leq s\,|\, \textstyle{\int_0^s Z_t^{-1}\,\dd t = a}\right).
\end{align*}
Here $\mathbbm{Z}^2_{s,c,b}$ denotes $\mathbbm{Z}^2_c(\cdot\,|\,Z_s = b)=\displaystyle{\lim_{\varepsilon\rightarrow\infty}\ZZ^2_c(\cdot\,|\,|Z_s-b|<\varepsilon)}$.  ${Q}_2(a,c,b,s)$ can now be generalised by letting $A\subseteq[0,\infty)$ be a measurable set and defining
\begin{align}
\tilde{Q}_2(A,c,b,s) &= \mathbbm{Z}^2_{s,c,b}\left(Z_t\leq 1\text{ for all } 0\leq t\leq s \text{ and } \textstyle{\int_0^s Z_t^{-1}\,\dd t \in A}\right) .\label{eq:Q2tilde} \\
\tilde{Q}_2(c,b,s) &= \tilde{Q}_2([0,\infty),c,b,s)=\mathbbm{Z}^2_{s,c,b}\left(Z_t\leq 1\text{ for all } 0\leq t\leq s \right) .\label{eq:Q2tilde__}
\end{align}
Hence if we suppose that $(Y_x)_{0\leq x\leq \mathrm{len}}$ is a $\mathrm{BES}Q_{ \mathrm{len}}^2(c,b)$ bridge of undetermined length, which is conditioned to have $\int_0^{ \mathrm{len}} Y_x\,\dd x = s$, then the probability that its length $ \mathrm{len}$ is in $A$ and that $Y_x\leq 1$ for all $0\leq x\leq  \mathrm{len}$ is given by $\tilde{Q}_2(A,c,b,s)$.  Here it is assumed that $c,b\in[0,1)$, $s> 0$ and $A\subseteq[0,\infty)$ are all fixed.

\begin{remark}
As we know that $\EE_T^\bullet=\{L_{-x}(T)\leq 1$ for all $x\geq 0\}\cap\{L_x(T)\leq 1$ for all $0\leq x\leq W_T\}\cap\{L_{W_T+x}(T)\leq 1$ for all $x\geq 0\}$, then when we condition on $L_0(T)=c$, $L_{W_T}(T)=b$, $S^-=\int_0^\infty L_{-x}(T)\,\dd x=s^-$ and $S^+=\int_0^\infty L_{W_T+x}(T)\,\dd x=s^+$ we get
\begin{align}
\WW(\EE^\bullet_T\,|\,c,b,s^-,s^+)=Q_0(c,s^-)\times \tilde{Q}_2(c,b,T-s^--s^+)\times Q_0(b,s^+).\label{eq:EETinQs}
\end{align}
$\WW(\EE_T^\bullet)$ can now be found by integrating (\ref{eq:EETinQs}) with respect to the joint law of $L_0(T)$, $L_{W_T}(T)$, $S^-$ and $S^+$.
\end{remark}

\begin{lemma}\label{lem:Q0}
$Q_0(c,s)$ is decreasing as a function of both $c\in(0,1)$ and $s> 0$.  Furthermore, for each fixed $c\in(0,1)$ we have
\begin{align}
\lim_{s\rightarrow \infty} \frac{1}{s}\log Q_0(c,s)= -2\pi^2. \label{eq:Q0est}
\end{align}
\end{lemma}

\begin{lemma}\label{lem:Q2}
If we assume that $c$ and $b$ are random variables supported on $[0,1)$, that $V_C\subseteq[1,\infty)$ is closed, and use $sV_C$ to denote $\{a:s^{-1}a\in V_C\}$, then
\begin{align}
\limsup_{s\rightarrow\infty} \frac{1}{s} \log \tilde{Q}_2(sV_C,c,b,s) \leq \displaystyle{-\inf_{v\in V_C} vJ\left(v^{-1}\right)},\label{eq:Y_closed}
\end{align}
where $J$ is given by Definition \ref{def:J}.  Conversely, if $V_O\subseteq[1,\infty)$ is an open set then
\begin{align}
\liminf_{s\rightarrow\infty} \frac{1}{s} \log \tilde{Q}_2(sV_O,c,b,s) \geq \displaystyle{-\inf_{v\in V_O} vJ\left(v^{-1}\right)}.\label{eq:Y_open}
\end{align}
\end{lemma}

\bp[Proof of Lemma \ref{lem:Q0}]
Since $(Y_x)_{x\geq 0}$ is a $\BESQ{0}{c}$ process conditioned on $\int_0^\infty Y_x\,\dd x=s$, then the auxiliary process  $(Z_t)_{0\leq t\leq s}$ has starting point $Z_0 = c$ and first hits 0 at $t=s$.
We have already noted that $Y_x\leq 1$ for all $x\geq 0$ if and only if $Z_t\leq 1$ for all $0\leq t\leq s$.  Therefore by coupling two $\mathrm{BES}Q^0$ processes $Y$ and $\tilde{Y}$ with respective starting points  $0<c\leq\tilde{c}<1$ in such a way that their auxiliary process $Z$ and $\tilde{Z}$ satisfy $Z_t\leq\tilde{Z}_t$ for all $0\leq t\leq s$, it is easy to deduce that
\begin{align*}
Q_0(c,s) = \mathbbm{Y}_c^0\left(Y_x\leq 1  \,|\,\textstyle{\int_0^\infty} Y_x\,\dd x = s\right) \geq \mathbbm{Y}_{\tilde{c}}^0\left(\tilde{Y}_x\leq 1 \,|\,\textstyle{\int_0^\infty} \tilde{Y}_x\,\dd x = s\right) = Q_0(\tilde{c},s) .
\end{align*}
Thus $Q_0(c,s)$ is decreasing as a function of $c\in(0,1)$.  The monotonicity of $Q_0(c,s)$ as a function of $s$ follows via a similar coupling argument.

Now define $E_C=\{\mu:\mathrm{support}(\mu)\subseteq(0,1]\}$ and $E_O=\{\mu:\mathrm{support}(\mu)\subseteq(0,1)\}$, and note that $E_C\subseteq\mathcal{P}((0,\infty))$ is weakly closed and $E_O\subseteq\mathcal{P}((0,\infty))$ is weakly open.  If we fix $c\in(0,1)$ and let $s\geq 1$, then we can use equation (\ref{eq:f_density}) to get the bounds
\begin{align*}
\frac{\ee^{-1}}{\sqrt{8\pi}} \frac{c}{s^{\frac{3}{2}}}  \leq \mathbbm{Y}^0_c\left(\textstyle{\int_0^\infty Y_x\,\dd x \in [s,s-1]}\right)
\quad\text{and}\quad
\mathbbm{Y}^0_c\left(\textstyle{\int_0^\infty Y_x\,\dd x \in [s,s+1]}\right)\leq \frac{c}{s^{\frac{3}{2}}} 
\end{align*}
for all fixed $c\in(0,1)$ and all $s\geq 1$.  Therefore because $Q_0(c,s)$ is decreasing as a function of $s$ we have
\begin{align*}
Q_0(c,s) &
%=\mathbbm{Y}_c^0\left(Y_x\leq 1\,|\,\textstyle{\int_0^\infty Y_x\,\dd x = s}\right)
\leq \mathbbm{Y}_c^0\left(Y_x\leq 1 \text{ for all } x\geq 0\,|\,\textstyle{\int_0^\infty Y_x\,\dd x \in [s-1,s]}\right)\\
&= \frac{\mathbbm{Y}^0_c\left(Y_x\leq 1 \text{ for all } x\geq 0 \text{ and } \int_0^\infty Y_x\,\dd x \in [s-1,s]\right)}{\mathbbm{Y}^0_c\left(\int_0^\infty Y_x\,\dd x \in [s-1,s]\right)} \\
&\leq \sqrt{8\pi}\ee\, s^{\frac{3}{2}}c^{-1}\, \mathbbm{Y}^0_c\left(Y_x\leq 1 \text{ for all } x\geq 0 \text{ and } \textstyle{\int_0^\infty Y_x\,\dd x \geq s-1}\right) \\
&= {\sqrt{8\pi}\ee} \,s^{\frac{3}{2}} c^{-1}\, \mathbbm{Y}^0_c\big(L((Y_x),\rho(s-1),\cdot)\in E_C\big),
\end{align*}
and similarly
\begin{align*}
Q_0(c,s) & \geq \mathbbm{Y}_c^0\left(Y_x\leq 1 \text{ for all } x\geq 0\,|\,\textstyle{\int_0^\infty Y_x\,\dd x \in [s,s+1]}\right) \\
& \geq s^{\frac{3}{2}} c^{-1}\, \mathbbm{Y}^0_c\big(L((Y_x),\rho(s),\cdot)\in E_O\big).
\end{align*}
As we know that $\displaystyle{\lim_{s\rightarrow\infty} s^{-1} \log s^{\frac{3}{2}}c^{-1} =0}$, 
then we can now get bounds for $\displaystyle{\lim_{s\rightarrow\infty}\frac{1}{s}}\log Q_0(c,s)$ by using Lemma \ref{lem:modifieddv}.
\begin{align*}
\limsup_{s\rightarrow\infty} \frac{1}{s}\log Q_0(c,s) &\leq \limsup_{s\rightarrow\infty} \frac{1}{s}\log \mathbbm{Y}_0^c(L((Y_x),\rho(s-1),\cdot)\in E_C) \\ &\leq -\inf_{\mu\in E_C} \frac{1}{\ex(\mu)}I_0(\mu) = -2\pi^2 ,\\
\liminf_{s\rightarrow\infty} \frac{1}{s}\log Q_0(c,s) &\geq \liminf_{s\rightarrow\infty} \frac{1}{s}\log \mathbbm{Y}_0^c(L((Y_x),\rho(s),\cdot)\in E_O) \\ &\geq -\inf_{\mu\in E_O} \frac{1}{\ex(\mu)}I_0(\mu) = -2\pi^2.
\end{align*}
Here the equality $\displaystyle{\inf_{\mu\in E_C} \frac{1}{\ex({\mu})}I_0(\mu)=\inf_{\mu\in E_O} \frac{1}{\ex({\mu})}I_0(\mu)=2\pi^2}$ comes from Lemma \ref{lem:I_0}.
\ep

\bp[Proof of Lemma \ref{lem:Q2}]
Let $E_{V_C} = \big\{\mu:\mathrm{support}(\mu)\subseteq[0,1]$ and $\int_{0}^\infty x^{-1}\,\dd \mu(x) \in V_C\big\}$ and note that if $L((Z_t),s,\cdot)\in E_{V_C}$ then we must have $Z_t\leq 1$ for all $0\leq t\leq s$.  Furthermore, because we also have that $\int^s_0 Z_t^{-1}\,\dd t = s\int_0^\infty x^{-1}\,\dd L((Z_t),s,x)$ then we see that the probability being estimated in (\ref{eq:Y_closed}) is exactly $\ZZ^2_{s,c,b}\left(L((Z_t),s,\cdot)\in E_{V_C}\right)$.  Since $E_{V_C}$ is closed with respect to the weak topology, then our aim is to estimate this probability by using the Donsker--Varadhan Theorem.  However, in order to do this we must first check that $\ZZ^2_{s,c,b}\left(L((Z_t),s,\cdot)\in E_{V_C}\right)=\displaystyle{\lim_{\varepsilon\rightarrow 0}}\,\ZZ^2_c\left(L((Z_t),s,\cdot)\in E_{V_C}\,|\,|Z_s - b|<\varepsilon\right)$ is comparable with $\ZZ^2_{c}\left(L((Z_t),s,\cdot)\in E_{V_C}\right)$.

By following the methods of Pinsky, \cite{pinsky1}, we know that when we condition 
on the occupation measure of a diffusion being contained within a collection of measures of uniformly bounded support (such as when we condition on $\{L((Z_t),s,\cdot)\in E_{V_C}\}_{s\geq 1}$)
%on a diffusion satisfying a sequence of  events of the form $\{L((Z_t),s,\cdot)\in E_{V_C}\}_{s\geq 1}$,
then the distribution of the end point, $Z_s$, will converge to an absolutely continuous random variable as $s\longrightarrow\infty$.  Therefore there must exists a $K_1<\infty$ (depending on $V_C$) such that
\begin{align}
\ZZ^2_c\big(|Z_s-b|<\varepsilon \,|\, L((Z_t),s,\cdot)\in E_{V_C}\big) <\varepsilon \, K_1
\end{align}
for all $s\geq 1$ and all $b\in[0,1]$.  By inspecting the infinitesimal generator of $Z$, as given by (\ref{eq:aux_gen}), we see that $(Z_s)_{s\geq 0}$ is equal in law to a (rescaled) Bessel process.  Therefore there must be constants $0<k_2<K_2<\infty$ such that
\begin{align}
k_2 < \frac{s^{\frac{3}{2}}}{\varepsilon (b+\varepsilon)^2}\mathbbm{Z}^2_c (|Z_{s}-b|<\varepsilon) <K_2, \label{eq:comparek2}
\end{align}
for all $\varepsilon>0$, $b\in[0,1)$ and $s\geq 1$.   Hence we can deduce that
\begin{align}
\lim_{\varepsilon\rightarrow 0}\,\ZZ^2_c (L((Z_t),s,\cdot)\in E_{V_C}\,|\,|Z_s-b|<\varepsilon\big) \leq \frac{K_1}{k_2} \frac{s^{\frac{3}{2}}}{b^2}\, \ZZ^2_{c}\left(L((Z_t),s,\cdot)\in E_{V_C}\right). \label{eq:comparedmeasures}
\end{align}
Since $\displaystyle{\lim_{s\rightarrow\infty}\frac{1}{s}\log\frac{K_1}{k_2}\frac{s^{\frac{3}{2}}}{b^2}=0}$, then by applying the Donsker--Varadhan Theorem we get
\begin{align*}
\limsup_{s\rightarrow\infty} \frac{1}{s}\log \tilde{Q}_2(sV_C,c,b,s) & \leq \limsup_{s\rightarrow\infty} \frac{1}{s}\log \ZZ^2_{c}\left(L((Z_t),s,\cdot)\in E_{V_C}\right) \\
&\leq - \inf_{\mu\in E_{V_C}} \tilde{I}_2(\mu) = -\inf_{\mu\in E_{V_C}} \ex(\psi(\mu))^{-1} I_2(\psi(\mu)) \\ &=  -\inf_{\mu\in \varphi(E_{V_C})} \ex(\mu)^{-1} I_2(\mu).
\end{align*}
Where $\psi$ and $\varphi$ are given by (\ref{eq:tiltedmeasure}) and (\ref{eq:tiltedinverse}).  It is easy to verify that $\varphi(E_{V_C}) = \{\mu:\mathrm{support}(\mu)\subseteq [0,1]$ and $\ex(\mu)^{-1}\in V_C\}$, and so 
%(\ref{eq:Y_closed}) follows from the definition of $J$.  See Definition \ref{def:J}.
\begin{align*}
-\inf_{\mu\in\varphi(E_{V_C})}\ex(\mu)^{-1}I_2(\mu) &= -\inf_{v\in V_C} v \, \inf\{I_2(\mu):\mathrm{support}(\mu)\subseteq[0,1]\text{ and }\ex(\mu)^{-1}=v\} \\ &= -\inf_{v\in V_C} v J\left(v^{-1}\right),
\end{align*}
proving (\ref{eq:Y_closed}).

Now set $E_{V_O}=\big\{\mu:\mathrm{support}(\mu)\subseteq[0,1)$ and $\ex(\mu)^{-1}\in V_O\big\}$ and note that if $L((Z_t),s,\cdot)\in E_{V_O}$ then $(Z_t)_{t\geq 0}$ must satisfy $Z_t\leq 1$ for all $0\leq t\leq s$ and $s^{-1}\int_0^s Z_t^{-1}\,\dd t \in V_O$.  In order to show (\ref{eq:Y_open}) we also note that for each $b\in [0,1)$ we can find a $k_3$ (depending on $b$ and $V_O$) such that
\begin{align}
\varepsilon \, k_3 < \ZZ^2_c\big(|Z_s-b|<\varepsilon \,|\, L((Z_t),s,\cdot)\in E_{V_O}\big)
\end{align}
for all $s\geq 1$ and all $\varepsilon$ sufficiently small.  By combining this with (\ref{eq:comparek2}) we get
\begin{align}
\lim_{\varepsilon\rightarrow 0}\,\ZZ^2_c (L((Z_t),s,\cdot)\in E_{V_O}\,|\,|Z_s-b|<\varepsilon\big) \geq \frac{k_3}{K_2} \frac{s^{\frac{3}{2}}}{b^2}\, \ZZ^2_{c}\left(L((Z_t),s,\cdot)\in E_{V_O}\right). \label{eq:comparedmeasures2}
\end{align}
Therefore, as $E_{V_O}$ is open with respect to the weak topology, (\ref{eq:Y_open}) follows by applying the Donsker--Varadhan Theorem.
\begin{align*}
\liminf_{s\rightarrow\infty} \frac{1}{s}\log \tilde{Q}_2(sV_O,c,b,s) & \geq \liminf_{s\rightarrow\infty} \frac{1}{s}\log \ZZ^2_{c}\left(L((Z_t),s,\cdot)\in E_{V_O}\right) \\
&\geq - \inf_{\mu\in E_{V_O}} \tilde{I}_2(\mu) = -\inf_{\mu\in E_{V_O}} \ex(\psi(\mu))^{-1} I_2(\psi(\mu)) \\ &=  -\inf_{\mu\in \varphi(E_{V_O})} \ex(\mu)^{-1} I_2(\mu) = - \inf_{v\in V_O} v J\left( v^{-1} \right). \qedhere
\end{align*}
\ep

At this point we can now show that for each $\varepsilon > 0$ there exists a $T_{\varepsilon}$ such that for all $T\geq T_{\varepsilon}$ we have
\begin{align}
\frac{\mathbbm{W}(\mathcal{E}^{\bullet}_T\cap\{|W_T - \gamma^{\bullet} T|\geq \varepsilon T\})}{\mathbbm{W}(\EE^{\bullet}_T)} < \varepsilon .\label{eq:speedjustT}
\end{align}
This is precisely (\ref{eq:speedratio}) in the case where $t=T$.  To do this we need an upper bound for $\mathbbm{W}(\mathcal{E}^{\bullet}_T\cap\{|W_T - \gamma^{\bullet} T|\geq \varepsilon T\})$ and a lower bound for $\mathbbm{W}(\EE^{\bullet}_T)$, and hence it suffices to prove the following two claims.

\begin{claim}\label{claim:claim2}
We have 
\begin{align*}
\displaystyle{\liminf_{T\rightarrow\infty} \frac{1}{T} \log \mathbbm{W}(\EE^{\bullet}_T) \geq - \Gamma^{\bullet}},
\end{align*}
where $\displaystyle{\Gamma^{\bullet} = \gamma^{\bullet}J\left((\gamma^{\bullet})^{-1}\right)}$ is the minimal value of $v J\left( v^{-1}\right)$.  See page \pageref{def:J}.
\end{claim}

\begin{claim}\label{claim:claim1}
For each $\varepsilon > 0$ there exists $k>0$ such that
\begin{align*}
\limsup_{T\rightarrow\infty} \frac{1}{T} \log \mathbbm{W}(\mathcal{E}^{\bullet}_T\cap\{|W_T - \gamma^{\bullet} T|\geq \varepsilon T\}) \leq -\Gamma^{\bullet} - k.
\end{align*}
\end{claim}

\begin{remark}
By combining Claim \ref{claim:claim2} and Claim \ref{claim:claim1} we get
\begin{align*}
\limsup_{T\rightarrow\infty} \frac{1}{T} \log \frac{\mathbbm{W}(\mathcal{E}^{\bullet}_T\cap\{|W_T - \gamma^{\bullet} T|\geq \varepsilon T\})}{\mathbbm{W}(\EE^{\bullet}_T)}  \leq - k,
\end{align*}
and so (\ref{eq:speedjustT}) follows immediately.
\end{remark}

\bp[Proof of Claim \ref{claim:claim2}]
Suppose we use $f_T(\cdot,\cdot)$ to denote the joint probability distribution of $L_0(T)$ and $L_{W_T}(T)$.  Observe that if we fix $T=1$ then, because $f_1(\cdot,\cdot)$ is continuous and strictly positive on $[0,\infty)^2$, by compactness there must exist $0<k_1<K_1<\infty$ such that $k_1\leq f_1(c,b)\leq K_1$ for all $c,b\in[0,1]$.  We also note that the scaling property of Brownian motion implies that $f_T(c,b) = T^{-1} f_1\left(cT^{-\frac{1}{2}}, b T^{-\frac{1}{2}}\right)$ for all $T>0$ and $c,b\in[0,\infty)$.  Therefore it must follow that
\begin{align}
\frac{k_1}{T}\leq f_T(c,b) \leq \frac{K_1}{T} \label{eq:localtimedensity}
\end{align}
for all $c,b\in[0,1]$ and all $T\geq 1$.  We can now observe that 
\begin{align}
\WW(\EE^\bullet_T)&= \int_{[0,\infty)^2}\WW(\EE^\bullet_T\,|\,L_0(T) = c, L_{W_T}(t) = b) f_T(c,b)\,\dd c\,\dd b \nonumber \\
&\geq \frac{k_1}{T} \int_{\left[\frac{1}{4},\frac{3}{4}\right]^2} \WW(\EE^\bullet_T\,|\, L_0(T) = c, L_{W_T}(T) = b)\,\dd c\,\dd b \nonumber\\
&\geq \frac{k_1}{4 T} \inf_{c,b\in\left[\frac{1}{4},\frac{3}{4}\right]} \WW(\EE^{\bullet}_T\,|\,L_0(T) = c, L_{W_T} = b) . \label{eq:EET_estimate_1}
\end{align}
From now on let $\hat{c}$ and $\hat{b}$ be the values of $c,b\in\left[\frac{1}{4},\frac{3}{4}\right]$ which minimise (\ref{eq:EET_estimate_1}), and use $g_{T,\hat{c},\hat{b}}(\cdot,\cdot)$ to denote the joint probability density of $S^- = \int_0^\infty L_{-x}(T)\,\dd x$ and $S^+ = \int_{0}^{\infty} L_{W_T + x}(T)\,\dd x$ with respect to $\WW(\cdot\,|\, L_0(T) = \hat{c},L_{W_T}(T) = \hat{b})$.  Brownian scaling also tells us that $g_{T,\hat{c}\sqrt{T},\hat{b}\sqrt{T}}(s^-,s^+) = T^{-2} g_{T,\hat{c},\hat{b}}(s^- T^{-1},s^+ T^{-1})$, therefore, since decreasing $\hat{c}$ and $\hat{b}$ would only make it more likely for $S^-$ and $S^+$ to be small, there must be some $k_2>0$ such that
\begin{align}
\frac{k_2}{T^2} < \WW(S^-\leq 1 \text{ and }S^+\leq 1\,|\, L_0(T)=\hat{c}, L_{W_T}(t)=\hat{b}) . \label{eq:S^-bound}
\end{align}
By combining this with (\ref{eq:EET_estimate_1}), and recalling that (\ref{eq:Q0}) and (\ref{eq:Q2tilde__}) give the probabilities of $(L_{-x}(T))_{x\geq 0}$, $(L_{W_T + x}(T))_{x\geq 0}$ and $(L_x(T))_{0\leq x\leq W_T}$ being bounded above by 1, we get
\begin{align*}
\WW(\EE^\bullet_T)&\geq\frac{k_1 T}{4 k_2}\WW(\EE_T^{\bullet}\,|\,L_0(T)=\hat{c},L_{W_T}(T)=\hat{b},S^-\leq 1 \text{ and }S^+\leq 1) \\
&\geq \frac{k_1 T}{4 k_2} \,\inf_{0\leq s^-\leq 1}Q_0(\hat{c},s^-)\times \inf_{0\leq s^+\leq 1}Q_0(\hat{b},s^+) \times \inf_{0\leq s\leq T} \tilde{Q}_2(\hat{c},\hat{b},s).
\end{align*}
Since Lemma \ref{lem:Q0} tells us that $Q_0(c,s)$ is decreasing as a function of both $c$ and $s$ then we must have $\displaystyle{ \inf_{0\leq s^-\leq 1}Q_0(\hat{c},s^-)}, \displaystyle{ \inf_{0\leq s^+\leq 1}Q_0(\hat{b},s^+)}\geq  \textstyle{Q_0\left(\frac{3}{4},1\right) > 0}$.  Therefore when we take the logarithm the only term to contribute is $\displaystyle{\inf_{0\leq s\leq T} \tilde{Q}_2(\hat{c},\hat{b},s)}$.  We can now complete the proof by using Lemma \ref{lem:Q2}.
\begin{align*}
&\liminf_{T\rightarrow \infty}\frac{1}{T}\log\WW(\EE^\bullet_T) \geq \liminf_{T\rightarrow\infty}\frac{1}{T}\log \tilde{Q}_2(\hat{c},\hat{b},T) \geq -\inf_{v\in[0,\infty)} v J\left(v^{-1}\right) = -\Gamma^{\bullet}. \qedhere
\end{align*}
\ep

\bp[Proof of Claim \ref{claim:claim1}]
Suppose we have conditioned on the event $\{W_T\geq 0\}$, and use $(L_x(T))_{x\in\RR}$ to denote the local time profile of $(W_t)_{t\geq 0}$ at time $T$.  Recall that $S = \int_{0}^{W_T} L_x(T)\,\dd x$, then in order to estimate $\mathbbm{W}(\mathcal{E}^{\bullet}_T\cap\{|W_T - \gamma^{\bullet} T|\geq \varepsilon T\})$ it will be useful for us to let $\eta = \eta(\varepsilon)>0$ be a positive constant -- which we shall determine later -- and consider the cases $S > T(1-\eta)$ and $S\leq T(1-\eta)$ separately.
By the Law of Total Probability we have
\begin{align*}
\WW\left(\EE^{\bullet}_T \cap\{| W_T-\gamma^{\bullet}T|\geq\varepsilon T\}\right) =& \WW\left(\EE^{\bullet}_T\cap\{| W_T-\gamma^{\bullet}T|\geq\varepsilon T\}\cap \{S> T\left(1-\eta\right)\}\right) \\ &\quad+ \WW\left(\EE^{\bullet}_T\cap \{| W_T-\gamma^{\bullet}T|\geq\varepsilon T\}\cap \{S\leq T\left(1-\eta\right)\}\right) ,
\end{align*}
and therefore we need to show that both of these terms are sufficiently small.

If we condition on $L_0(T)=c$, $L_{W_T}(T)=b$ and $\int_0^{W_T}L_x(T)\,\dd x = s$ then $(Z_t)_{t\geq 0}$, the auxiliary process to $(L_x(T))_{0\leq x\leq W_T}$, has  $Z_0=c$, $Z_s=b$ and $\int^s_0 Z_t^{-1}\,\dd t = W_T$.  Therefore if $T(1-\eta)<s\leq T$ and $|W_T-\gamma^\bullet_T|\geq \varepsilon T$ then
\begin{align}
\left|\frac{W_T}{T}-\gamma^\bullet\right|=\left|\frac{1}{T}\int^s_0 Z_t^{-1}\,\dd t - \gamma^\bullet\right|>\varepsilon\quad\Longrightarrow\quad\left|\frac{1}{s}\int^s_0 Z_t^{-1}\,\dd t - \gamma^\bullet\right|>\frac{\varepsilon}{2} \label{eq:boundawayfromgamma}
\end{align}  
provided $s$ is sufficiently close to $T$, i.e. provided $\eta$ is sufficiently small.  Assume that $\eta$ is small enough for (\ref{eq:boundawayfromgamma}) to hold and set $V_O=\big\{v\in[0,\infty) : |v-\gamma^\bullet|>\frac{1}{2}\varepsilon\big\}$.  Since the event $\EE^{\bullet}_T$ implies that $Z_t\leq 1$ for all $0\leq t\leq s$, then we can use Lemma \ref{lem:Q2} to estimate $\WW\left(\EE^{\bullet}_T\cap\{| W_T-\gamma^{\bullet}T|\geq\varepsilon T\}\cap \{S> T\left(1-\eta\right)\}\right)$.
\begin{align*}
&\limsup_{T\rightarrow\infty}\frac{1}{T}\log\WW\left(\EE^{\bullet}_T\cap\{| W_T-\gamma^{\bullet}T|\geq\varepsilon T\}\cap \{S> T\left(1-\eta\right)\}\right) \\
&\quad \leq\limsup_{T\rightarrow\infty} \frac{1}{T}\log \left(\sup_{\substack{T(1-\eta)<s<T \\b,c\in [0,1)}}\tilde{Q}_2(s V_O,c,b,s) \right) \leq - (1-\eta) \inf_{v \in V_O} v J\left(v^{-1}\right)  .
\end{align*}
From Lemma \ref{lem:J_props} we also know that $v J\left(v^{-1}\right)$ attains a unique minimum value of $\Gamma^{\bullet}$ at $v = \gamma^{\bullet}$.  Therefore, by continuity, there must exist some $\eta>0$ and some ${k_1}>0$ such that
\begin{align*}
- (1-\eta) \inf_{v \in V_O} v J\left(v^{-1}\right) \leq - \Gamma^{\bullet} -{k_1} .
\end{align*}
Hence we have controlled the first term.   We now fix such a value of $\eta$ and move on to considering the case where $0\leq S\leq T\left(1-\eta\right)$.  Recall that we defined $S^- = \int_{-\infty}^0 L_x(T)\,\dd x$ and $S^+ = \int_{W_T}^{\infty} L_x(T)\,\dd x$, and that when we condition on $L_0(T)=c$, $L_{W_T}(T)=b$, $S^- = s^-$ and $S^+ = s^+$ then $(L_{-x}(T))_{x\geq 0}$ and $(L_{W_T+x}(T))_{x\geq 0}$ are equal in law to $\BESQ{0}{c}$ and $\BESQ{0}{b}$ processes conditioned to have integrals equal to $s^-$ and $s^+$ respectively.
The probabilities that these processes are bounded above by 1 are given by $Q_0(c,s^-)$ and $Q_0(b,s^+)$.  Thus since we know that $S^- + S + S^+ = T$ we get 
\begin{align*}
&\WW(\EE^\bullet_T\cap\{|W_T-\gamma^\bullet T|\geq\varepsilon T\}\cap\{S\leq T(1-\eta)\}) \leq \WW(\EE^\bullet_T\cap\{S\leq T(1-\eta)\}) \\
&\quad\leq\sup_{\substack{0\leq s\leq T(1-\eta) \\ s^-+s+s^+ = T}}\left(\sup_{c\in[0,1)}Q_0(c,s^-)\times \sup_{b\in[0,1)}Q_0(b,s^+) \times \sup_{c,b\in[0,1)}\tilde{Q}_2([0,\infty),c,b,s)\right) .
\end{align*}
From  Lemma \ref{lem:Q0} and Lemma \ref{lem:Q2} we know that
\begin{align*}
\limsup_{s\rightarrow\infty} \frac{1}{s}\log\sup_{c,b\in[0,1)}\tilde{Q}_2([0,\infty),c,b,s) \leq - \Gamma^\bullet\quad\text{and}\quad  \limsup_{s \rightarrow\infty}\frac{1}{s}\log\sup_{c\in[0,1)} Q_0(c,s) -2\pi^2 . 
\end{align*}
Therefore, because $\Gamma^\bullet<2\pi^2$, we get
\begin{align*}
\limsup_{T\rightarrow \infty} \frac{1}{T}\log\WW(\EE^\bullet_T\cap\{S\leq T(1-\eta)\})&\leq\sup_{0\leq \beta\leq 1-\eta} - \beta \,\Gamma^{\bullet} - (1-\beta) \,2\pi^2 = - \Gamma^{\bullet} - k_2
\end{align*}
for some $k_2 > 0$.  The claim is now proved by setting $k = \min \{k_1,k_2\}$.
\ep

Having shown that (\ref{eq:speedjustT}) holds for all $T\geq T_{\varepsilon}$ we can now show (\ref{eq:speedratio}) by noting that
\begin{align*}
&\frac{\mathbbm{W}(\mathcal{E}^{\bullet}_T\cap\{|W_t - \gamma^{\bullet}t|\geq \varepsilon t\})}{\mathbbm{W}(\EE^{\bullet}_T)} = \frac{\mathbbm{W}(\EE^{\bullet}_T\,|\, \EE_t^{\bullet}\cap\{|W_t - \gamma^{\bullet}t|\geq \varepsilon t\})}{\mathbbm{W}(\EE^{\bullet}_T\,|\, \EE_{t}^{\bullet})}\frac{\mathbbm{W}(\mathcal{E}^{\bullet}_t\cap\{|W_t - \gamma^{\bullet} t|\geq \varepsilon t\})}{\mathbbm{W}(\EE^{\bullet}_t) } ,
\end{align*}
and then checking that both $\mathbbm{W}(\EE^{\bullet}_T\,|\, \EE_t^{\bullet}\cap\{W_t - \gamma^{\bullet}t|\geq \varepsilon t\})$ and $\mathbbm{W}(\EE^{\bullet}_T\,|\, \EE_{t}^{\bullet})$ are comparable with $\WW(\EE_{T-t}^{\bullet})$.  This is the content of the following two lemmas.

\begin{lemma}\label{claim:claim3}
Let $t\geq 0$ and suppose $A$ is any event which is measurable with respect to $\mathcal{F}_t$.  We must then have $\WW(\EE^{\bullet}_T\,|\,A)\leq \WW(\EE^{\bullet}_{T-t})$ for all $T\geq t$.
\end{lemma}

\begin{lemma}\label{claim:claim4}
There exists $\eta > 0$ such that
\begin{align}
\eta \,\WW(\EE_{T-t}^{\bullet})\leq\WW(\EE_T^{\bullet}\,|\,\EE_t^\bullet) \label{eq:EETeta}
\end{align}
for all $0\leq t\leq T$.
\end{lemma}

\begin{remark}
Having proved these two lemmas we then have that
\begin{align}
\frac{\mathbbm{W}(\EE^{\bullet}_T\cap\{|W_t - \gamma^{\bullet}t|\geq \varepsilon t\})}{\mathbbm{W}(\EE^{\bullet}_T)} \leq \frac{1}{\eta} \frac{\mathbbm{W}(\EE^{\bullet}_t\cap\{|W_t - \gamma^{\bullet}t|\geq \varepsilon t\})}{\mathbbm{W}(\EE^{\bullet}_t)}  .
\end{align}
Because we can conclude from Claim \ref{claim:claim1} and Claim \ref{claim:claim2} that the ratio on the right hand side becomes arbitrarily small as $t\longrightarrow\infty$, then it follows that we can always find a $T_{\varepsilon}$ such that for all $T_{\varepsilon} \leq t \leq T$ we have
\begin{align*}
\WW\left(\left|\frac{W_t}{t} - \gamma^{\bullet}\right|\geq \varepsilon\,|\,\EE_T^{\bullet}\right) < \varepsilon ,
\end{align*}
proving (\ref{eq:speed2}).  Therefore, once we have proved Lemma \ref{claim:claim3} and Lemma \ref{claim:claim4}, we will have shown that the limiting process is ballistic with speed $\gamma^{\bullet}$.
\hfill$\square$
\end{remark}

\bp[Proof of Lemma \ref{claim:claim3}]
If $(W_s)_{s\geq 0}$ is a Brownian motion and $t>0$ is fixed, then we know that $(W_{s+t}-W_t)_{s\geq 0}$ is also a Brownian motion which is independent of $(W_s)_{0\leq s\leq t}$.  Therefore $(L_{W_t + x}(T) - L_{W_t + x}(t))_{x\in\RR}$ is independent of $(L_x(t))_{x\in\RR}$ and is equal in law to $(L_x(T-t))_{x\in\RR}$.  Now if we are given $(W_s)_{0\leq s\leq t}$ for some $(W_s)_{s\geq 0}$ in $A$ then the event $(W_s)_{s\geq 0}\in\EE_T^{\bullet}$ is precisely the event that $L_{W_T + x}(T) - L_{W_T + x}(t) \leq 1 - L_{W_T + x}(t)$ for all $x\in\RR$.  Because $1-L_{W_T+x}(t)\leq 1$ for all $x\in\RR$ then we get the upper bound
\begin{align*}
\WW(\EE_T^\bullet\,|\,A) &\leq\WW(L_{W_T+x}(T) - L_{W_T+x}(t)\leq 1 \text{ for all }x\in\RR) = \WW(\EE^\bullet_{T-t}). \qedhere
\end{align*}
\ep

The proof of Lemma \ref{claim:claim4} is somewhat more tricky and uses the following.

\begin{lemma}\label{lem:below0}
Given a Brownian motion $(W_t)_{t\geq 0}$ and a fixed $T>0$ we let $S^-=\int_0^\infty L_{-x}(T)\,\dd x$ denote the amount of time that $(W_t)_{0 \leq t\leq T}$ spends below 0.  Suppose we condition  $(W_t)_{t\geq 0}$ on $\EE_{T}^{\bullet}$, then $S^-$ has an exponential tail whose exponent does not depend on $T$.
In other words there are universal constants $\xi>0$ and $\Xi<\infty$ such that for each $T>0$ and all $a\geq 0$ we have
\begin{align}
\WW\big(S^-=\textstyle{\int^\infty_0L_{-x}(T)\,\dd x > a\,|\,\EE_{T}^{\bullet}}\big)&<\Xi\,\ee^{-\xi a} . \label{eq:tail_time}
\end{align}
\end{lemma}
\bp
Recall that $\tilde{Q}_2(c,b,s) $ gives us the probability that the auxiliary process $(Z_t)_{t\geq 0}$, conditioned to have $Z_0=c$ and $Z_s=b$, satisfies $Z_t\leq 1$ for all $0\leq t\leq s$.  Lemma \ref{lem:Q2} allows us to understand how $\tilde{Q}_2(c,b,s)$ behaves as $s\longrightarrow\infty$, but we would also like to understand how $\tilde{Q}_2(c,b,s+\Delta s)$ relates to $\tilde{Q}_2(c,b,s)$.  Consider the measure
\begin{align*}
\nu_{S,\Delta S,c,b} = \mathbbm{Z}^2_c(Z_s\in\cdot\,|\,Z_t\leq 1\text{ for all }0\leq t\leq s+\Delta s \text{ and }Z_{s+\Delta s}= b).  
\end{align*}
By standard coupling arguments one can show that the measure $\nu=\displaystyle{\lim_{c,b\rightarrow 1}}\nu_{1,1,c,b}$ stochastically dominates $\nu_{S,\Delta S, c,b}$ for all $c,b\in[0,1)$ and each $s,\Delta s\geq 1$.  It is also simple to use coupling arguments to show that $\tilde{Q}_2(c,b,s) $ is a decreasing function of $c$, $b$ and $s$.  Therefore by using the fact that $(Z_s)_{s\geq 0}$ is a Markov process and applying Chebyshev's sum inequality we get
\begin{align*}
&\tilde{Q}_2(c,b,s+\Delta s) = \int_0^1 \tilde{Q}_2(c,y,\Delta s) \tilde{Q}_2(y,b,s) \,\dd\nu_{s,\Delta s,c,b}(y) 
\geq \int_0^1 \tilde{Q}_2(c,y,\Delta s)\tilde{Q}_2(y,b,s) \,\dd\nu(y) \\
&\qquad\geq \int_0^1 \tilde{Q}_2(c,y,\Delta s) \,\dd\nu(y) \times \int_0^1 \tilde{Q}_2(y,b,s) \,\dd\nu(y) \geq k_1 \, \tilde{Q}_2(c,b,s)\times \int_0^1 \tilde{Q}_2(c,y,\Delta s)\,\dd\nu(y) ,
\end{align*}
for some $k_1$ which is independent of $c$, $b$ and $s$.  Since Lemma \ref{lem:Q2} tells us how $\tilde{Q}_2(c,y,\Delta s)$ behaves (for large $\Delta s$), therefore we can deduce that
\begin{align}
\liminf_{\Delta s\rightarrow\infty} \frac{1}{\Delta s} \inf_{\substack{s\geq 1 \\ c,b\in[0,1)}} \log \frac{\tilde{Q}_2(c,b,s+\Delta s)}{\tilde{Q}_2(c,b,s)} \geq -\inf_{v\in[0,\infty)} v J\left(v^{-1} \right) = -\Gamma^{\bullet} . \label{eq:DeltaS}
\end{align}

Now suppose we condition on $L_0(T)=c$, $L_{W_T}(T)=b$ and $ (L_{W_T+x}(T))_{x\geq 0}=(Y^+_x)_{x\geq 0}$ for some $c,b\in[0,1]$ and some $(Y^+_x)_{x\geq 0}$ with $Y^+_x\leq 1$ for all $x\geq 0$ and $\int_0^\infty Y^+_x\,\dd x = s^+$, $0\leq s^+\leq T$.  Use $\mu_{T,s^+,c,b}$ to denote the measure of $\WW\big(S^- = \int_{0}^{\infty}L_{-x}(T)\,\dd x \in \cdot \,|\, W_T>0, L_0(T)=c, L_{W_T}(T)=b$ and $S^+ = \int_0^\infty L_{W_T+x}(T)\,\dd x  = s^+\big)$, then by arguing in the same way as we did for (\ref{eq:S^-bound}) we can find a $k_2>0$ such that
\begin{align}
\mu_{T,s^+,c,b}([0,1]) > \frac{k_2}{T}
\end{align}
for all $c,b\in[0,1)$ and all $0\leq s^+\leq T$.  We then have
\begin{align*}
\WW(\EE_T^\bullet\,|\,c,b,Y^+) &= \int_0^{T-s^+} Q_0(s,c)\, \tilde{Q}_2(T-s^+-s,c,b) \, \dd \mu_{T,s^+,c,b}(s) \\
&\geq \frac{k_2}{T} \inf_{s\in[0,1]}Q_0(s,c)\, \tilde{Q}_2(T-s^+,c,b),
\end{align*}
and
\begin{align*}
\WW(\EE_T^\bullet\cap\{S^->a\}\,|\,c,b,Y^+) &= \int_a^{T-S^+} Q_0(s,c)\, \tilde{Q}_2(T-s^+-s,c,b) \, \dd \mu_{T,s^+,c,b}(s)  \\
&\leq \sup_{s> a} \,Q_0(s,c)\, \tilde{Q}_2(T-s^+-a,c,b) .
\end{align*}
Lemma \ref{lem:Q0} tells us that $\displaystyle{\limsup_{s\rightarrow\infty} \frac{1}{s}\log \sup_{c\in[0,1) }Q_0(s,c)\leq -2\pi^2}$, and so by combining this with (\ref{eq:DeltaS}) we deduce that there is some $\Xi<\infty$ such that
\begin{align}
\frac{\WW(\EE_T^\bullet\cap\{S^->a\}\,|\,c,b,Y^+) }{\WW(\EE_T^\bullet\,|\,c,b,Y^+) } &\leq \frac{T}{k_2}\frac{\sup_{s>a} Q_0(s,c)}{\inf_{s\in[0,1]} Q_0(s,c)} \frac{\tilde{Q}_2(T-s^+-a,c,b)}{\tilde{Q}_2(T-s^+,c,b)} \label{eq:exptail} \\
&\leq \Xi \,\frac{\exp\{-(2\pi^2-1)a\}}{\exp\{-(\Gamma^\bullet+1)a\}} = \Xi \,\ee^{-\xi a} \nonumber
\end{align}
for all $c,b\in[0,1)$, $Y^+$, $T\geq 1$ and all $a\geq 0$.  Here $\xi = 2\pi^2-\Gamma^\bullet-2>0$.  (\ref{eq:tail_time}) now follows by integrating over $c$, $b$ and $Y^+$.
\ep

\bp[Proof of Lemma \ref{claim:claim4}]
Let $m_t = \inf\{W_s : 0\leq s\leq t\}$, $M_t = \sup\{W_s : 0\leq s\leq t\}$, and define the events
\begin{align*}
\EE_t^{\bullet -} &= \{m_t\geq -1 , L_x(t)\leq \tfrac{1}{2} \text{ for all }x\in[-1,1] \text{ and }L_x(t)\leq 1 \text{ for all }x\in\RR\} \\
\EE_t^{\bullet +} &= \{M_t - W_t\leq 1 , L_x(t)\leq \tfrac{1}{2} \text{ for all }x-W_t\in[-1, 1] \text{ and }L_x(t)\leq 1 \text{ for all }x\in\RR\} .
\end{align*}
Observe that if $(W_s)_{0\leq s\leq t}\in\EE_{t}^{\bullet +}$ and $(W_{t+s}-W_t)_{s\geq 0}\in\EE^{\bullet -}_{T-t}$ then $(W_s)_{s\geq 0}\in\EE_T^\bullet$.  Therefore since $(W_s)_{0\leq s\leq t}$ and $(W_{t+s}-W_t)_{s\geq 0}$ are independent we get
\begin{align*}
\WW(\EE_T^{\bullet}\,|\,\EE_t^{\bullet}) = \frac{\WW(\EE_T^{\bullet})}{\WW(\EE_t^{\bullet})}\geq \frac{\WW(\EE^{\bullet + }_t)\WW(\EE_{T-t}^{\bullet -})}{\WW(\EE_t^{\bullet})} ,
\end{align*}
and so if we can find a $k\geq 0$ such that
\begin{align}
\WW(\EE^{\bullet - }_t)= \WW(\EE_t^{\bullet +}) >k \,\WW(\EE^{\bullet}_t) \label{eq:EEpmc}
\end{align}
for all $t\geq 0$, then we would have shown (\ref{eq:EETeta}) with $\eta= k^2$.

From Lemma \ref{lem:below0} we know that when we condition on $\EE^\bullet_T$ then $S^-$ and (by symmetry) $S^+$ must both be (uniformly) small.  Therefore, since $S=T-s^--s^+$, there must be some universal constant $K_1$ such that
\begin{align*}
\WW\big(S=\textstyle{\int_0^{W_T}L_x(T)\,\dd x < T-K_1\,|\, \EE^\bullet_T, L_0(T)=c, L_{W_T}(T) = b}\big)\geq \displaystyle{\frac{1}{2}}
\end{align*}
for all $T\geq 1$ and $c,b\in [0,1)$.  From (\ref{eq:DeltaS}) we also know that there must be some universal $K_2$ such that
\begin{align*}
\sup_{T-K_1\leq s\leq T}\,{\tilde{Q}_2(c,b,s)}< {K_2}\,\tilde{Q}_2(c,b,T)
\end{align*}
for all $T\geq 1$ and $c,b\in[0,1)$.  Now write $f_T(\cdot,\cdot)$ for the join density of $L_0(T)$ and $L_{W_T}(T)$, and recall that $\tilde{Q}_2(c,b,s)$ represent the probability of $\{L_x(T)\leq 1$ for all $0\leq x\leq W_T\}$ conditionally on $L_0(T)=c$, $L_{W_T}(T)=b$ and $\int_0^{W_T}L_x(T)\,\dd x=s$.  We must then have
\begin{align}
\WW(\EE_T^\bullet)&\leq 2 K_2 \int^1_0 \int_0^1 \tilde{Q}_2(c,b,T) f_T(c,b)\,\dd c\,\dd b \leq 32 K_2 \int^{\frac{1}{4}}_0 \int_0^\frac{1}{4} \tilde{Q}_2(c,b,T) f_T(c,b)\,\dd c\,\dd b \label{eq:EETupper}.
\end{align}
Here the second inequality follows because both $\tilde{Q}_2(c,b,T)$ and $f_T(c,b)$ are decreasing as functions of $c$ and $b$.

On the other hand, since one can check that the proof of Lemma \ref{lem:below0} also holds when we condition on $\EE^{\bullet -}_T$ rather than $\EE^{\bullet}_T$, we see that there must also be a universal constant $K_3$ such that
\begin{align*}
\WW( S^-< K_3 \text{ and }S^+<K_3\,|\,\EE^{\bullet -}_T, L_0(T)=c, L_{W_T}(T)=b) \geq\frac{1}{2}
\end{align*}
for all $T\geq 1$, $c,b\in\left[0,\frac{1}{4}\right)$.  Therefore if we are given $c,b\in\left[0,\frac{1}{4}\right)$ we can lower bound $\WW(\EE^{\bullet -}_T\,|\,L_0(T)=c, L_{W_T}(T)=b)$ by
\begin{align*}
&\WW(\EE^{\bullet -}_T\,|\,L_0(T)=c, L_{W_T}(T)=b)\geq \frac{1}{2}\, \WW(\EE^{\bullet -}_T\,|\,L_0(T)=c, L_{W_T}(T)=b\text{ and }S^-,S^+<K_3)\\
&\quad\geq \frac{1}{2}\, \inf_{0\leq s\leq K_3}\mathbbm{Y}^0_c\big(Y_x\leq\tfrac{1}{2} \,|\,\textstyle{\int_0^\infty Y_x\,\dd x = s}\big) \times \displaystyle{\inf_{0\leq s\leq K_3}} \mathbbm{Y}^0_b\big(Y_x\leq 1 \,|\,\textstyle{\int_0^\infty Y_x\,\dd x s}\big) \times \tilde{Q}_2(c,b,T).
\end{align*}
We can now find constants $k_4,k_5>0$ such that 
\begin{align*}
\displaystyle{\inf_{c\in\left[0,\frac{1}{4}\right)}  \inf_{0\leq s\leq K_3}}\mathbbm{Y}^0_c\big(Y^-_x\leq \tfrac{1}{2}\text{ for all }x\geq 0 \,|\,\textstyle{\int_0^\infty Y_x\,\dd x =s}\big) &\geq k_4\\ \text{and}\quad \displaystyle{\inf_{c\in\left[0,\frac{1}{4}\right)}  \inf_{0\leq s\leq K_3}}\mathbbm{Y}^0_c\big(Y^-_x\leq 1 \text{ for all }x\geq 0 \,|\,\textstyle{\int_0^\infty Y_x\,\dd x =s}\big) &\geq k_5 ,
\end{align*}
Therefore by integrating and comparing with (\ref{eq:EETupper}) we get
\begin{align*}
\WW(\EE^{\bullet -}_T)&\geq \int^{\frac{1}{4}}_0 \int_0^\frac{1}{4} \WW(\EE^{\bullet -}_T\,|\,L_0(T)=c, L_{W_T}(T)=b) f_T(c,b)\,\dd c\,\dd b  \\
&\geq \frac{ k_4 k_5}{2} \int^{\frac{1}{4}}_0 \int_0^\frac{1}{4} \tilde{Q}_2(c,b,T) f_T(c,b)\,\dd c\,\dd b \geq \frac{k_4k_5}{64 K_2} \WW(\EE^\bullet_T) . %\qedhere
\end{align*}
Hence we have proved (\ref{eq:EETeta}) with $\eta = \displaystyle{\left(\frac{k_4k_5}{64 K_2}\right)^2}$.
\ep

\subsection{The weak convergence of $\WW(\cdot\,|\,\EE_T^{\bullet})$}\label{sec:conv}

To prove that $\WW(\cdot\,|\,\mathcal{E}_T^{\bullet})$ has a unique weak limit as $T\longrightarrow\infty$ we shall follow the techniques developed by Benjamini and Berestycki in\cite{b&b1}.  By doing so we shall in fact prove that the measures $\WW(\cdot\,|\,\mathcal{E}_T^{\bullet})$ form a Cauchy sequence with respect to the total variation distance on sets restricted to $\mathcal{F}_{R}$ (for each fixed $R>0$).  It turns out that this is a stronger condition than that of weak convergence, and so we start by presenting a lemma which shows that controlling the total variation distance (for each fixed $R>0$) does indeed imply weak convergence in the Skorokhod topology.  It will then suffice to prove that $\WW(\cdot\,|\,\mathcal{E}_T^{\bullet})$ is a Cauchy sequence in this sense.

Given $R> 0$ and probability measures $\pr$ and $\QQ$ on $(\Omega,\mathcal{F}_{R})$ we define the total variation metric, $d_R$, by
\begin{align*} d_{R}(\pr,\QQ) = \sup_{A\in\mathcal{F}_{R}}|\pr(A)-\QQ(A)| .
\end{align*}
We now have the following lemma from \cite[Lemma 6]{b&b1}.

\begin{lemma}\label{lem:weakconvergence}
Let $\{\pr_T\}_{T>0}$ be a sequence of probability measures on $\mathcal{F}$ which satisfy the following two conditions.
\vspace{-0.5\baselineskip}
\begin{description}
\item[Condition 1]
For every $R>0$ the restrictions of $\pr_T$ to $\mathcal{F}_{R}$ form a Cauchy sequence for the distance $d_R$.  I.e. for every $\varepsilon>0$ there exists $T_{R,\varepsilon}$ such that for all $T, {T^\prime} \geq T_{R,\varepsilon}$ we have $d_R(\pr_T,\pr_{{T^\prime}})<\varepsilon$.
\item[Condition 2] 
For each fixed $R>0$
\begin{align}\label{eq:escape}
\lim_{k\rightarrow\infty} \lim_{T\rightarrow\infty} \pr_T\left(\sup_{0\leq t\leq R}|W_t|\geq k\right) = 0.
\end{align}
\end{description}
Then there exists a unique probability measure $\pr$ such that $\pr_T\longrightarrow\pr$ weakly in the Skorokhod topology as $T\longrightarrow\infty$.
\end{lemma}

To prove that $\WW(\cdot\,|\,\mathcal{E}^{\bullet}_T)$ converges to a unique weak limit it now suffices to show that $\{\WW(\cdot\,|\,\mathcal{E}^{\bullet}_T)\}_{T>0}$ satisfies Condition 1 and Condition 2.  The second condition is required to show that $W_t$ does not escape to infinity in finite time -- and thus the limit of $\WW(\cdot\,|\,\mathcal{E}^{\bullet}_T)$ is non-trivial.  Since we already have the tools to show that $\WW(\cdot\,|\,\mathcal{E}^{\bullet}_T)$ satisfies Condition 2 then this is where we shall start.

\myparagraph{Proof that $\WW(\cdot\,|\,\mathcal{E}^{\bullet}_T)$ satisfies Condition 2}

Suppose that $R>0$ is fixed, then for each $\varepsilon>0$ we need to find an $M$ such that
\begin{align}
\WW\left(\sup_{0\leq t\leq R}|W_t|\geq M \,|\,\EE^{\bullet}_T\right) \leq \varepsilon \label{eq:noescape}
\end{align}
whenever $T\geq R$.
From Lemma \ref{claim:claim3} and Lemma \ref{claim:claim4} we see that there is a universal constant $\eta>0$ such that if $A$ is $\mathcal{F}_R$-measurable then
\begin{align}
\WW(A\,|\,\EE_T^\bullet) \leq \frac{1}{\eta}\,\WW(A\,|\,\EE_R^\bullet)
\end{align}
for all $T\geq R$.  $\displaystyle{\sup_{0\leq t\leq R}|W_t|}$ has a Guassian tail and so we can find an $M$ such that
\begin{align*}
\WW\left(\sup_{0\leq t\leq R}|W_t|\geq M\right) < \varepsilon\, \eta \,\WW(\EE^\bullet_R) .
\end{align*}
Hence
\begin{align*}
\WW\left(\sup_{0\leq t\leq R}|W_t|\geq M \,|\,\EE^{\bullet}_T\right) &\leq  \frac{1}{\eta}\, \WW\left(\sup_{0\leq t\leq R}|W_t|\geq M \,|\,\EE^{\bullet}_R\right) \leq \frac{1}{\eta} \frac{\WW\left(\sup |W_t|\geq M\right) }{\WW(\EE_R^\bullet)} <\varepsilon 
\end{align*}
for all $T\geq R$.
\hfill$\square$

\myparagraph{Proof that $\WW(\cdot\,|\,\mathcal{E}^{\bullet}_T)$ satisfies Condition 1}

Fix $R>0$ and $\varepsilon>0$.  Our aim is to find a constant $T_{R,\varepsilon}$ such that 
\begin{align}
|\WW(A\,|\,\EE^{\bullet}_T)-\WW(A\,|\,\EE^{\bullet}_{{T^\prime}})|<\varepsilon , \label{eq:Cauchy}
\end{align}
for every $A\in\mathcal{F}_R$ and all $T,{T^\prime}\geq T_{R,\varepsilon}$.

As a first step we shall decompose the event $A$ into the disjoint union of well behaved and badly behaved parts, $A = (A\cap B_{R,\varepsilon}) \sqcup (A\setminus B_{R,\varepsilon})$.  Provided we can show the probability of $B_{R,\varepsilon} \in\mathcal{F}_{R}$ is small, $\WW(B_{R,\varepsilon}\,|\,\mathcal{E}_{T}^{\bullet})<\tfrac{1}{2}\varepsilon$ for all $T\geq R$ say, then we can use the identity
\begin{align}
\WW(A\,|\,\mathcal{E}_{T}^{\bullet}) = \WW(A\,|\,\mathcal{E}_{{T^\prime}}^{\bullet}) \frac{\WW(\mathcal{E}_{{T^\prime}}^{\bullet}\,|\,A\cap\mathcal{E}_{T}^{\bullet})}{\WW(\mathcal{E}_{{T^\prime}}^{\bullet}\,|\,\mathcal{E}_{T}^{\bullet})} .
\end{align}
to show (\ref{eq:Cauchy}).  In particular if we have
\begin{align}
\left|1 - \frac{\WW(\mathcal{E}_{{T^\prime}}^{\bullet}\,|\,A\cap\mathcal{E}_{T}^{\bullet})}{\WW(\mathcal{E}_{{T^\prime}}^{\bullet}\,|\,\mathcal{E}_{T}^{\bullet})} \right|<\frac{1}{2}\,\varepsilon ,\label{eq:Cauchy2}
\end{align}
for all $A\in\mathcal{F}_R$ with $A\cap B_{R,\varepsilon}=\emptyset$ and all $T,{T^\prime}\geq T_{R,\varepsilon}$, then the result is proved.

\begin{claim}\label{claim:claim21}
Recall that we defined $m_R = \inf\{W_t : 0\leq t\leq R\}$ and $M_R = \sup\{W_t : 0\leq t\leq R\}$.  Suppose we set $\ell_{\max}(R)=\sup\{L_x(R):x\in\RR\}$, then for each $\varepsilon>0$ we can find $M>0$ and $k<1$ (depending on $R$ but not on $T$) for which the bad event
\begin{align}
B_{R,\varepsilon} = \{m_R\leq -M\}\cup\{M_R\geq M\}\cup\{\ell_{\max}(R)\geq k\}
\end{align}
has probability $\WW(B_{R,\varepsilon}\,|\,\EE_{T}^{\bullet})<\tfrac{1}{2}\varepsilon$ for all $T\geq R$.
\end{claim}
\bpc
From (\ref{eq:noescape}) we know that when $M$ is sufficiently large we have
\begin{align*}
\WW(\{m_R\leq -M\}\cup\{M_R\geq M\}\,|\,\EE_T^\bullet)  =\WW\left(\sup_{0\leq t\leq R}|W_t|\geq M\,|\,\EE_T^\bullet\right) <\frac{1}{4}\,\varepsilon
\end{align*}
for all $T\geq R$.  Having fixed such a value of $M$ we now observe that when we condition on $\EE_{R}^\bullet$ then then $\ell_{\max}(R)$ is a random variable supported on $[0,1)$.  From Lemma \ref{claim:claim3} and Lemma \ref{claim:claim4} of Section \ref{sec:ballistic} we have
\begin{align*}
\WW(\ell_{\max}(R)\geq y\,|\,\EE_{T}^{\bullet})&= \WW(\ell_{\max}(R)\geq y\,|\,\EE_{R}^{\bullet}) \frac{\WW(\EE_T^{\bullet}\,|\,\EE_{R}^{\bullet}\cap\{\ell_{\max}(R)\geq y\})}{\WW(\EE^{\bullet}_T\,|\,\EE_R^{\bullet})}\\&\leq\frac{1}{\eta} \,\WW(\ell_{\max}(R)\geq y\,|\,\EE_{R}^{\bullet}) ,
\end{align*}
where $\eta>0$ is independent of $T\geq R$.  Therefore there must exist $k\in[0,1)$ such that $\WW(\ell_{\max}(R)\geq k\,|\,\EE_{T}^{\bullet})\leq\frac{1}{4}\varepsilon$ for all $T\geq R$.  Summing these probabilities completes the proof.
\epc

By removing the behaviour of the bad event $B_{R,\varepsilon}$ we can ensure that the Brownian motion does not become trapped to the left of $M_R$, and therefore at large times the effect of the event $A$ becomes negligible.  We now aim to couple two processes $(W_t)_{t\geq 0}$ and $(\tilde{W}_t)_{t\geq 0}$ with respective laws $\WW(\cdot\,|\,A\cap\EE_{T}^{\bullet})$ and $\WW(\cdot\,|\, \EE_{T}^{\bullet})$ in such a way that their local times agree on a large region to the right of $M_R$.  More explicitly, suppose we are given two processes $(W_t)_{t\geq 0}$ and $(\tilde{W}_t)_{t\geq 0}$ with laws $\WW(\cdot\,|\,A\cap\EE^\bullet_T)$ and $\WW(\cdot\,|\,\EE^\bullet_T)$.  Let their respective local time profiles be $(L_x(t))_{x\in\RR,t\geq 0}$ and $(\tilde{L}_x(t))_{x\in\RR,t\geq 0}$, then we seek a pair of levels $X,\tilde{X}\in\RR$ such that
\begin{align}
&X>M_R = \sup\{W_t : 0\leq t\leq R\}, \quad\tilde{X}>0, \quad
L_X(T) = \tilde{L}_{\tilde{X}}(T) \nonumber \\&\qquad\text{and}\quad
\int_{-\infty}^{X}L_x(T)\,\dd x = \int_{-\infty}^{\tilde{X}} \tilde{L}_x(T)\,\dd x . \label{eq:XXconditions}
\end{align}
Our coupling will then replace $(\tilde{W}_t)_{t\geq 0}$ with a new process $(\hat{W}_t)_{t\geq 0}$ whose local time $(\hat{L}_x(t))_{x\in\RR ,t\geq 0}$ is equal in law to  $(\tilde{L}_x(t))_{x\in\RR ,t\geq 0}$ and is such that
\begin{align}
\hat{L}_x(T) = \left\{ \begin{array}{ll} \tilde{L}_x(T) &\text{for all }x\leq \tilde{X} \\
L_{X-\tilde{X}+x}(T) &\text{for all }x\geq \tilde{X} \end{array}\right. \label{eq:coupling_local} .
\end{align}

\begin{definition}[Coupling]\label{def:jointlaw}
Fix $A\in\mathcal{F}_R$ with $A\cap B_{R,\varepsilon}=\emptyset$, and let $T\geq R$.  Suppose $(W_t)_{t\geq 0}$ is a process with law $\WW(\cdot\,|\,A\cap\EE_{T}^{\bullet})$ and $(\tilde{W}_t)_{t\geq 0}$ is an independent process with law $\WW(\cdot\,|\, \EE_{T}^{\bullet})$.  Define the random variables $X$, $\tilde{X}$ and $\Delta_T$ by
\begin{align*}
X &= \inf\big\{z>M_R : \exists\, \tilde{z}>0\text{ with } L_z(T) = L_{\tilde{z}}(T)\text{ and } \textstyle{\int_{-\infty}^{z}{L}_{x}(T)\,\dd x = \int_{-\infty}^{\tilde{z}}\tilde{L}_{x}(T)\,\dd x} \big\} \\
\tilde{X} &= \inf\big\{ \tilde{z}>0 :  \textstyle{\int_{-\infty}^{\tilde{z}}\tilde{L}_{x}(T)\,\dd x = \int_{-\infty}^{X}L_{x}(T)\,\dd x}\big\}\\
\Delta_T &= \left\{\begin{array}{ll} W_T - X&\text{if }W_T> X \text{ and }\tilde{W}_T>\tilde{X} \\ -\infty &\text{otherwise}\end{array} \right. .
\end{align*}
We now define a new process $(\hat{W}_t)_{t\geq 0}$ by considering two cases.
\vspace{-0.5\baselineskip}\begin{enumerate}
\item If $X\geq W_T$ or $\tilde{X}\geq \tilde{W}_T$ then set $(\hat{W}_t)_{t\geq 0}=(\tilde{W}_t)_{t\geq 0}$.
\item Otherwise we have $X< W_T$ and $\tilde{X}< \tilde{W}_T$.  It is known from It\^o's theory of excursions that a Brownian motion started at $0$ can be decomposed into a sequence of positive and negative excursions from $0$, and that the excursions are indexed by the local time at level $0$.  For a reference to excursion theory see \cite[Chapter XII]{r&y}.  Therefore, because the local times at level $0$ of $(W_{t} - X)_{0\leq t\leq T}$ and $(\tilde{W}_{ t }- \tilde{X})_{0\leq t\leq T}$ are equal at time $T$, then we can form $(\hat{W}_{t}-\tilde{X})_{0\leq t\leq T}$ by taking the negative excursions of $(\tilde{W}_{ t} - \tilde{X})_{0\leq t\leq T}$ and combining with the positive excursions of $(W_{t} - X)_{0\leq t\leq T}$.  By doing this we ensure that $\hat{L}_x(T) = L_x(T)$ for all $x\leq \tilde{X}$ and $\hat{L}_x(T) = L_{X-\tilde{X}+x}(T)$ for all $x\geq \tilde{X}$, and so (\ref{eq:coupling_local}) is satisfied.
\end{enumerate}
Write $\mathbbm{J}_{T,A}$ for the joint law of $(W_t)_{0\leq t\leq T}$ and $(\hat{W}_t)_{0\leq t\leq T}$.  It is clear from our construction that the first marginal satisfies
\begin{align}\mathbbm{J}_{T,A}^{(1)}\big((L_x(T))_{x\in\RR}\in \cdot\big)=\WW\big((L_x(T))_{x\in\RR}\in \cdot\,|\,A\cap\EE^\bullet_T\big). \label{eq:marginal1}
\end{align}
Our construction also gives
\begin{align*}
\mathbbm{J}_{T,A}^{(2)}\big((\hat{L}_{\tilde{X}-x}(T))_{x\geq 0}\in \cdot\big)=\WW\big((\tilde{L}_{\tilde{X}-x}(T))_{x\geq 0}\in \cdot\,|\,\EE^\bullet_T\big).
\end{align*}
Therefore we now fix $(\hat{L}_{\tilde{X}-x}(T))_{x\geq 0}=(Y_x)_{x\geq 0}$ and consider
\begin{align*}
&\mathbbm{J}_{T,A}^{(2)}\big((\hat{L}_{\tilde{X}+x}(T))_{x\geq 0}\in \cdot\,|\, (\hat{L}_{\tilde{X}-x}(T))_{x\geq 0}=(Y_x)_{x\geq 0}\big) \\
&\qquad = \WW\big((L_{X+x}(T))_{x\in\RR}\in \cdot\,|\,A\cap\EE^\bullet_T, ({L}_{X-x}(T))_{x\geq 0}=(Y_x)_{x\geq 0} \big).
\end{align*}
From Theorem \ref{thm:ray_knight_fixed_t}, and the fact that a square Bessel process is Markovian, we see that when we condition on $W_T > X$ then the distribution of $(L_{X+x}(T))_{x\geq 0}$ depends only on $L_X(T)$ and $\int^X_{-\infty} L_x(T)\,\dd x$.  We also observe that because  $A$ is $\mathcal{F}_R$ measurable and $M_R<X$ then, as far as the distribution of $(L_{X+x}(T))_{x\in\RR}$ is concerned, conditioning on $\{A\cap\EE_T^\bullet, ({L}_{X-x}(T))_{x\geq 0}=(Y_x)_{x\geq 0}  \}$ is equivalent to conditioning on $\{\EE_T^\bullet, ({L}_{X-x}(T))_{x\geq 0}=(Y_x)_{x\geq 0} \}$.
Hence if we know that $({L}_{X-x}(T))_{x\geq 0}$ has $L_X(T)=y$ and $\int^X_{-\infty} L_x(T)\,\dd x = r$ then
\begin{align*}
&\mathbbm{J}_{T,A}^{(2)}\big((\hat{L}_{\tilde{X}+x}(T))_{x\geq 0}\in \cdot\,|\, (\hat{L}_{\tilde{X}-x}(T))_{x\geq 0}=(Y_x)_{x\geq 0}\big) \\
&\qquad = \WW\big((L_{X+x}(T))_{x\in\RR}\in \cdot\,|\,A\cap\EE^\bullet_T, L_X(T)=y, \textstyle{\int^X_{-\infty} L_x(T)\,\dd x = r} \big) \\
&\qquad = \WW\big((\tilde{L}_{\tilde{X}+x}(T))_{x\in\RR}\in \cdot\,|\,\EE^\bullet_T, \tilde{L}_{\tilde{X}}(T)=y, \textstyle{\int^{\tilde{X}}_{-\infty} \tilde{L}_x(T)\,\dd x = r} \big) \\
&\qquad = \WW\big((\tilde{L}_{\tilde{X}+x}(T))_{x\in\RR}\in \cdot\,|\,\EE^\bullet_T, (\tilde{L}_{\tilde{X}-x}(T))_{x\geq 0}=(Y_x)_{x\geq 0} \big).
\end{align*}
Since this holds for all choices of $(Y_x)_{x\geq 0}$ then by interating we can conclude that \begin{align}
\mathbbm{J}_{T,A}^{(2)}\big((\tilde{L}_x(T))_{x\in\RR}\in \cdot\big)=\WW\big((\tilde{L}_x(T))_{x\in\RR}\in \cdot\,|\,\EE^\bullet_T\big), \label{eq:marginal2}
\end{align}
as desired.
\end{definition}

Having constructed this coupling our next step is to show that $\Delta_T$ is large with high $\mathbbm{J}_{T,A}$-probability.  This would mean that there is a large region of space on which the local time profiles of $(W_t)_{0\leq t\leq T}$ and $(\hat{W}_t)_{0\leq t\leq T}$ agree.  Later we shall use this property to show that the ratio between $\WW(\EE_{{T^\prime}}^\bullet\,|\,(W_t)_{0\leq t\leq T})$ and $\WW(\EE_{{T^\prime}}^\bullet\,|\,(\hat{W}_t)_{0\leq t\leq T})$ becomes arbitrarily close to 1 as $T,{T^\prime}\longrightarrow\infty$.

\begin{claim}\label{claim:claim22}
Given $A\in\mathcal{F}_R$ with $A\cap B_{R,\varepsilon}=\emptyset$, let $\{\mathbbm{J}_{T,A}\}_{T\geq R}$ be the sequence of joint laws defined by Definition \ref{def:jointlaw}.  We then have that
\begin{align}
\mathbbm{J}_{T,A}\left(\Delta_T \geq \tfrac{1}{2} T \right) \longrightarrow 1 \quad \text{as}\quad T\longrightarrow \infty . \label{eq:Deltaconv}
\end{align}
What is more, this convergence is uniform over all events $A\in\mathcal{F}_R$ with $A\cap B_{R,\varepsilon}=\emptyset$.
\end{claim}
\bpc
Given $(W_t)_{t\geq 0}\in A\cap\EE_T^\bullet$ and $(\tilde{W}_t)_{t\geq 0}\in\EE_T^\bullet$ we define 
\begin{align*}
\Sigma_{R,T} = \max \left\{\textstyle{\int_{-\infty}^{M_R}L_x(T)\,\dd x , \int_{-\infty}^{0}\tilde{L}_x(T)\,\dd x}\right\} .
\end{align*}
Now let $\delta>0$.  Our first task is to show that we can always find a $K_1$ (depending on $R$, but not on $T$ or $A$ such that
\begin{align}
\mathbbm{J}_{T,A}(\Sigma_{R,T} \geq K_1)< \frac{1}{3}\,\delta . \label{eq:sigmaest}
\end{align}
Let $M$ be as given by Claim \ref{claim:claim21}, and define $r = \int_{-\infty}^{-M} L_x(T)\,\dd x$.  Because $A\cap B_{R,\varepsilon} = \emptyset$ then we must have $M_R \leq M$.  Therefore as  $L_x(T)\leq 1$ for all $x\in\RR$ it must follow that $\textstyle{\int_{-\infty}^{M_R}L_x(T)\,\dd x \leq 2 M +r}$.  Now suppose we condition on $(W_t)_{0\leq t\leq R} = (V_t)_{0\leq t\leq R}$ and $(L_x(T))_{-M\leq x\leq M} = (Y_x)_{-M\leq x\leq M}$.  If we assume that $\displaystyle{\sup_{0\leq t\leq R}|V_t|\leq M}$  then we can see from Theorem \ref{thm:ray_knight_fixed_t} that the law of $(L_x(T))_{x\in\RR\setminus [-M,M]} $ depends only on the values of $Y_{-M}$, $Y_M$ and $\int_{-M}^M Y_x\,\dd x$ (and not on $(V_t)_{0\leq t\leq R}$).  Therefore, if we also condition on $(L_{W_T+x}(T))_{x\geq 0} = (Y^+_x)_{x\geq 0}$ with $\int_0^\infty Y^+_x\,\dd x = s^+ $, then we can follow the argument of Lemma \ref{lem:below0} and replace (\ref{eq:exptail}) by
\begin{align*}
\frac{\WW(\{r>a\}\cap\EE_{T}^{\bullet}\,|\, V,Y,Y^+)}{\WW(\EE_{T}^{\bullet}\,|\,  V,Y,Y^+)} 
&\leq \frac{\WW(\{r>a\}\cap\{L_x(T)\leq 1 \text{ for all }x\in\RR\setminus[-M,M]\}\,|\, V,Y,Y^+)}{\WW(\EE_{T}^{\bullet}\,|\,  V,Y,Y^+)} \\
&\leq \frac{T}{k_2} \sup_{s> a} \sup_{c\in[0,1)}\frac{\sup_{s> a} Q_0(c,s)}{\inf_{s\in[0,1]}Q_0(c,s)}\sup_{c,b\in[0,1)}\frac{\tilde{Q}_2(T-2M-s^+ -a,c,b)}{\tilde{Q}_2(T-s^+,c,b)} .
\end{align*}
Since our assumption on the event $A$ implies that $\displaystyle{\sup_{0\leq t\leq R}|W_t|\leq M}$ we can now integrate over $A$ and deduce that there exists a universal constant $\Xi\leq K_2<\infty$ such that
\begin{align*}
\WW(\{r>a\}\cap\EE_T^{\bullet}\,|\,A)\leq K_2 \,\ee^{-\xi a}\,\WW(\EE_T^\bullet), 
\end{align*} 
for all $A\in\mathcal{F}_R$ with $A\cap B_{R,\varepsilon}=\emptyset$, all $T\geq R$ and all $a\geq 0$.  Here $\xi>0$ is as given by Lemma \ref{lem:below0}.  On the other hand, our assumption that $A\cap B_{R,\varepsilon}=\emptyset$ means that $L_x(R)\leq k <1$ for all $x\in\RR$.  Therefore if we let $\tau_{M+1}$ be the first hitting time of $M+1$ (and assume $k\geq\tfrac{1}{2}$), then the intersection of the events 
\begin{align*}
&A\cap \{L_x(\tau_{M+1})-L_x(R)\leq 1-k\,\forall\,x\in\RR\} 
\cap\{L_x(T)-L_x(\tau_{M+1})\leq 1\,\forall\,x\geq M+1,\\
&\qquad  L_x(T)-L_x(\tau_{M+1})\leq \tfrac{1}{2} \,\forall\,M<x<M+1 \text{ and }L_x(T)-L_x(\tau_{M+1})=0\text{ otherwise}\}
\end{align*}
is contained within $\EE^{\bullet}_T$.
%Observe that these three events are measurable with respect to $\mathcal{F}_R$, $\mathcal{F}_{\tau_{M+1}}\setminus\mathcal{F}_R$ and $\mathcal{F}_T\setminus\mathcal{F}_{\tau_{M+1}}$ and therefore they are independent.
The probability of $\{L_x(\tau_{M+1})-L_x(R)\leq 1-k$ for all $x\in\RR\} $ is minimised when $W_R$ is minimal.  Therefore, since our conditions on $A$ imply that $W_R\geq -M$, we find that there must be a positive quantity $\zeta(M,k)>0$ such that
\begin{align*}
\WW(L_x(\tau_{M+1})-L_x(R)\leq 1-k\text{ for all }x\in\RR\,|\,A)&\\
\geq \WW(L_x(\tau_{2M+1}) \leq 1-k\text{ for all }x\in\RR\,|\,A)&\geq\zeta(M,k)>0.
\end{align*}
Note that since $M$ and $k$ depend only on $R$ and $\varepsilon$ then so too does $\zeta(M,k)$.  The probability of the third event is equal to $\WW(\EE^{\bullet -}_{T-\tau_{M+1}})$, where $\EE^{\bullet -}_{T-\tau_{M+1}}$ is defined in the proof of Lemma \ref{claim:claim4}.  As we know that $\WW(\EE^{\bullet -}_{T-\tau_{M+1}})\geq \WW(\EE^{\bullet -}_{T})\geq \eta\, \WW(\EE_T^{\bullet})$ then we get 
\begin{align*}
\WW(\{r > a\}\,|\,A\cap\EE_T^\bullet) &=\frac{\WW(\{r>a\}\cap\EE_{T}^{\bullet})\,|\,A)}{\WW(\EE_T^{\bullet}\,|\,A)} \leq \frac{K_2\,\ee^{-\xi a}\,\WW(\EE_T^\bullet)}{\zeta(M,k)\,\eta\,\WW(\EE_T^\bullet)} \leq K_3 \,\ee^{-\xi a},
\end{align*}
for some $K_3<\infty$ which depends only on $R$ and $\varepsilon$.  From this we conclude that there must exist a $K_4<\infty$ such that
\begin{align*}
\mathbbm{J}^{(1)}_{T,A} \left(\textstyle{\int_{-\infty}^{M_R}{L}_x(T)\,\dd x \geq K_4}\right) < \frac{1}{6}\,\delta ,
\end{align*}
For all $A\in\mathcal{F}_R$ with $A\cap B_{R,\varepsilon}=\emptyset$ and all $T\geq R$.  Lemma \ref{lem:below0} also tells us that we can find a $K_5$ such that
\begin{align*}
\mathbbm{J}^{(2)}_{T,A} \left(\textstyle{\int_{-\infty}^{0}\tilde{L}_x(T)\,\dd x \geq K_5}\right) < \frac{1}{6}\,\delta ,
\end{align*}
for all $A$ and $T\geq R$.  Thus we can conclude that (\ref{eq:sigmaest}) holds for $K_1 = \max\{K_4,K_5\}$.

%%%%%

Using symmetry we can see that $S^+ = \int_{W_T}^\infty L_x(T)\,\dd x$ is equal in law to $S^- = \int_{-\infty}^0 L_x(T)\,\dd x$.  Therefore it is a simple corollary of Lemma \ref{lem:below0} to find a $K_6<\infty$ such that
\begin{align*}
\mathbbm{J}^{(1)}_{T,A} (S^+ >K_6)<\frac{1}{6}\,\delta \quad \text{and} \quad \mathbbm{J}^{(2)}_{T,A} (\tilde{S}^+ >K_6)<\frac{1}{6}\,\delta
\end{align*}
for all $A$ and $T\geq R$.  By combining this with \ref{eq:sigmaest}) we see that
\begin{align}
\mathbbm{J}_{T,A}\left(\textstyle{\Sigma_{R,T} < \frac{1}{6} T \text{ and } \min\left\{ \int_{-\infty}^{W_T}L_x(T)\,\dd x , \int_{-\infty}^{\tilde{W}_T}\tilde{L}_x(T)\,\dd x\right\} >\frac{5}{6} T} \right)\leq \frac{2}{3}\,\delta \label{eq:XXconditions1}
\end{align}
for all $A$ and all $T\geq\max\{R, 6K_1,6K_6\}$.

%%%%

Recall that in Definition \ref{def:auxiliary} we constructed $(Z_t)_{t\geq 0}$ as an axillary process to $(L_x(T))_{x\geq 0}$ and that $(Z_t)_{t\geq 0}$ satisfies $Z_t = L_{\rho(t)}(T)$, where $\rho(t)$ is such that $\int_0^{\rho(t)} L_x(T)\,\dd x = t$.  Let $(\tilde{Z}_t)_{t\geq 0}$ be the axillary process to $(\tilde{L}_x(T))_{x\geq 0}$, $S^- = \int_{-\infty}^{0}L_x(T)\,\dd x$ and $\tilde{S}^- = \int_{-\infty}^{0}\tilde{L}_x(T)\,\dd x$.  Suppose we can find a 
$t\geq \Sigma_{R,T}$
%$t\geq \max\{\Sigma_{R,T} - S^-, \Sigma_{R,T} - \tilde{S}^-\}$
such that $Z_{t-S^-}=\tilde{Z}_{t-\tilde{S}^-}$.  If we put $X=\rho(t-S^-)$ and $\tilde{X} = \tilde{\rho}(t-\tilde{S}^-)$ then we can check that
$X>M_R$, $\tilde{X}>0$, $L_X(T) = \tilde{L}_{\tilde{X}}(T)$ and $\int_{-\infty}^{X}L_x(T)\,\dd x = \int_{-\infty}^{\tilde{X}} \tilde{L}_x(T)\,\dd x$.  Thus $X$ and $\tilde{X}$ satisfy (\ref{eq:XXconditions}).
%Moreover if $t\leq \min\{T-S^+-S^-,T-\tilde{S}^+-\tilde{S}^-\}$ then we also have $X$

If we assume that $\Sigma_{R,T} < \frac{1}{6} T $ and $\min\left\{ \int_{-\infty}^{W_T}L_x(T)\,\dd x , \int_{-\infty}^{\tilde{W}_T}\tilde{L}_x(T)\,\dd x\right\} >\frac{5}{6} T$, then $(Z_{t-S^-})_{\Sigma_{R,T}\leq t\leq \frac{7}{8}T}$ and $(\tilde{Z}_{t-\tilde{S}^-})_{\Sigma_{R,T}\leq t\leq \frac{5}{6}T}$ are both Markov processes which are conditioned to stay in $[0,1]$, and with infinitesimal generator given by (\ref{eq:aux_gen}).   For a fixed $\Sigma_{R,T}\leq t\leq \frac{5}{6}T-1$ consider $\mathbbm{J}_{T,A}(Z_{t-S^{-}+s}=\tilde{Z}_{t-{\tilde{S}^-}+s}$ for some $s\in[0,1])$.  This probability will depend on $Z_{t-S^-}$ and $\tilde{Z}_{t-\tilde{S}^-}$, but must be positive for all $(Z_{t-S^-},\tilde{Z}_{t-\tilde{S}^-})\in[0,1]^2$.  Therefore, by compactness, we find a $k_7>0$ such that
\begin{align*}
\mathbbm{J}_{T,A}(Z_{t+s}=\tilde{Z}_{t+s}\text{ for some } s\in[0,1])\geq k_7
\end{align*}
for all possible $t$, $Z_{t-S^-}$ and $\tilde{Z}_{t-\tilde{S}^-}$.   If $K_7$ is large enough then $(1-k_7)^{K_7} < \frac{1}{3} \delta$.  Therefore we can use that that $Z$ and $\tilde{Z}$ are Markovian to deduce that
\begin{align}
\mathbbm{J}_{T,A} \left(\text{there exists } \Sigma_{R,T}< t<\tfrac{1}{3}T \text{ such that } Z_{t-S^-}=\tilde{Z}_{t-\tilde{S}^-}\right) > 1 - \frac{1}{3}\,\delta , \label{eq:XXconditions2}
\end{align}
for all $T\geq 6K_7$.  If we assume that the events given by (\ref{eq:XXconditions1}) and (\ref{eq:XXconditions2}) both hold then we must have $\tilde{W}_T>\tilde{X}$ and $\int^{W_T}_{X} L_x(T)\,\dd x \geq \frac{1}{2}T$.  Since $L_x(T)\leq 1$ for all $x\in\RR$ we must then also have $\Delta_T = W_T - X > \frac{1}{2}T$ and so 
\begin{align*}
\mathbbm{J}_{T,A}\left(\Delta_T\geq \tfrac{1}{2}T\right) > 1-\delta
\end{align*}
for all $A\in\mathcal{F}_R$ with $A\cap B_{R,\varepsilon}=\emptyset$ and all $T\geq\max\{R, 6K_1,6K_6,6K_7\}$.  As $\delta$ was arbitrary then (\ref{eq:Deltaconv}) is now proved.
\epc

We now have the tools we need to estimate (\ref{eq:Cauchy2}).
%Recall that the marginals of $\mathbbm{J}_{T,A}$ satisfy $\mathbbm{J}^{(1)}_{T,A} = \WW(\cdot\,|\,A\cap\EE_T^{\bullet})$ and $\mathbbm{J}^{(2)}_{T,A} = \WW(\cdot\,|\,\EE_T^{\bullet})$, and that the way we have constructed $\mathbbm{J}_{T,A}$ ensures that the distribution of $(L_{W_T+x}(T))_{x\geq-\Delta_T}$ with respect to $\mathbbm{J}_{T,A}^{(1)}$ is equal to the distribution of $(\tilde{L}_{\tilde{W}_T+x}(T))_{x\geq-\Delta_T}$ with respect to $\mathbbm{J}_{T,A}^{(2)}$.
Given $T\leq {T^\prime}$, define $m_{T,{T^\prime}}=\inf\{W_t-W_T:T\leq t\leq {T^\prime}\}$, and observe that we can partition the event $\EE^{\bullet}_{{T^\prime}}$ as $\EE^{\bullet}_{{T^\prime}} = (\EE^{\bullet}_{{T^\prime}} \cap \{m_{T,{T^\prime}} \geq - \Delta_T\}) \sqcup (\EE^{\bullet}_{{T^\prime}} \cap \{m_{T,{T^\prime}} < - \Delta_T\})$.  Therefore (\ref{eq:Cauchy2}) becomes
\begin{align}
&\left|1 - \frac{\WW(\mathcal{E}_{{T^\prime}}^{\bullet}\,|\,A\cap\mathcal{E}_{T}^{\bullet})}{\WW(\mathcal{E}_{{T^\prime}}^{\bullet}\,|\,\mathcal{E}_{T}^{\bullet})} \right|= \left|1 - \frac{\mathbbm{J}_{T,A}^{(1)}(\EE^{\bullet}_{{T^\prime}})}{\mathbbm{J}_{T,A}^{(2)}(\EE^{\bullet}_{{T^\prime}}) } \right|\nonumber \\
&\qquad= \left|1 - \frac{\mathbbm{J}_{T,A}^{(1)}(\EE^{\bullet}_{{T^\prime}}\cap\{ m_{T,{T^\prime}} \geq - \Delta_T\}) +  \mathbbm{J}_{T,A}^{(1)}(\EE^{\bullet}_{{T^\prime}}\cap\{ m_{T,{T^\prime}} < - \Delta_T\})}{\mathbbm{J}_{T,A}^{(2)}(\EE^{\bullet}_{{T^\prime}}\cap\{ m_{T,{T^\prime}} \geq - \Delta_T\}) +  \mathbbm{J}_{T,A}^{(2)}(\EE^{\bullet}_{{T^\prime}} \cap\{ m_{T,{T^\prime}} < - \Delta_T\}) } \right| \label{eq:Cauchy3}.
\end{align}
We now claim that $\mathbbm{J}_{T,A}^{(1)}(\EE^{\bullet}_{{T^\prime}}\cap\{ m_{T,{T^\prime}} \geq - \Delta_T\})$ and $\mathbbm{J}_{T,A}^{(2)}(\EE^{\bullet}_{{T^\prime}}\cap\{ m_{T,{T^\prime}} \geq - \Delta_T\})$ are equal.  Indeed, suppose we are given $\Delta_T>0$ and $(W_t)_{0\leq t\leq T}\in A\cap\EE^\bullet_T$ (or $(\hat{W}_t)_{0\leq t\leq T}\EE^\bullet_T$).  The probability of the event $\EE^{\bullet}_{{T^\prime}}\cap\{ m_{T,{T^\prime}} \geq - \Delta_T\}$ is then exactly equal to the probability that an independent Brownian motion $(\bar{W}_t)_{0\leq t\leq {T^\prime}-T}$ has $\bar{L}_x({T^\prime}-T)\leq 1-L_{W_t + x}(T)$ (or $\bar{L}_x({T^\prime}-T)\leq 1-\hat{L}_{\hat{W}_t + x}(T)$) for all $x>-\Delta_T$ and $\bar{L}_x({T^\prime}-T)=0$ otherwise.  Because of the way that $\mathbbm{J}_{T,A}$ was constructed in Definition \ref{def:jointlaw} we know that the distribution of $(L_{W_T+x}(T))_{x\geq-\Delta_T}$ with respect to $\mathbbm{J}_{T,A}^{(1)}$ is equal to the distribution of $(\hat{L}_{\hat{W}_T+x}(T))_{x\geq-\Delta_T}$ with respect to $\mathbbm{J}_{T,A}^{(2)}$.  Thus the equality of $\mathbbm{J}_{T,A}^{(1)}(\EE^{\bullet}_{{T^\prime}}\cap\{ m_{T,{T^\prime}} < \Delta_T\})$ and $\mathbbm{J}_{T,A}^{(2)}(\EE^{\bullet}_{{T^\prime}}\cap\{ m_{T,{T^\prime}} < \Delta_T\})$ must hold.
From equation (\ref{eq:Cauchy3}) we now get
\begin{align*}
\left|1 - \frac{\WW(\mathcal{E}_{{T^\prime}}^{\bullet}\,|\,A\cap\mathcal{E}_{T}^{\bullet})}{\WW(\mathcal{E}_{{T^\prime}}^{\bullet}\,|\,\mathcal{E}_{T}^{\bullet})} \right|\leq 
&\left| \frac{\mathbbm{J}_{T,A}^{(1)}(\EE^{\bullet}_{{T^\prime}}\cap\{ m_{T,{T^\prime}} < - \Delta_T\}) }{  \mathbbm{J}_{T,A}^{(1)}(\EE^{\bullet}_{{T^\prime}}\cap\{ m_{T,{T^\prime}} \geq - \Delta_T\})}-\frac{\mathbbm{J}_{T,A}^{(2)}(\EE^{\bullet}_{{T^\prime}}\cap\{ m_{T,{T^\prime}} < - \Delta_T\}) }{  \mathbbm{J}_{T,A}^{(2)}(\EE^{\bullet}_{{T^\prime}} \cap\{ m_{T,{T^\prime}} \geq - \Delta_T\}) } \right| \\
\leq & \left| \mathbbm{J}_{T,A}^{(1)}( m_{T,{T^\prime}}< - \Delta_T\,|\,\EE_{{T^\prime}}^{\bullet}) - \mathbbm{J}_{T,A}^{(2)}( m_{T,{T^\prime}}< - \Delta_T\,|\,\EE_{{T^\prime}}^{\bullet})\right| .
\end{align*} 
%When we condition on $\EE^{\bullet}_{T}$ then $-m_T$ has an exponential tail (uniformly in $T$).  By applying Lemma \ref{lem:b&b} we find that the law of $m_{T,{T^\prime}}$ given $\EE_{{T^\prime}}^{\bullet}$ is stochastically dominated by that of $-m_{{T^\prime}-T}$, and so $m_{T,{T^\prime}}$ must also have a (uniform) exponential tail.
By using Lemma \ref{lem:below0} it can be shown that the tail of $m_{T,{T^\prime}}$ can be uniformly bounded over all $R\leq T\leq {T^\prime}$.  Since Claim \ref{claim:claim22} tells us that $\Delta_T$ becomes arbitrarily large as $T\longrightarrow\infty$ we deduce that $\mathbbm{J}_{T,A}^{(1)}( m_{T,{T^\prime}}< - \Delta_T\,|\,\EE_{{T^\prime}}^{\bullet})$ and $\mathbbm{J}_{T,A}^{(2)}( m_{T,{T^\prime}}< -\Delta_T\,|\,\EE_{{T^\prime}}^{\bullet})$ must converge to 0 as $T\longrightarrow\infty$.  Therefore we can find a $T_{R,\varepsilon}$ such that $\mathbbm{J}_{T,A}^{(1)}( m_{T,{T^\prime}}< - \Delta_T\,|\,\EE_{{T^\prime}}^{\bullet})<\frac{1}{4}\varepsilon$ and $\mathbbm{J}_{T,A}^{(2)}( m_{T,{T^\prime}}< - \Delta_T\,|\,\EE_{{T^\prime}}^{\bullet})<\frac{1}{4}\varepsilon$ for all $A\in\mathcal{F}_R$ and all $T_{R,\varepsilon}\leq T\leq {T^\prime}$, and so  (\ref{eq:Cauchy2}) is satisfied.
\hfill$\square$

Since both Condition 1 and Condition 2 are satisfied then Lemma \ref{lem:weakconvergence} tells us that $\WW(\cdot\,|\,\EE_{T}^{\bullet})$ must converge weakly as $T\longrightarrow\infty$, and so Theorem \ref{thm:reduced_speed} is proved.
\hfill$\square$
%\nl\vspace{-2\baselineskip}

\section{Proof of Theorem \ref{thm:epsilon3a}}\label{sec:proof2}

In Benjamini and Berestycki's paper, \cite{b&b1}, it is shown that when we condition on $\EE_{a}^*$ then $(L_x(\tau_a))_{x\geq 0}$ converges to a process $(L_x(\infty))_{x\geq 0}$ with stationary distribution $\mu^*$ as $a\longrightarrow\infty$.  As a first step to proving Theorem \ref{thm:epsilon3a} we shall now show that $(L_x(\infty))_{x\geq 0}$ also has a stationary distribution, $\mu^{\bullet}$, with respect to $\QQ^{\bullet}=\displaystyle{\lim_{T\rightarrow \infty}\WW(\cdot\,|\,\EE^\bullet_T})$.  Having done this we then complete the proof using calculations of Radon--Nikodym derivatives.

\subsection{The stationary distribution of $(L_x(\infty))_{x\geq 0}$}

As a consequence of Lemma \ref{lem:mu_bullet} we know that there is a unique measure $\mu^{\bullet}\in E_C=\{\mu\in\mathcal{P}(\RR):\mathrm{support}(\mu)\subseteq[0,1]\}$ which minimises $\ex(\mu)^{-1} I_2(\mu)$ over all $\mu\in E_C$.  By using the Donsker--Varadhan Theorem we shall now show that $\mu^{\bullet}$ gives the stationary distribution of $(L_x(\infty))_{x\geq 0}=\displaystyle{\lim_{T\rightarrow\infty}}(L_x(T))_{x\geq 0}$ with respect to $\QQ^{\bullet}$.

\begin{lemma}\label{lem:Yconv}
Let $\mu^{\bullet}$ be as defined above, suppose $(Y_x)_{x\geq 0}$ is a $\BESQ{2}{y}$ process for some $c\in[0,1)$, and recall that we defined $\rho(s) = \inf\{u:\int_0^u Y_x\,\dd x\geq s\}$.  If we now condition on the event $\{Y_x\leq 1$ for all $0\leq x\leq\rho(s)\}$ then $L((Y_x),\rho(s),\cdot)$ converges in $\mathbbm{Y}_c$-probability to $\mu^{\bullet}$ as $s\longrightarrow\infty$.
\bp
Let $U\subseteq\mathcal{P}(\RR)$ be any open set (with respect to the weak topology) containing $\mu^{\bullet}$.  It now suffices to prove that
\begin{align}
&\liminf_{s\rightarrow\infty} \frac{1}{s}\log \mathbbm{Y}_c(L((Y_x),\rho(c),\cdot)\in U\cap E_C)\nonumber \\
&\qquad\qquad - \limsup_{s\rightarrow\infty} \frac{1}{s}\log \mathbbm{Y}_c(L((Y_x),\rho(s),\cdot)\in U^{\textsc{c}}\cap E_C) > 0 .\label{eq:mubulletU}
\end{align}
This would imply that the ratio between $\mathbbm{Y}(L((Y_x),\rho(s),\cdot)\in U\,|\,Y_x\leq 1$ for all $ x\in[0,\rho(s)])$ and $\mathbbm{Y}(L((Y_x),\rho(s),\cdot)\in U^{\textsc{c}}\,|\,Y_x\leq 1 $ for all $ x\in[0,\rho(s)])$ tends to infinity as $s\longrightarrow \infty$.  Since this holds for arbitrary $U\ni \mu^{\bullet}$ then $L((Y_x),\rho(s),\cdot)$ must converge to $\mu^{\bullet}$ in $\QQ^\bullet$-probability.

Define $E_O=\{\mu\in\mathcal{P}(\RR):\mathrm{support}(\mu)\subseteq [0,1)\}\subseteq E_C$.  Since both $E_O$ and $U\cap E_O$ are open and so we can apply Lemma \ref{lem:modifieddv} and Lemma \ref{lem:mu_bullet} to get
\begin{align}
&\liminf_{s\rightarrow\infty}\frac{1}{s}\log \mathbbm{Y}_c(L((Y_x),\rho(s),\cdot)\in U\cap E_C) \nonumber \\ &\quad\geq \liminf_{s\rightarrow\infty}\frac{1}{s}\log \mathbbm{Y}_c(L((Y_x),\rho(s),\cdot)\in U\cap E_O)  \geq - \inf_{\mu\in E_O\cap U}\frac{I_2(\mu)}{\ex(\mu)} = - \frac{I_2(\mu^{\bullet})}{\ex(\mu^{\bullet})}.\label{eq:mubulletU1}
\end{align}
Likewise $E_C\cap U^{\textsc{c}}$ is closed and so from Lemma \ref{lem:modifieddv}  and Lemma \ref{lem:mu_bullet} we have
\begin{align}
\limsup_{s\rightarrow\infty} \frac{1}{s}\log \mathbbm{Y}_c(L((Y_x),\rho(s),\cdot)\in U^{\textsc{c}}\cap E_C) \leq  -\inf_{\mu\in E_C\cap U^{\textsc{c}}}\frac{I_2(\mu)}{\ex(\mu)} < - \frac{I_2(\mu^{\bullet})}{\ex(\mu^{\bullet})}.\label{eq:mubulletU2}
\end{align}
Here there is a strict inequality because the minimiser of $\ex(\mu)^{-1} I_2(\mu)$ over $\mu\in E_C$ is unique.  By combining (\ref{eq:mubulletU1}) and (\ref{eq:mubulletU2}) we can deduce (\ref{eq:mubulletU}) and so the lemma is proved.
\ep
\end{lemma}

Now use $(Z_t)_{t\geq 0}$ to denote the auxiliary process to $(Y_x)_{x\geq 0}$.  From the results of \cite{pinsky1} we know that if we condition on $Z_t\leq 1$ for all $0\leq t\leq s$ then as $s\longrightarrow\infty$ $(Z_t)_{t\geq 0}$ will converge to some stationary process.  Consequentially if we condition on $Y_x\leq 1$ for all $0\leq x\leq \rho(s)$ then $(Y_x)_{x\geq 0}$ must also converge to some stationary process as $s\longrightarrow\infty$.
Furthermore, because we know that the occupation measure of $(Y_x)_{x\geq 0}$ converges to $\mu^{\bullet}$, then $\mu^{\bullet}$ must also be the limiting stationary distribution of $(Y_x)_{x\geq 0}$.

In our proof of ballisticity (Section \ref{sec:ballistic}) we showed that if $(L_x(T))_{x\in\RR}$ is the local time of a Brownian motion, $(W_t)_{t\geq 0}$, conditioned on $\EE_T^{\bullet}$ then $(L_x(T))_{0\leq x\leq W_T}$ is equal in law to $(Y_x)_{0\leq x\leq W_T}$.  Since $W_T\longrightarrow\infty$ as $T\longrightarrow\infty$ it must follow that, with respect to the law $\QQ^{\bullet}$, $(L_x(\infty))_{x\geq 0}$ has an invariant distribution $\mu^{\bullet}$.

\subsection{The tails of $\mu^*$ and $\mu^{\bullet}$}

Theorem \ref{thm:epsilon3a} now follows fairly easily by showing that there exists constants $C^{*}$ and $C^{\bullet}$ with
\begin{align}
\mu^{*}((1-\varepsilon,1])\sim C^{*}\varepsilon^3\quad\text{and}\quad\mu^{\bullet}((1-\varepsilon,1])\sim C^{\bullet}\varepsilon^3 ,
\end{align}
as $\varepsilon\longrightarrow 0$.  This is the content of the lemma below.

\begin{lemma}\label{lem:eps3}
Let $\alpha\in(0,1)$, and suppose $\mu_{\alpha}$ is the unique probability measure supported on $[0,1]$ which has $\ex(\mu_{\alpha})=\alpha$ and $I_2(\mu_{\alpha})=J(\alpha)$, then there exists $C_{\alpha}$ such that
\begin{align}
\mu_{\alpha}((1-\varepsilon,1])\sim C_{\alpha}\varepsilon^3 ,
\end{align}
as $\varepsilon\longrightarrow 0$.
\bp
We prove this by analysing the Radon--Nikodym derivative of $\mu_\alpha$.  The measure $\mu_\alpha$ is defined to be the minimiser of $I_2(\mu)$ over $\{\mu : \mathrm{support}(\mu)\subseteq[0,1]$ and $\ex(\mu)=\alpha\}$.  Therefore we can use the integral form for $I_2$ given by (\ref{eq:Imu_int}), and observe that minimising $I_2(\mu)$ over $\{\mu : \mathrm{support}(\mu)\subseteq[0,1]$ and $\ex(\mu)=\alpha\}$ is equivalent to finding a $g_{\alpha}\in C^1([0,1])$ with $\Vert g_{\alpha} \Vert_2 = 1$, $\Vert \sqrt{x} g_{\alpha}(x) \Vert_2 = \alpha$ and
\begin{align}
\int_{0}^{1} 2x \left(\frac{\mathrm{d}}{\mathrm{d}x}g_{\alpha}(x)\right)^2\dd x = \inf \left\{ \int_{0}^{1} 2x \left(\frac{\mathrm{d}}{\mathrm{d}x}g(x)\right)^2\dd x : \Vert g \Vert_2=1 \text{ and }\Vert \sqrt{x}g(x) \Vert_2 = \alpha \right\}.
\end{align}

As in the proof of Lemma \ref{lem:lemma9} we can do this using the Euler--Lagrange equation.  To include the twin constraints $\Vert g \Vert_2 = 1$ and $\Vert \sqrt{x}g(x) \Vert_2 = \alpha$ we must also include the Lagrangian multipliers
\begin{align}
\lambda\left(\int_{0}^1 g(x)^2\,\dd x - 1 \right) \quad \text{and} \quad \nu\left(\int_{0}^1 x g(x)^2\,\dd x - \alpha \right).
\end{align}
Therefore we get $\mathcal{F}[x,g(x),g^{\prime}(x),\lambda,\nu] = 2x \,g'(x)^2 -2\lambda(g(x)^2-1) -2 \nu(x g(x)^2 - \alpha)$.  Putting this into the Euler--Lagrange equation
\begin{align*}\dfrac{\partial\mathcal{F}}{\partial g}-\dfrac{\dd}{\dd x}\left(\dfrac{\partial\mathcal{F}}{\partial g^{\prime}}\right)=0
\end{align*}
now gives
\begin{align}
x \frac{\mathrm{d}^2}{\mathrm{d}x^2}g(x)+\frac{\mathrm{d}}{\mathrm{d}x}g(x)-(\lambda+\nu x) g(x) = 0, \label{eq:de2}
\end{align}
for $x\in [0,1]$.  The particular $\lambda$ and $\nu$ will  depend on $\alpha$.  However, even without knowing these we can now deduce that $g_{\alpha}$ is twice differentiable on $(0,1)$ -- since it satisfies (\ref{eq:de2}), and that $g_{\alpha}'(1) \neq 0$.  This second claim follows because the Radon-Nikodym derivative of $\mu$ must be continuous on $[0,\infty)$ (else $I_2(\mu)=\infty$) and therefore $g_{\alpha}(1)=0$.  By inspecting (\ref{eq:de2}) we see that if we also had $g_{\alpha}'(1)=0$, then $g_{\alpha}$ would be the trivial solution, $g_{\alpha}(x)=0$.  Since this can not be the case we must therefore have $g'_{\alpha}\neq 0$.

By taking the Taylor expansion of $g_\alpha(x) $ at $x=1$ we now have
\begin{align}
\mu_{\alpha}((1-\varepsilon,1]) &= \int_{1-\varepsilon}^1 \frac{\dd\mu_{\alpha}}{\dd x} \,\dd x = \int_{1-\varepsilon}^1 \left(g_{\alpha}(x)\right)^2\,\dd x = \int_0^{\varepsilon}\left( g_{\alpha}^{\prime}(1)x+\mathcal{O}\left(x^2\right)\right)^2 \, \dd x  
\\ &= \frac{1}{3}g_{\alpha}^{\prime}(1)^2\varepsilon^3 + \mathcal{O}\left(\varepsilon^4\right) \sim C_{\alpha}\varepsilon^3, \nonumber
\end{align}
as required.
\ep
\end{lemma}

Because we know that $\mu^{*}=\mu_{\alpha^*}$ with $\alpha^*=(\gamma^*)^{-1}$ and $\mu^{\bullet}=\mu_{\alpha^{\bullet}}$ with $\alpha^{\bullet}=(\gamma^{\bullet})^{-1}$, then this completes the proof of Theorem \ref{thm:epsilon3a}.
\hfill$\square$

\begin{remark}
For this proof it was not necessary to try and solve (\ref{eq:de2}) explicitly.  However, doing so provides a good way of calculating the function $J$ numerically.
See Definition \ref{def:J} and Figure \ref{fig:J}.
%See Section \ref{sec:I(v)}.
Solving  (\ref{eq:de2}) would also enable us to compute $C^*$ and $C^\bullet$.
\end{remark}

\section{A general framework}\label{sec:discussion}

We conclude this chapter by considering a more general framework in which a Brownian motion $(W_t)_{t\geq 0 }$ can be conditioned to have its local time bounded by 1.  In this section we give three conjectures about the existence and behaviour of these limiting processes, and give non-rigorous explanations of why we believe these conjectures to be true.

%To describe our more general framework
Consider $(W_t,t)_{t\geq 0}$ as a process on $\RR\times[0,\infty)$, and let $U\subseteq\RR\times[0,\infty)$ be an open set containing $(0,0)$.  For each $a\in(0,\infty)$ we write $aU=\{(x,t):(a^{-1}x,a^{-1}t)\in U\}$, and define $\tau_a^U$ by $\tau_a^U=\inf\{t\geq 0:(W_t,t)\notin aU\}$.  From this we then obtain a collection of events
\begin{align}
\EE_a^U=\{L_x(\tau_a^U)\leq 1 \text{ for all } x\in\RR \}.
\end{align}

Observe that since $U$ is open then $\bigcup_{a>0}aU = \RR\times[0,\infty),$ and thus $\bigcap_{a>0}\EE_a^U = \EE\{L_x(t)\leq 1$ for all $x$ and $t\}$.  Note also that for $U=\{(x,t)\in\RR\times[0,\infty):x<1\}$ this definition gives $\EE_a^U=\EE^*_a$, where $\EE^*$ is given by (\ref{eq:Ea*}), and when $U=\{(x,t)\in\RR\times[0,\infty):t<1\text{ and }x\geq 0\}$ then $\EE_a^U=\tilde{\EE}^\bullet_a$, where $\tilde{\EE}^\bullet_a$ is given by (\ref{eq:Eabu0}).
%Therefore considering conditioned measures of the form $\WW(\cdot\,|\,\EE^U_a)$ gives a nice general framework for studying Wiener processes with bounded local time.
We believe that, subject to certain conditions on the set $U$, the measures $\WW(\cdot\,|\,\EE^U_a)$ will weakly converge as $a\longrightarrow\infty$.  Furthermore the behaviour of the limiting process with measure $\mathbbm{Q}^U=\displaystyle{\lim_{a\rightarrow\infty}}\WW(\cdot\,|\,\EE^U_a)$, can be deduced from the set $U$.

\subsection{Understanding $\WW(\cdot\,|\,\EE^U_a)$ via the theory of large deviations}\label{sec:discussion1}

One of the main principles of the theory of large deviations is that when we condition on a process satisfying a sequence of increasingly (exponentially) unlikely events, then -- in the limit -- the path taken by the process is the one which is least unlikely.  Suppose an event $\Theta$ has $\WW(\Theta)=\ee^{-\vartheta}$, then we say the \emph{cost} of $(W_t)_{t\geq 0}$ satisfying $\Theta$ is $\vartheta$.  Now consider a set $U\subseteq\RR\times[0,\infty)$.  Using the language of cost we can then say that if we condition on $\EE_a^U$ then -- in the limit -- the path taken by $a^{-1}(W_t)_{t\geq_0}$ will be the one which is least \emph{expensive}.

It is well know that
\begin{align}
\lim_{T\rightarrow\infty}\frac{1}{T}\log\WW(W_T\geq v T) = -\frac{v^2}{2}, \label{eq:speedcost}
\end{align}
and so it is clearly costly for a Brownian motion to be ballistic with a high speed.  However, in order to satisfy $L_{x}(\tau_a^U)\leq 1$ for all $x\in\RR$, a Brownian motion must travel at a speed of at least 1.  Furthermore, the faster a Brownian motion travels, the easier it is for $(W_t)_{t\geq 0}$ to satisfy the condition on its local time.  Therefore we end up with a trade-off between the cost of travelling quickly and the cost of satisfying $L_{x}(\tau_a^U)\leq 1$.

When a Brownian motion is ballistic then we can estimate its speed though a point $x$ by calculating $\displaystyle{\lim_{T\rightarrow\infty}}\ex(L_x(T))^{-1}$.  Since $L_x(T)$ can be described in terms of a $\mathrm{BES}Q^2$ process, then the cost of a Brownian motion travelling a unit distance at speed $v$, whilst ensuring its local time is bounded by 1, is given by $J(v^{-1})$.  Recall Definition \ref{def:J}.

In the case where $U=aU^* =\{(x,t):x<a\}$,  $(W_t)_{t\geq 0}$ simply has to travel $a$ units of distance whilst maintaining a bounded local time.  Therefore $(W_t)_{t\geq 0}$ will travel at a speed for which the cost $J(v^{-1})$ is minimal.  Lemma \ref{lem:J_props} tells us that $J(v^{-1})$ is minimised when $v=\gamma^*$, and thus we can use the theory of large deviations to confirm the result of Benjamini and Berestycki, \cite{b&b1}.
However, in the  the case where $U=aU^{\bullet} =\{(x,t):t<a\}$ then the limiting process has a speed which is least expensive per unit time.  Because a particle travelling at speed $v$ will cover $v$ units of distance in each unit of time, then we can deduce that the unit time cost of a satisfying $L_x(t)\leq 1$ is $v J(v^{-1})$.  Lemma \ref{lem:J_props} also tells us that $v J(v^{-1})$ has a unique minimum at $\gamma^{\bullet}<\gamma^*$, and so we see why the limiting measures $\QQ^*$ and $\QQ^{\bullet}$ should be different.

\begin{figure}[p]
\centering
\includegraphics[width=.9\textwidth]{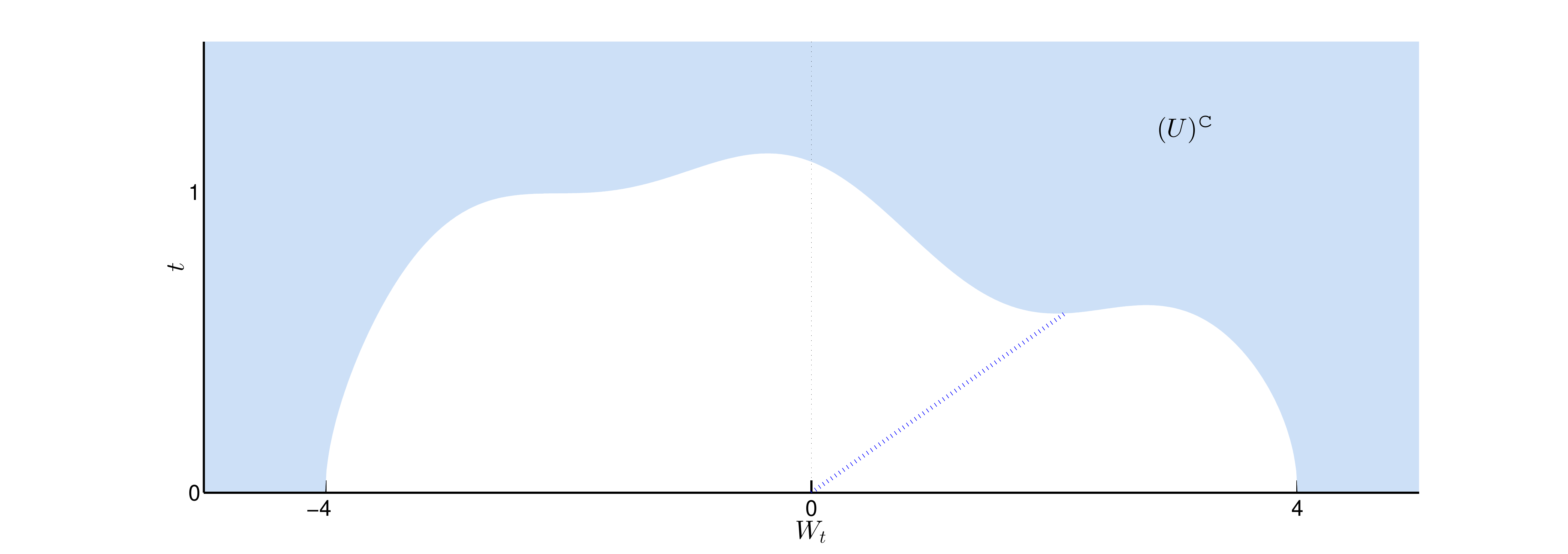}
\includegraphics[width=.9\textwidth]{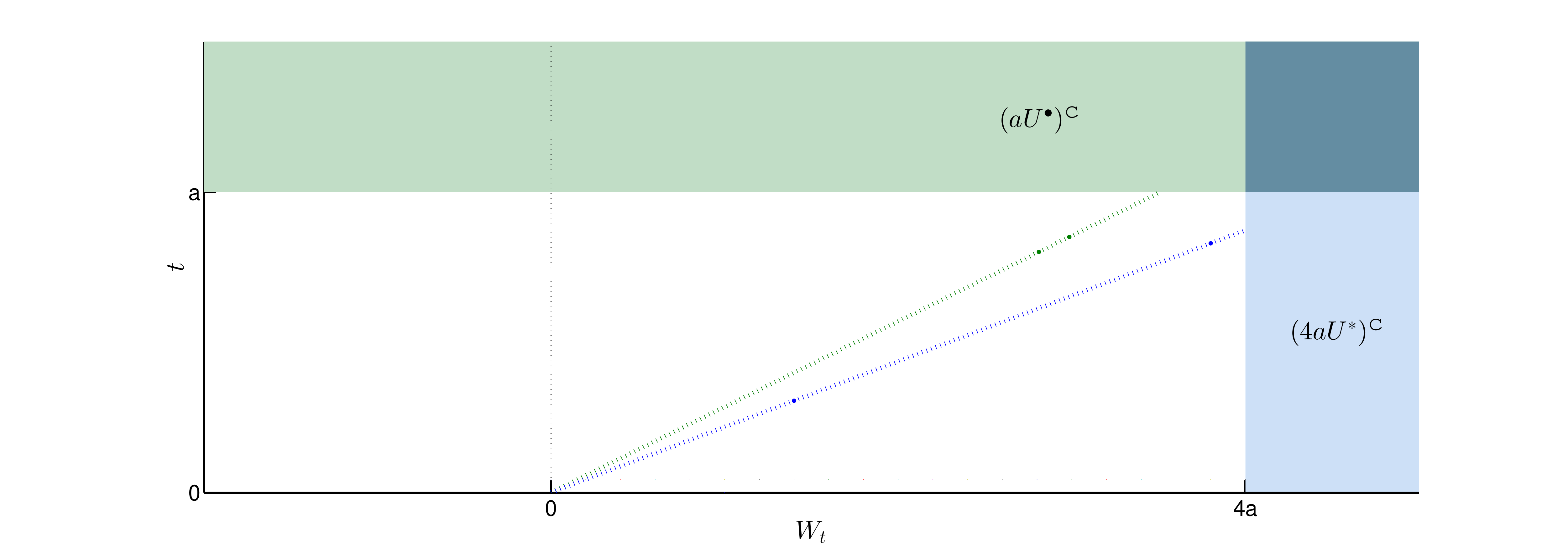}
\includegraphics[width=.9\textwidth]{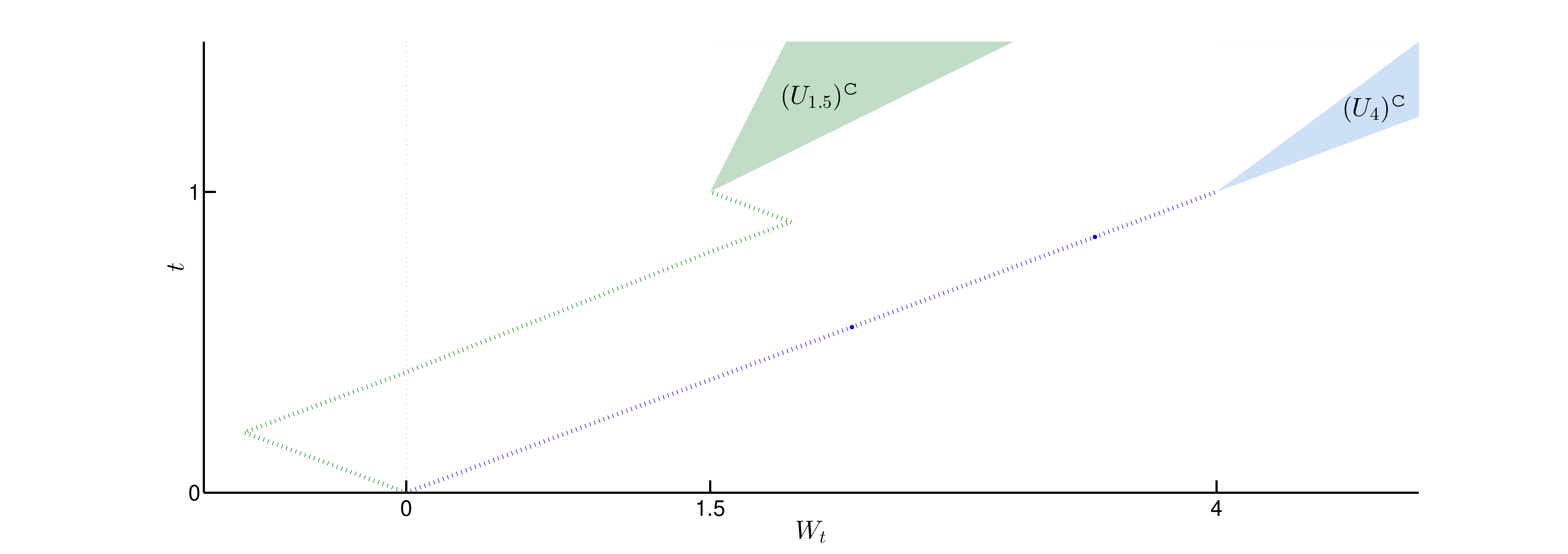}
\caption{The first graph shows the process $(W_t,t)_{t\geq 0}$ taking the path $P$ as it leaves a set $U$.  This need not necessarily be the shortest, but will be the path with lowest cost.  What is more, if the least expensive path $P$ is unique then it must be a straight line.  Obviously the gradient of the line will depend on the set $U$, and this is why different conditionings can lead to Brownian motion with bounded local time having different ballistic rates.
\nl
The second graph shows the most likely paths taken by $(W_t,t)_{t\geq 0}$ as it leaves $aU^{\bullet}$ (green path) and $4aU^{*}$ (blue path). Since $(W_t,t)$ is able to leave $4aU^{*}$ sooner by travelling at a quicker speed we see that $\WW(\cdot\,|\,\EE_{a}^{*})$ and $\WW(\cdot\,|\,\EE_{a}^{\bullet})$ will converge to processes with different ballistic rates.
\nl
The third graph shows that for a given $v$ we can set $U_v$ to be $\RR\times[0,\infty)$ minus a wedge with its point at $(v,1)$.  Provided $v$ is sufficiently large this will give a unique less expensive path for $(W_t,t)_{t\geq 0}$.  However, if $v$ is too small then it ends up being easier for $W_t$ to go back on itself than to go forwards at a slow speed.  Therefore there can be no unique least expensive path, and so the measures $\WW(\cdot\,|\,\EE_{a}^{U_v})$ need not converge.}\label{fig:U}
\end{figure}

Now for a general set $U$ there may be many possible paths $P$ for $a^{-1}(W_t,t)_{t\geq 0}$ to take as it leaves $U$, and for each of these there will be a cost for $a^{-1}(W_t,t)$ to stay close to $P$ and maintain $L_{x}(\tau_a^U) \leq 1$.  However, if we know that  there is is unique least expensive path, then conditionally on $\EE_a^U$ the limiting route taken by $a^{-1}(W_t,t)_{t\geq 0}$ will converge to this least expensive path as $a\longrightarrow\infty$.

\subsection{The cost of {$a^{-1}(W_t,t)_{t\geq 0}$} following a path}

Assume that a path $P$ can be parametrised on an interval $[0,T]$ by $(f(t),t)$ for some piecewise differentiable increasing function $f$.  The cost of $a^{-1}(W_t,t)_{t\geq 0}$ staying close to $P$ will then be asymptotically equal to $a\times\mathrm{cost}(P)$, where
\begin{align}
\mathrm{cost}(P) = \int_0^T J\left(\frac{1}{f'(t)}\right) \,\dd t .
\end{align}
Note that because $J$ is convex then $\mathrm{cost}(P)$ can only be minimised when $P$ is a straight line.

In the case where $P$ is a path whose first parameter is not strictly increasing (or decreasing) then $(W_t)_{t\geq 0}$ must revisit regions where it has already been.  From Lemma \ref{lem:Q0} we can deduce that
\begin{align}
\limsup_{T\rightarrow\infty}\frac{1}{T}\log\WW\left(W_T\in[-1,1] \text{ and }L_x(T)\leq 1 \text{ for all } x\in\RR\right) \leq -2\pi^2 , \label{eq:return_cost1}
\end{align}
and from this we get a lower bound on the cost of $a^{-1}(W_t,t)_{t\geq 0}$ following a path which comes back on itself.  Using (\ref{eq:return_cost1}) as lower bound it is possible to show that for certain $U$ there is a unique path $P$ from $(0,0)$ to $\partial U$ which minimises $\mathrm{cost}(P)$.  We now make the following conjecture.

\begin{conjecture}\label{con:limit}
Suppose $(0,0)\in U\subseteq\RR\times[0,\infty)$ is an open set, and assume that there is a unique path $P$ from $(0,0)$ to $\partial U$ with $\mathrm{cost}(P)< \mathrm{cost}(\tilde{P})$ for all other paths $\tilde{P}$ from $(0,0)$ to $\partial U$.  Then
\vspace{-0.5\baselineskip}\begin{itemize}
\item $P$ can be parametrised by $(vt,t)$ for some constant $v$.
\item $\WW(\cdot\,|\,\EE^U_a)$ converges to a weak limit, $\QQ^U$, as $a\longrightarrow\infty$.
\item $\QQ^U$ is such that $\dfrac{W_t}{t}\longrightarrow v$ in $\QQ^U$-probability.
\end{itemize}
\end{conjecture}

%\subsection{The lowest possible speed}

Following on from Conjecture \ref{con:limit} we now ask for which values of $v$ can we find a set $U_v$ such that $\QQ^{U_v}$ exists and has ballistic rate $v$.  For simplicity we shall restrict our attention to the case where $v\geq 0$. 

Clearly it is not possible for us to have $v<1$ as this would imply that $L_x(T)> 1$ for some $x$ and $t$.  What is more, Lemma \ref{lem:J_props} tells us that $J\left(v^{-1}\right)\longrightarrow\infty$ as $v\searrow 1$ and so we see that it is very expensive for $(W_t)_{t\geq 0}$ to maintain $L_x(t)\leq 1$ whilst travelling slowly.  Consequentially, if we want $(W_t)_{t\geq 0}$ to be near $vT$ at time $T$ and satisfy $L_x(T)\leq 1$ for all $x\in\RR$, then when $v$ is small the least expensive way for this to happen is if $(W_t)_{t\geq 0}$ goes off at some speed greater than $v$ and then changes direction in order to come back to $vT$.

From (\ref{eq:return_cost1}) we can see that the asymptotic cost of $(W_t)_{t\geq 0}$ returning to the interval $[0,vT]$ after $\lambda T$ units of time is at least $2\pi^2\lambda T$.  In fact without too much difficulty one can show that this lower bound is sharp.

We also know that the cost of $(W_t)_{0\leq t\leq T}$ spending $(1-\lambda)T$ units of time in the interval $[0,vT]$ whilst maintaining $L_x(T)\leq 1$ is asymptotically equal to $vT J\left((1-\lambda)v^{-1}\right)$.  Therefore we see that staying close to $(vt,t)_{0\leq t\leq T}$ proves to be the least expensive way for $(W_t,t)_{t\geq 0}$ to end up near $(vT,T)$ if and only if
\begin{align}
v J\left(v^{-1}\right) < \lambda 2\pi^2 + (1-\lambda)v J\left( (1-\lambda) v^{-1}\right) \label{eq:v_min}
\end{align}
for all $0<\lambda\leq 1$.  Although we do not include details, it can be shown by studying the properties of $J$ that there is a critical value $1<\gamma^{\circ}<\gamma^{\bullet}$, equal to the minimal root of $v J\left(v^{-1}\right)=2\pi^2$, such that (\ref{eq:v_min}) is satisfied for all $v>\gamma^{\circ}$ and for no $v<\gamma^{\circ}$.  See Figure \ref{fig:J} and Figure \ref{fig:U}.

\begin{conjecture}\label{con:rate}
There exists $1<\gamma^{\circ}<\gamma^{\bullet}$ such that for each $v>\gamma^{\circ}$ there is an open set $U_v$ for which $\QQ^{U_v}=\displaystyle{\lim_{a\rightarrow\infty} \WW(\cdot\,|\,\EE^{U_v}_a)}$ exists and is such that
\begin{align}
\lim_{t\rightarrow\infty} \frac{W_t}{t} = v \quad\text{in }\QQ^{U_v}\text{-probability} \label{eq:ballistic_}.
\end{align}
What is more, for each $v<\gamma^{\circ}$ there is no open set $U_v$ for which (\ref{eq:ballistic_}) holds. 
\end{conjecture}

%In proving Theorem \ref{thm:reduced_speed} and Theorem \ref{thm:epsilon3a} we introduce many (although not all) of the tools that would be required to prove the two conjectures above.  However, for the sake of expediency and space we content ourselves to give only the preceding discussion of why they should in fact be true.

\subsection{The universal exponent}\label{sec:rem_gen}

Suppose the measure $\displaystyle{\QQ^v=\lim_{a\rightarrow\infty}\WW(\cdot\,|\,\EE_a^{U_v})}$ has $\displaystyle{\lim_{t\rightarrow\infty} \frac{W_t}{t} = v}$ in $\QQ^{v}$-probability, for some $v>\gamma^{\circ}$.  Provided we knew that $(L_x(\tau_{a}^{U_v}))_{x\geq 0}$ converged to a stationary process as $a\longrightarrow\infty$, then because the occupation measure of $(L_x(\infty))_{x\geq 0}$ must converge to a measure which minimises $I_2(\mu)$ over all $\mu\in\{\mu\in\mathcal{P}(\RR):\mathrm{support}(\mu)\subseteq [0,1]$ and $\ex(\mu)=v^{-1}\}$, we would be able to deduce that the stationary measure of $(L_x(\infty))_{x\geq 0}$ with respect to $\QQ^v$ is $\mu_{v^{-1}}$.   Lemma \ref{lem:eps3}  then tells us that there is a constant $C_v>0$ for which we have
\begin{align}
\mu_{v^{-1}}((1-\varepsilon,1])\sim C_v \varepsilon^3,
\end{align}
and so we can make the following conjecture.

\begin{conjecture}\label{con:epsilon}
For each $v>\gamma^{\circ}$ and each measure $\displaystyle{\QQ^{U_v} = \lim_{a\rightarrow\infty}\WW(\cdot\,|\,\EE_{a}^{U_v})}$ there exists a constant $C_v>0$ with
\begin{align}
\lim_{x\rightarrow\infty}\QQ^{v}(L_x(\infty)>1-\varepsilon)\sim C_{v} \varepsilon^3 .
\end{align}
\end{conjecture}

\section*{Acknowledgements}

I would like to give thanks to Nathana\"el Berestycki for introducing me to this problem and for helpful and stimulating discussion.
This work has been supported by the UK Engineering and Physical Sciences Research Council (EPSRC) grant EP/H023348/1.

\bibliographystyle{plain}

\end{document}